\documentclass[11pt,leqno]{amsart}
\usepackage{amsmath, amssymb}
\usepackage{graphicx,color,hyperref}
\usepackage{verbatim,enumitem}

\def\smallskip{\vskip3pt plus1pt minus1pt}
\def\medskip{\vskip5pt plus2pt minus1pt}

\def\plainsubsection#1|{%
  \par\vskip0.25cm\penalty -100
  \centerline{{\sc #1}}
  \vskip3pt plus 1pt minus 0pt
  \penalty 500}

\long\def\claim#1|#2\endclaim{\par\vskip 4pt\noindent 
{\bf #1.}\ {\em #2}\par\vskip 5pt}

\def\plainproof{\noindent{\em Proof}}

\def\today{\ifcase\month\or
January\or February\or March\or April\or May\or June\or July\or August\or
September\or October\or November\or December\fi \space\number\day,
\number\year}

\catcode`\@=11
\newcount\@tempcnta \newcount\@tempcntb 
\def\timeofday{{%
\@tempcnta=\time \divide\@tempcnta by 60 \@tempcntb=\@tempcnta
\multiply\@tempcntb by -60 \advance\@tempcntb by \time
\ifnum\@tempcntb > 9 \number\@tempcnta:\number\@tempcntb
\else\number\@tempcnta:0\number\@tempcntb\fi}}

\def\openup{\afterassignment\@penup\dimen@=}
\def\@penup{\advance\lineskip\dimen@
  \advance\baselineskip\dimen@
  \advance\lineskiplimit\dimen@}
\newdimen\jot \jot=3pt
\newskip\plaincentering \plaincentering=0pt plus 1000pt minus 1000pt
\def\ialign{\everycr{}\tabskip\z@skip\halign}
\def\eqalign#1{\null\,\vcenter{\openup\jot\m@th
  \ialign{\strut\hfil$\displaystyle{##}$&$\displaystyle{{}##}$\hfil
      \crcr#1\crcr}}\,}
\newif\ifdt@p
\def\displ@y{\global\dt@ptrue\openup\jot\m@th
  \everycr{\noalign{\ifdt@p \global\dt@pfalse \ifdim\prevdepth>-1000\p@
      \vskip-\lineskiplimit \vskip\normallineskiplimit \fi
      \else \penalty\interdisplaylinepenalty \fi}}}
\def\@lign{\tabskip\z@skip\everycr{}} 
\def\displaylines#1{\displ@y \tabskip\z@skip
  \halign{\hbox to\displaywidth{$\@lign\hfil\displaystyle##\hfil$}\crcr
    #1\crcr}}
\def\eqalignno#1{\displ@y \tabskip\plaincentering
  \halign to\displaywidth{\hfil$\@lign\displaystyle{##}$\tabskip\z@skip
    &$\@lign\displaystyle{{}##}$\hfil\tabskip\plaincentering
    &\llap{$\@lign##$}\tabskip\z@skip\crcr
    #1\crcr}}
\def\leqalignno#1{\displ@y \tabskip\plaincentering
  \halign to\displaywidth{\hfil$\@lign\displaystyle{##}$\tabskip\z@skip
    &$\@lign\displaystyle{{}##}$\hfil\tabskip\plaincentering
    &\kern-\displaywidth\rlap{$\@lign##$}\tabskip\displaywidth\crcr
    #1\crcr}}
\def\plaincases#1{\left\{\,\vcenter{\normalbaselines\m@th
    \ialign{$##\hfil$&\quad##\hfil\crcr#1\crcr}}\right.}
\def\plainmatrix#1{\null\,\vcenter{\normalbaselines\m@th
    \ialign{\hfil$##$\hfil&&\quad\hfil$##$\hfil\crcr
      \mathstrut\crcr\noalign{\kern-\baselineskip}
      #1\crcr\mathstrut\crcr\noalign{\kern-\baselineskip}}}\,}

\catcode`\@=12

\def\Bibitem#1&#2&#3&#4&%
{\hangindent=2cm\hangafter=1
\noindent\rlap{\hbox{\bf #1}}\kern2cm{\rm #2}{\it #3}{\rm #4.}} 

\let\em=\it

\font\twelvebf=cmbx10 at 12pt
\font\tenbf=cmbx10
\font\eightbf=cmbx10 at 8pt
\font\sevenbf=cmbx10 at 7pt
\font\sixbf=cmbx10 at 6pt
\font\fivebf=cmbx10 at 5pt

\font\eighti=cmmi8
\font\sixi=cmmi6

\font\eightrm=cmr10 at 8pt
\font\sixrm=cmr10 at 6pt

\font\eightsy=cmsy8
\font\sixsy=cmsy6

\font\eightit=cmti10 at 8pt
\font\eightsl=cmsl10 at 8pt
\font\eighttt=cmtt10 at 8pt

\font\twelvebsy=cmbsy10 at 12pt
\font\tenbsy=cmbsy10
\font\eightbsy=cmbsy8
\font\sevenbsy=cmbsy7
\font\sixbsy=cmbsy6
\font\fivebsy=cmbsy5

\font\tenmsa=msam10

\font\sevenmsa=msam7
\font\fivemsa=msam5
\newfam\msafam
  \textfont\msafam=\tenmsa
  \scriptfont\msafam=\sevenmsa
  \scriptscriptfont\msafam=\fivemsa

\font\tenmsb=msbm10
\font\eightmsb=msbm8
\font\sevenmsb=msbm7
\font\fivemsb=msbm5
\newfam\msbfam
  \textfont\msbfam=\tenmsb
  \scriptfont\msbfam=\sevenmsb
  \scriptscriptfont\msbfam=\fivemsb
\def\Bbb{\fam\msbfam\tenmsb}

\font\tenCal=eusm10
\font\sevenCal=eusm7
\font\fiveCal=eusm5
\newfam\Calfam
  \textfont\Calfam=\tenCal
  \scriptfont\Calfam=\sevenCal
  \scriptscriptfont\Calfam=\fiveCal
\def\Cal{\fam\Calfam\tenCal}

\font\teneuf=eusm10
\font\teneuf=eufm10
\font\seveneuf=eufm7
\font\fiveeuf=eufm5
\newfam\euffam
  \textfont\euffam=\teneuf
  \scriptfont\euffam=\seveneuf
  \scriptscriptfont\euffam=\fiveeuf

\font\twelvebfit=cmmib10 at 12pt
\font\tenbfit=cmmib10
\font\eightbfit=cmmib8
\font\sevenbfit=cmmib7
\font\sixbfit=cmmib6
\font\fivebfit=cmmib5
\newfam\bfitfam
  \textfont\bfitfam=\tenbfit
  \scriptfont\bfitfam=\sevenbfit
  \scriptscriptfont\bfitfam=\fivebfit

\catcode`\@=11
\def\eightpoint{%
  \textfont0=\eightrm \scriptfont0=\sixrm \scriptscriptfont0=\fiverm
  \def\rm{\fam\z@\eightrm}%
  \textfont1=\eighti \scriptfont1=\sixi \scriptscriptfont1=\fivei
  \def\oldstyle{\fam\@ne\eighti}%
  \textfont2=\eightsy \scriptfont2=\sixsy \scriptscriptfont2=\fivesy
  \textfont\itfam=\eightit
  \def\it{\fam\itfam\eightit}%
  \textfont\slfam=\eightsl
  \def\sl{\fam\slfam\eightsl}%
  \textfont\bffam=\eightbf \scriptfont\bffam=\sixbf
  \scriptscriptfont\bffam=\fivebf
  \def\bf{\fam\bffam\eightbf}%
  \textfont\ttfam=\eighttt
  \def\tt{\fam\ttfam\eighttt}%
  \textfont\msbfam=\eightmsb
  \def\Bbb{\fam\msbfam\eightmsb}%
  \abovedisplayskip=9pt plus 2pt minus 6pt
  \abovedisplayshortskip=0pt plus 2pt
  \belowdisplayskip=9pt plus 2pt minus 6pt
  \belowdisplayshortskip=5pt plus 2pt minus 3pt
  \smallskipamount=2pt plus 1pt minus 1pt
  \medskipamount=4pt plus 2pt minus 1pt
  \bigskipamount=9pt plus 3pt minus 3pt
  \normalbaselineskip=9.6pt
  \setbox\strutbox=\hbox{\vrule height7pt depth2pt width0pt}%
  \let\bigf@ntpc=\eightrm \let\smallf@ntpc=\sixrm
  \normalbaselines\rm}
\catcode`\@=12

\def\eightpointbf{%
 \textfont0=\eightbf   \scriptfont0=\sixbf   \scriptscriptfont0=\fivebf
 \textfont1=\eightbfit \scriptfont1=\sixbfit \scriptscriptfont1=\fivebfit
 \textfont2=\eightbsy  \scriptfont2=\sixbsy  \scriptscriptfont2=\fivebsy
 \eightbf
 \baselineskip=10pt}

\def\tenpointbf{%
 \textfont0=\tenbf   \scriptfont0=\sevenbf   \scriptscriptfont0=\fivebf
 \textfont1=\tenbfit \scriptfont1=\sevenbfit \scriptscriptfont1=\fivebfit
 \textfont2=\tenbsy  \scriptfont2=\sevenbsy  \scriptscriptfont2=\fivebsy
 \tenbf}
        
\def\twelvepointbf{%
 \textfont0=\twelvebf   \scriptfont0=\eightbf   \scriptscriptfont0=\sixbf
 \textfont1=\twelvebfit \scriptfont1=\eightbfit \scriptscriptfont1=\sixbfit
 \textfont2=\twelvebsy  \scriptfont2=\eightbsy  \scriptscriptfont2=\sixbsy
 \twelvebf
 \baselineskip=14.4pt}

\newdimen\srdim \srdim=\hsize
\newdimen\irdim \irdim=\hsize
\def\NOSECTREF#1{\noindent\hbox to \srdim{\null\dotfill ???(#1)}}
\def\SECTREF#1{\noindent\hbox to \srdim{\csname REF\romannumeral#1\endcsname}}
\def\INDREF#1{\noindent\hbox to \irdim{\csname IND\romannumeral#1\endcsname}}
\newlinechar=`\^^J
  
\newbox\plaintitlebox   \setbox\plaintitlebox\hbox{\hfil}
\newbox\plainsectionbox \setbox\plainsectionbox\hbox{\hfil}
\def\folio{\ifnum\pageno=1 \hfil \else \ifodd\pageno
           \hfil {\eightpoint\copy\plainsectionbox\kern8mm\number\pageno}\else
           {\eightpoint\number\pageno\kern8mm\copy\plaintitlebox}\hfil \fi\fi}

\def\plaintitlerunning#1{\setbox\plaintitlebox\hbox{\eightpoint #1}}

\def\plainsection#1{%
  \par\vskip0.5cm\penalty -100
  \vbox{\baselineskip=14.4pt\noindent{{\twelvepointbf #1}}}
  \vskip2pt
  \penalty 500}

\def\plainsubsection#1|{%
  \par\vskip0.25cm\penalty -100
  \vbox{\noindent{{\tenpointbf #1}}}
  \vskip1pt
  \penalty 500}

\newdimen\plainitemindent \plainitemindent=6mm

\def\plainitem#1{\par\noindent\hangindent\plainitemindent%
            \rlap{#1}\kern\plainitemindent\ignorespaces}
\def\plainitemitem#1{\par\noindent\hangindent2\plainitemindent%
            \kern\plainitemindent\rlap{#1}\kern\plainitemindent\ignorespaces}
\def\plainitemitemitem#1{\par\noindent\hangindent3\plainitemindent%
            \kern2\plainitemindent\rlap{#1}\kern\plainitemindent\ignorespaces}

\long\def\claim#1|#2\endclaim{\par\vskip 5pt\noindent 
{\tenpointbf #1.}\ {\em #2}\par\vskip 5pt}

\def\plainproof{\noindent{\em Proof}}

\def\today{\ifcase\month\or
January\or February\or March\or April\or May\or June\or July\or August\or
September\or October\or November\or December\fi \space\number\day,
\number\year}

\catcode`\@=11
\newcount\@tempcnta \newcount\@tempcntb 
\def\timeofday{{%
\@tempcnta=\time \divide\@tempcnta by 60 \@tempcntb=\@tempcnta
\multiply\@tempcntb by -60 \advance\@tempcntb by \time
\ifnum\@tempcntb > 9 \number\@tempcnta:\number\@tempcntb
  \else\number\@tempcnta:0\number\@tempcntb\fi}}
\catcode`\@=12

\def\bibitem#1&#2&#3&#4&%
{\hangindent=1.85cm\hangafter=1
\noindent\rlap{\hbox{\eightpointbf #1}}\kern1.85cm{\rm #2}{\it #3}{\rm #4.}} 

\def\bB{{\Bbb B}}
\def\bC{{\Bbb C}}
\def\bD{{\Bbb D}}
\def\bE{{\Bbb E}}
\def\bG{{\Bbb G}}

\def\bN{{\Bbb N}}
\def\bP{{\Bbb P}}
\def\bQ{{\Bbb Q}}
\def\bR{{\Bbb R}}
\def\bS{{\Bbb S}}

\def\bZ{{\Bbb Z}}
\def\bOne{{\mathchoice {\rm 1\mskip-4mu l} {\rm 1\mskip-4mu l}
    {\rm 1\mskip-4.5mu l} {\rm 1\mskip-5mu l}}}

\font\sixbf=cmbx10 at 6pt
\font\fourbf=cmbx10 at 4pt
\font\threebf=cmbx10 at 3pt
\def\bDelta{{\mathchoice
    {\Delta\rlap{\kern-6.5pt\raise1.8pt\hbox{\sixbf/}}}
    {\Delta\rlap{\kern-6.5pt\raise1.8pt\hbox{\sixbf/}}}    
    {\Delta\rlap{\kern-4.9pt\raise1pt\hbox{\fourbf/}}}
    {\Delta\rlap{\kern-4pt\raise1pt\hbox{\threebf/}}}}}    
  
\def\cE{{\Cal E}}

\def\cI{{\Cal I}}
\def\cJ{{\Cal J}}

\def\cL{{\Cal L}}

\def\cO{{\Cal O}}

\def\cX{{\Cal X}}
\def\cY{{\Cal Y}}

\def\rS{{\kern0.18ex{\rm S}}}

\def\\{\hfil\break}
\def\ii{\mathop{\rm i}}

\def\ldotsp{\mathrel{.\kern1pt.\kern1pt.}}
\def\bu{{\scriptstyle\bullet}}
\def\ort{\mathop{\hbox{\kern1pt\vrule width4pt height0.4pt depth0pt
      \vrule width0.4pt height7pt depth0pt\kern3pt}}}

\def\square{\hbox{%
\vrule height 1.5ex  width 0.1ex  depth 0ex\kern-0.1ex
\vrule height 1.5ex  width 1.5ex  depth -1.4ex\kern-1.5ex
\vrule height 0.1ex  width 1.5ex  depth 0ex\kern-0.1ex
\vrule height 1.5ex  width 0.1ex  depth 0ex}}
\def\qed{\strut~\hfill\square\vskip6pt plus2pt minus1pt}

\def\hexnbr#1{\ifnum#1<10 \number#1\else
 \ifnum#1=10 A\else\ifnum#1=11 B\else\ifnum#1=12 C\else
 \ifnum#1=13 D\else\ifnum#1=14 E\else\ifnum#1=15 F\fi\fi\fi\fi\fi\fi\fi}
\def\msatype{\hexnbr\msafam}
\def\msbtype{\hexnbr\msbfam}
\mathchardef\restriction="3\msatype16   
\mathchardef\compact="3\msatype62
\mathchardef\complement="0\msatype7B
\mathchardef\smallsetminus="2\msbtype72   \let\ssm\smallsetminus
\mathchardef\subsetneq="3\msbtype28
\mathchardef\supsetneq="3\msbtype29
\mathchardef\leqslant="3\msatype36   \let\le\leqslant
\mathchardef\geqslant="3\msatype3E   \let\ge\geqslant
\mathchardef\ltimes="2\msbtype6E
\mathchardef\rtimes="2\msbtype6F

\let\ol=\overline

\let\wt=\widetilde
\let\wh=\widehat
\def\swt#1|{\smash{\widetilde#1}}
\def\swh#1|{\smash{\widehat#1}}
\def\build#1|#2|#3|{\mathrel{\mathop{\null#1}\limits^{#2}_{#3}}}
\def\buildo#1^#2{\mathrel{\mathop{\null#1}\limits^{#2}}}
\def\buildu#1_#2{\mathrel{\mathop{\null#1}\limits_{#2}}}

\mathchardef\rsa"3\msatype20
\def\vlra#1|{\hbox to#1mm{\rightarrowfill}}
\def\vlhra#1|{\lhook\joinrel\hbox to#1mm{\rightarrowfill}}

\def\lraww{\mathrel{\rlap{$\longrightarrow$}\kern-1pt\longrightarrow}}
\def\hdashpiece{\hbox{\vrule height2.45pt depth-2.15pt width2.3pt\kern1.2pt}}
\def\dashto{\mathrel{\hdashpiece\hdashpiece\kern-0.5pt\hbox{\tenmsa K}}}
\def\dasharrow{\mathrel{\hdashpiece\hdashpiece\hdashpiece
    \kern-0.3pt\hbox{\tenmsa K}}}
\def\vdashpiece{\smash{\hbox{\vrule height2.3pt depth0pt width0.3pt}}}
\def\vdashto{\mathrel{\smash{\hbox{\vbox{
  \bgroup\baselineskip=3.5pt
  \vdashpiece\vdashpiece\smash{\raise2pt\hbox{\kern-3.1pt\tenex y}}
  \egroup}}}}}
\def\vdasharrow{\mathrel{\smash{\hbox{\vbox{
  \bgroup\baselineskip=3.5pt
  \vdashpiece\vdashpiece\vdashpiece
  \smash{\raise2pt\hbox{\kern-3.1pt\tenex y}}
  \egroup}}}}}
\catcode`\@=11
\newdimen\@rrowlength \@rrowlength=6ex
\def\ssrelbar{\vrule width\@rrowlength height0.64ex depth-0.56ex\kern-4pt}
\def\llra#1{\@rrowlength=#1\ssrelbar\rightarrow}
\catcode`\@=12

\def\lcm{\mathop{\rm lcm}\nolimits}

\def\Re{\mathop{\rm Re}\nolimits}
\def\Im{\mathop{\rm Im}\nolimits}
\def\Arg{\mathop{\rm Arg}\nolimits}

\def\Hom{\mathop{\rm Hom}\nolimits}

\def\SM{\mathop{\rm SM}\nolimits}
\def\Herm{\mathop{\rm Herm}\nolimits}

\def\Tr{\mathop{\rm Tr}\nolimits}

\def\Pic{\mathop{\rm Pic}\nolimits}

\def\Proj{\mathop{\rm Proj}\nolimits}

\def\Vol{\mathop{\rm Vol}\nolimits}

\def\Span{\mathop{\rm Span}\nolimits}

\def\rank{\mathop{\rm rank}\nolimits}
\def\div{\mathop{\rm div}\nolimits}

\def\dbar{{\overline\partial}}
\def\ddbar{{\partial\overline\partial}}

\def\Id{{\rm Id}}
\def\Sing{{\rm Sing}}
\def\orb{{\rm orb}}
\def\FS{{\rm FS}}
\def\GG{{\rm GG}}
\def\loc{{\rm loc}}

\title[On the existence of logarithmic and orbifold jet differentials]
{On the existence of logarithmic\vskip3pt and orbifold jet differentials}

\author[Fr\'ed\'eric Campana, Lionel Darondeau, Jean-Pierre Demailly,
Erwan Rousseau]{Fr\'ed\'eric Campana, Lionel Darondeau,\vskip3pt
Jean-Pierre Demailly,  Erwan Rousseau}

\begin{document}

\begin{abstract} We introduce the concept of directed orbifold, namely
triples $(X,V,D)$ formed by a directed algebraic or analytic variety
$(X,V)$, and a ramification divisor $D$, where $V$ is a coherent subsheaf of
the tangent bundle~$T_X$. In~this context, we introduce an algebra of
orbifold jet differentials and their sections. These jet sections can
be seen as algebraic differential
operators acting on germs of curves, with meromorphic coefficients,
whose poles are supported by $D$ and multiplicities are bounded by
the ramification indices of the components of $D$. We estimate
precisely the curvature tensor of the corresponding directed
structure $V\langle D\rangle$ in the general orbifold case -- with
a special attention to the compact case $D=0$ and to the
logarithmic situation where the ramification indices are infinite.
Using holomorphic Morse inequalities on the
tautological line bundle of the projectivized orbifold
Green-Griffiths bundle, we finally obtain effective sufficient
conditions for the existence of global orbifold jet differentials.

\vskip3pt\noindent
{\sc Keywords.} Projective variety, directed variety, orbifold,
ramification divisor, entire curve, jet differential, Green-Griffiths
conjecture, algebraic differential operator, holomorphic Morse
inequalities, Chern curvature, Chern form.

\vskip3pt\noindent
{\sc MSC classification 2020.} 32Q45, 32H30, 14F06

\vskip3pt\noindent
{\sc Funding.} The third author is supported by the Advanced ERC grant
ALKAGE, no 670846 from September 2015, attributed by the
European Research Council.
\end{abstract}

\maketitle

\plainsection{0. Introduction and main definitions}

The present work is concerned primarily with the existence of logarithmic and
orbifold jet differentials on projective varieties. For the sake of generality,
and in view of potential applications to the case of foliations, we work
throughout this paper in the category of directed varieties, and generalize
them by introducing the concept of directed orbifold.

\claim 0.1. Definition|Let $X$ be a complex manifold or variety.
A directed structure $(X,V)$ on $X$ is defined to be a subsheaf
$V\subset\cO(T_X)$ such that
$\cO(T_X)/V$ is torsion free. A~morphism of directed varieties
$\Psi:(X,V)\to (Y,W)$ is a holomorphic map $\Psi:X\to Y$ such that
$d\Psi(V)\subset\Psi^*W$. We say that $(X,V)$ is non singular if
$X$ is non singular and $V$ is locally free, i.e., is a holomorphic
subbundle of $T_X$.
\endclaim

We refer to the {\it absolute case} as being the situation when $V=T_X$,
the {\it relative case} when $V=T_{X/S}$ for some fibration $X\to S$, and
the {\it foliated case} when $V$ is integrable, i.e.\ $[V,V]\subset V$,
that is, $V$ is the tangent sheaf to a holomorphic foliation. We now combine
these concepts with orbifold structures in the sense
of Campana [Cam04].

\claim 0.2. Definition|A directed orbifold is a triple $(X,V,D)$
where $(X,V)$ is a directed variety and $D=\sum(1-{1\over\rho_j})\Delta_j$
an effective real divisor, where $\Delta_j$ is an irreducible
hypersurface and $\rho_j\in{}]1,\infty]$ an associated ``ramification number''.
We denote by $\lceil D\rceil=\sum\Delta_j$ the
corresponding reduced divisor, and by $|D|=\bigcup \Delta_j$ its support.
\vskip2pt
\plainitem{\rm(a)} We will say that $(X,V,D)$ is non singular if
$(X,V)$ is non singular and
$D$ is a simple normal crossing divisor such that $D$ is transverse
to $V$. If $r=\rank(V)$, we mean by this that there are at most $r$ components
$\Delta_j$ meeting at any point $x\in X$, and that for any $p$-tuple
$(j_1,\ldots,j_p)$ of indices, $1\le p\le r$, we have
$\dim V_x\cap \bigcap_{j=1}^pT_{\Delta_{j_\ell},x}=r-p$ at
any point $x\in\bigcap_{j=1}^p\Delta_{j_\ell}$.
\vskip2pt
\plainitem{\rm(b)} If $(X,V,D)$ is non singular, the canonical divisor of
$(X,V,D)$ is defined to be
$$
K_{V,D}=K_V+ D
$$
$($in additive notation$)$, where $K_V=\det V^*$.
\vskip2pt
\plainitem{\rm(c)} The so called logarithmic case corresponds to all multiplicities
$\rho_j=\infty$ being taken infi\-nite, so that $D=
\sum\Delta_j=\lceil D\rceil$.\vskip0pt
\endclaim

In case $V=T_X$, we recover the concept of orbifold introduced in
[Cam04], except possibly for the fact that we allow here
$\rho_j>1$ to be real or $\infty$, (even though the case where
$\rho_j$ is in $\bN\cup\{\infty\}$ is of greater interest).
It would certainly be interesting to investigate the case
when $(X,V,D)$ is singular, by allowing singularities
in $V$ and tangencies between $V$ and~$D$, and to study whether the
results discussed in this paper can be extended in some way, e.g.\ by
introducing suitable multiplier ideal sheaves taking care of singularities,
as was done in [Dem15] for the study of directed varieties $(X,V)$.
For the sake of technical simplicity, we will refrain to do so here, and will
therefore leave for future work the study of singular directed orbifolds.

\claim 0.3. Definition|Let $(X,V,D)$ be a singular directed orbifold.
We say that $f:\bC\to X$ is an orbifold entire curve if $f$ is a non
constant holomorphic map such that$\;:$\vskip2pt
\plainitem{\rm(a)} $f$ is tangent to $V$
$($i.e.\ $f'(t)\in V_{f(t)}$ at every point,
or equivalently $f:(\bC,T_\bC)\to (X,V)$ is a morphism of directed
varieties$\,;$
\vskip2pt
\plainitem{\rm(b)} $f(\bC)$ is not identically contained in $|D|\,;$
\vskip2pt
\plainitem{\rm(c)} at every point $t_0\in\bC$ such that $f(t_0)\in\Delta_j$,
  $f$ meets $\Delta_j$ with ramification number${}\ge\rho_j$, i.e., if
$\Delta_j=\{z_j=0\}$ near $f(t_0)$, then $z_j\circ f(t)$ vanishes with
multiplicity${}\ge\rho_j$ at $t_0$.
\vskip2pt\noindent
In the case of a logarithmic component $\Delta_j$
$(\rho_j=\infty)$, condition {\rm(c)} is to be replaced by the assumption
\vskip2pt
\plainitem{$({\rm c}')$} $f(\bC)$ does not meet~$\Delta_j$.\vskip0pt
\endclaim

\noindent
One can now consider a category of directed orbifolds as follows.

\claim 0.4. Definition|Consider directed orbifolds
$(X,V,D)$, $(Y,W,D')$ with
$$
D=\sum\Big(1-{1\over\rho_i}\Big)\Delta_i,\qquad
D'=\sum\Big(1-{1\over\rho'_j}\Big)\Delta'_j.
$$
A morphism $\Psi:(X,V,D)\to(Y,W,D')$ is a morphism
$\Psi:(X,V)\to(Y,W)$ of directed varieties satisfying the
additional following  properties {\rm(a,b,c)}.
\vskip2pt
\plainitem{\rm(a)} for every component $\Delta'_j$, $\Psi^{-1}(\Delta'_j)$
consists of a union of components $\Delta_i$, $i\in I(j)$,
eventually after adding a number of extra components $\Delta_i$
with $\rho_i=1\;;$
\vskip2pt
\plainitem{\rm(b)} in case $\rho'_j<\infty$, for every $i\in I(j)$
and $z\in\Delta_i$, 
the derivatives $d^\alpha\Psi(z)$ of $\Psi$ at $z$, computed
in suitable local coordinates on $X$ and $Y$, vanish for all multi-indices
$\alpha\in\bN^n$ with $0<|\alpha|<\rho'_j/\rho_i\;;$
\vskip2pt
\plainitem{\rm(c)} if $\Delta'_j$ is a logarithmic component
$(\rho'_j=\infty)$, then $\Phi^{-1}(\Delta'_j)=
\bigcup_{i\in I(j)}\Delta_i$ where the $(\Delta_i)_{i\in I(j)}$
consist of logarithmic components $(\rho_i=\infty)$.\vskip0pt
\endclaim

\noindent
It is easy to check that, if the image of the composed morphism is not contained in the support of the divisor on the target space, the composite of directed orbifold morphisms
is actually a directed orbifold morphism, and that the composition
of an orbifold entire curve $f:\bC\to(X,V,D)$ with a
directed orbifold morphism $\Psi:(X,V,D)\to(Y,W,D')$ 
produces an orbifold entire curve $\Psi\circ f:\bC\to(Y,W,D')$ (provided that $\Psi\circ f(\bC) \not \subset |D'|$).
One of our main goals is to investigate the following orbifold
generalization of the Green-Griffiths conjecture.

\claim 0.5. Conjecture|Let $(X,V,D)$ be a non singular directed orbifold
of general type, in the sense that the canonical divisor $K_V+ D$
is big. Then then exists an algebraic subvariety $Y\subsetneq X$
containing all orbifold entire curves $f:\bC\to(X,V,D)$.
\endclaim

\noindent
As in the absolute case ($V=T_X$, $D=0$), the idea is to show, at least
as a first step towards the conjecture, that orbifold entire curves must satisfy
suitable algebraic differential equations. In section~1, we introduce graded
algebras
$$
\bigoplus_{m\in\bN}E_{k,m}V^*\langle D\rangle\leqno(0.6)
$$
of sheaves of ``orbifold jet differentials''. These sheaves correspond to
algebraic differential operators $P(f;f',f'',\ldots,f^{(k)})$ acting on
germs of $k$-jets of curves that are tangent to $V$ and satisfy the
ramification conditions prescribed by~$D$. The strategy relies on the
following orbifold version of the vanishing theorem, whose proof is 
sketched in the appendix.

\claim 0.7. Proposition|Let $(X,V,D)$ be a projective non singular
directed orbifold, and $A$~an ample divisor on $X$. Then, for every
orbifold entire curve $f:\bC\to(X,V,D)$ and every global
jet differential operator $P\in H^0(X,E_{k,m}V^*\langle D\rangle
\otimes\cO_X(-A))$, we have~\hbox{$P(f;f',f'',\ldots,f^{(k)})=0$}.
\endclaim

\noindent
The next step consists precisely of finding sufficient conditions that
ensure the existence of global sections
$P\in H^0(X,E_{k,m}V^*\langle D\rangle\otimes\cO_X(-A))$. Recall that
it has been shown in [CDR20, Proposition 5.1] that the general type
assumption is not a sufficient condition for the existence of global
jet differentials.

\noindent
Among more general results, we obtain

\claim 0.8. Theorem|Let $D=\sum_j(1-{1\over\rho_j})\Delta_j$
a simple normal crossing orbifold divisor on~$\bP^n$ with
$\deg\Delta_j=d_j$. Then there exist non zero jet differentials of order $k$ 
and large degree $m$ on $\bP^n\langle D\rangle$, with a small negative
twist $\cO_{\bP^n}(-m\tau)$, $\tau>0$, under any of the 
following two sufficient conditions~$:$\vskip6pt
\plainitem{\rm(a)}
$\,k\ge n,\quad N\ge 1,\quad \rho_j\ge \rho>n$~~ and\vskip-8pt
$$
\sum_j d_j\cdot
\min\bigg(\!\min_j\bigg({\rho_j\over d_j}\bigg),{1\over 2}\bigg)
\,\prod_{s=1}^n\Big(1-{s\over\rho}\Big)>c_n
$$
where
$$
c_n:=n(n^2+n-1)\,n!\,\bigg(\sum_{s=1}^n{1\over s}+{1\over n^3}\bigg)^{n-1}
\sim(2\pi)^{1/2}\,
n^{n+7/2}\,e^{-n}(\gamma+\log n)^{n-1}.
$$
\plainitem{\rm(b)}
$\,k\ge 1,\quad N\ge n,\quad \rho_j\ge \rho>1$~~ and
for~ $t=\max(\max(d_j/\rho_j),2),$\vskip-10pt
$$
\sum_{J\subset\{1,\ldots,N\},\,|J|=n}~
\prod_{j\in J}d_j\bigg(1-{1\over \rho_j}\bigg)>
(2n-1)\,t\,
\Big(n\,t-n-1+\sum_jd_j(1-1/\rho_j)\Big)^{n-1}.
$$

\noindent
When all components $(\Delta_j)_{1\le j\le N}$ possess
the same degrees $d_j=d\ge 1$ and ramification numbers $\rho_j\ge\rho$,
we get the following simpler sufficient conditions~$:$\vskip6pt
\plainitem{{\rm(${\rm a}'$)}}
$\displaystyle~\,k\ge n,~~N\ge 1,~~\rho>n,\quad
N\min(\rho,d)\,\prod_{s=1}^n\Big(1-{s\over\rho}\Big)>2c_n,$
\vskip6pt
\plainitem{{\rm(${\rm b}'$)}}
$\displaystyle~\,k\ge 1,~~N\ge n,~~\rho>1,\quad
N\min(\rho,d)\,\Big(1-{1\over\rho}\Big)^n>2^n\,(2n-1)\,n^n.$
\endclaim

\noindent
Let us recall some related results previously obtained in this
orbifold setting. In the case of orbifold surfaces 
$\big(\bP^2,\big(1-{1\over\rho}\big)C\big)$ where $C$ is a smooth curve
of degree $d$, such existence results have been obtained in [CDR20]
for $k=2$, $d\ge 12$ and $\rho \ge 5$ depending on $d$. 
In [DR20], the existence of jet differentials is obtained for 
orbifolds $\big(\bP^n,\sum_{i=1}^d\big(1-{1\over\rho}\big)H_i\big)$
in any dimension for $k=1$, $\rho \ge 3$ along an arrangement of 
hyperplanes of degree $d\ge 2n\big({2n \over \rho-2} + 1\big)$.
In [BD18], it is established that the orbifold
$\big(\bP^n,\big(1-{1\over d}\big)D\big)$, where $D$ is a general
smooth hypersurface of degree $d$, is hyperbolic i.e.\ there is no
non-constant orbifold entire curve $f:\bC \to\big(\bP^n,
\big(1-{1\over d}\big)D\big)$, if $d\ge (n+2)^{n+3}(n+1)^{n+3}$.

\noindent
The proof of Theorem 0.8 depends on a number of ingredients and on
rather extensive curvature calculations. The first point is
that the curvature tensor of the orbifold directed structure
$V\langle D\rangle$ can be controlled in a precise manner.
This is detailed in \S$\,$6.A.

\claim 0.9. Theorem|Assume that $X$ is projective. Given an
an ample line bundle $A$ on $X$, let $\gamma_V$ be the infimum
of real numbers $\gamma\ge 0$ such that $\gamma\,\Theta_A\otimes\Id_V
-\Theta_V$ is positive in the sense of Griffiths, for suitable $C^\infty$
smooth hermitian metrics on $V$. Assume that $D=\sum_j(1-1/\rho_j)\Delta_j$
is transverse to $V$, and select $d_j\ge 0$ such that $d_jA-\Delta_j$
is~nef. Then for $\gamma>\gamma_{V,D}:=\max(\max(d_j/\rho_j),\gamma_V)\ge 0$
and for suitable hermitian metrics on $A$,~$V$, $\cO_X(\Delta_j)$,
the ``orbifold metric''
\vskip5pt\noindent
{\rm(a)}~~$\displaystyle
|u|_{h_{V\langle D\rangle,\varepsilon}}^2:=
|u|^2_{h_V}+\sum_{1\le j\le N}\varepsilon_j\,|\sigma_j|^{-2+2/\rho_j}\,
|\nabla_j\sigma_j(u)|^2_{h_j},~~u\in V,~~\sigma_j\in H^0(X,\cO_X(\Delta_j))$
\vskip5pt\noindent
yields a curvature tensor $\gamma\,\Theta_A\otimes\Id-\Theta_{V\langle D\rangle}$
such that the associated quadratic form $Q_{V\langle D\rangle,\gamma,\varepsilon}$ on
$T_X\otimes V$ satisfies for
$\varepsilon_N\ll\varepsilon_{N-1}\ll\cdots\ll\varepsilon_1\ll 1$
the curvature estimate
$$
\leqalignno{\kern10pt Q_{V\langle D\rangle,\gamma,\varepsilon}
&(z)(\xi\otimes u)\simeq\gamma\,\Theta_A(\xi,\xi)\,|u|^2
-\langle\Theta_V(\xi,\xi)\cdot u,u\rangle&{\rm(b)}\cr
+\sum_j&\varepsilon_j\,|\sigma_j|^{-2+2/\rho_j}\,
\big(\gamma\,\Theta_A(\xi,\xi)-\rho_j^{-1}\,\Theta_{\Delta_j}(\xi,\xi)\big)
\,|\nabla_j\sigma_j(u)|^2\cr
+\sum_j&{\varepsilon_j\,|\sigma_j|^{-2+2/\rho_j}\over
1+\varepsilon_j\,|\sigma_j|^{-2+2/\rho_j}\,
|\nabla_j\sigma_j|^2}\,\big|\,\nabla_j^2\sigma_j(\xi,u)-(1-1/\rho_j)\,
\sigma_j^{-1}\,\nabla_j\sigma_j(\xi)\nabla_j\sigma_j(u)\big|^2.\cr}
$$
Here, the symbol $\simeq$ means that the ratio of the left and right hand
sides can be chosen in $[1-\alpha,1+\alpha]$ for any
$\alpha>0$ prescribed in advance.
\endclaim

\noindent The next argument is the observation that the sheaf
$\cO_X(E_{k,m}V^*\langle D\rangle)$ is the direct image of
a certain tautological rank $1$ sheaf $\cO_{X_k(V\langle D\rangle)}(m)$
on the ``orbifold $k$-jet bundle'' \hbox{$X_k(V\langle D\rangle)\to X$}.
Choosing hermitian metrics according to Theorem 0.9, one then gets a
hermitian metric on  $\cO_{X_k(V\langle D\rangle)}(1)$ associated with
an ``orbifold Finsler metric'' on the bundle $J_kV$ of $k$-jets of
holomorphic curves $f:(\bC,0)\to(X,V)$. In normalized coordinates
$(z_1,\ldots,z_n)$ on $X$ and on $V$, the latter can be expressed as
$$
\Bigg(\sum_{s=1}^k\varepsilon_s^{2b}\bigg(\sum_{j=1}^p
|f_j|^{-2(1-s/\rho_j)_+}|f_j^{(s)}|^2+\sum_{j=p+1}^r|f_j^{(s)}|^2\bigg)^{2b/s}\Bigg)^{1/b},\quad f\in J^kV,~~f(0)=x,
\leqno(0.10)
$$
at any point $x\in X$ where $\Delta_j=\{z_j=0\}$, $1\le j\le p$, $r=\rank V$.
An application of holomorphic
Morse inequalities ([Dem85], see also \S$\,2,3,4$) then provides asymptotic 
estimates of the dimensions of the cohomology groups
$$
H^q(X,E_{k,m}V^*\langle D\rangle\otimes\cO_X(-A))\simeq
H^q(X_k(V\langle D\rangle),
\cO_{X_k(V\langle D\rangle)}(m)\otimes\pi_k^*\cO_X(-A)).
\leqno(0.11)
$$
This is done in several steps. Section \S$\,$4 expresses the Morse integrals
that need to be computed. Section \S$\,$5 establishes some general estimates
of Chern forms related to the curvature tensor $\Theta_{E,h}$ of a given
hermitian vector bundle $(E,h)$, under suitable positivity assumptions.
More precisely, Proposition~5.13 gives upper and lower bounds of
integrals of the form
$$
\int_{u\in S(E)}|\ell_1(u)|^2\ldots|\ell_k(u)|^2\,
\langle\Theta_{E,h}(u),u\rangle_h^{p-k}\,d\mu(u)
\leqno(0.12)
$$
in terms of $\Tr_E\Theta_{E,h}=\Theta_{\det E,\det h}$,
where $\mu$ is the unitary invariant probability measure on the unit
sphere bundle $S(E)$, and the $\ell_j$ are linear forms. As far as we
know, these estimates seem to be new. Sections \S$\,$6.B and \S$\,$7
then proceed with the detailed calculations of the orbifold and
logarithmic Morse integrals involved in the problem.  It is remarkable
that a large part of the calculations use Chern forms and are non
cohomological, although the final bounds are purely cohomological.
At~this point, we do not have a complete explanation of this
``transcendental'' phenomenon.

\plainsection{1. Logarithmic and orbifold jet differentials}

\plainsubsection 1.A. Directed varieties and associated jet differentials|

Let $(X,V)$ be a non singular directed variety. We set $n=\dim_\bC X$,
$r=\rank_\bC V$, and following the exposition of [Dem97], we
denote by $\pi_k:J^kV\to X$ the bundle
of $k$-jets of holomorphic curves tangent to
$V$ at each point. The canonical bundle of $V$ is defined to be
$$
K_V=\det(V^*)=\Lambda^rV^*.\leqno(1.1)
$$
If $f:(\bC,0)\to X$, $t\mapsto f(t)$
is a germ of holomorphic curve tangent to $V$, we denote
by $f_{[k]}(0)$ its $k$-jet at ~$t=0$. For $x_0\in X$ given, we take a
coordinate system $(z_1,\ldots,z_n)$ centered at $x_0$ such that
$V_{x_0}=\Span({\partial\over\partial z_\mu})_{1\le \mu\le r}$.
Then there exists a neighborhood $U$ of $x_0$ such that
$V_{|U}$ admits a holomorphic frame $(e_\mu)_{1\le\mu\le r}$ of the form
$$
e_\mu(z)={\partial\over\partial z_\mu}+\sum_{r+1\le \lambda\le n}
a_{\lambda\mu}(z){\partial\over\partial z_\lambda},\quad
1\le\mu\le r,\leqno(1.2)
$$
with $a_{\lambda\mu}(0)=0$. Germs of curves $f:(\bC,0)\to X$ tangent to $V_{|U}$
are obtained by integrating the system of ordinary differential equations
$$
f'_\lambda(t)=\sum_{1\le\mu\le r}a_{\lambda\mu}(f(t))\,f'_\mu(t),\quad
r+1\le \lambda\le n,\leqno(1.3)
$$
when we write $f=(f_1,\ldots,f_n)$ in coordinates. Therefore any such germ of
curve $f$ is uniquely determined by its initial point $z=f(0)$ and its
projection $\tilde f=(f_1,\ldots,f_r)$ on the first $r$ coordinates. By
definition, every $k$-jet $f_{[k]}\in J^kV_z=\pi_k^{-1}(z)$ is
uniquely determined
by its initial point $f(0)=z\simeq(z_1,\ldots,z_n)$ and the Taylor expansion
of order $k$
$$
\tilde f(t)-\tilde f(0)=t\xi_1+{1\over 2!}t^2\xi_2+\cdots+{1\over k!}
t^k\xi_k+O(t^{k+1}),\quad t\in\bD(0,\varepsilon),~\xi_s\in\bC^r,~1\le s\le k.
\leqno(1.4)
$$
Alternatively, we can pick an arbitrary local holomorphic connection $\nabla$
on $V_{|U}$ and represent the $k$-jet $f_{[k]}(0)$ by
$(\xi_1,\ldots,\xi_k)$, where
$\xi_s=\nabla^sf(0)\in V_z$ is defined inductively 
by $\nabla^1 f=f'$ and $\nabla^sf=\nabla_{f'}(\nabla^{s-1}f)$. This
gives a local biholomorphic trivialization of $J^kV_{|U}$ of the form
$$
J_kV_{|U}\to V_{|U}^{\oplus k},\qquad
f_{[k]}(0)\mapsto(\xi_1,\ldots,\xi_k)=(\nabla f(0),\ldots,\nabla f^k(0))\,;
\leqno(1.5)
$$
the particular choice of the ``trivial connection'' $\nabla_0$ of $V_{|U}$
that turns $(e_\mu)_{1\le\mu\le r}$ into a parallel frame precisely yields the
components $\xi_s\in V_{|U}\simeq\bC^r$ appearing in (1.4). We could of
course also use a $C^\infty$ connection $\nabla=\nabla_0+\Gamma$ where
$\Gamma\in C^\infty(U,T^*_X\otimes\Hom(V,V))$, and
in this case, the corresponding trivialization (1.5) is just a
$C^\infty$ diffeomorphism; the advantage, though, is that we can always
produce such a global $C^\infty$ connection $\nabla$ by using a partition of
unity on~$X$, and then (1.5) becomes a global $C^\infty$ diffeomorphism.
Now, there is a global holomorphic $\bC^*$ action on $J^kV$ given at the level
of germs by $f\mapsto\alpha\cdot f$ where $\alpha\cdot f(t):=f(\alpha t)$,
$\alpha\in\bC^*$. With respect to our trivializations (1.5), this is
the weighted $\bC^*$ action defined by
$$
\alpha\cdot(\xi_1,\xi_2,\ldots,\xi_k)=(\alpha\xi_1,
\alpha^2\xi_2,\ldots,\alpha^k\xi_k),\quad \xi_s\in V.\leqno(1.6)
$$
We see that $J^kV\to X$ is an algebraic fiber bundle
with typical fiber $\bC^{rk}$, and that the projectivized $k$-jet bundle 
$$
X_k(V):=(J^kV\ssm\{0\})/\bC^*,\qquad \pi_k:\smash{X_k(V)}\to X
\leqno(1.7)
$$
is a $\bP(1^{[r]},2^{[r]},\ldots,k^{[r]})$ weighted projective 
bundle over $X$, of total dimension
$$
\dim X_k(V)=n+kr-1.
\leqno(1.8)
$$

\claim 1.9. Definition|We define $\cO_X(E_{k,m} V^*)$ to be the sheaf over
$X$ of holomorphic functions $P(z\,;\,\xi_1,\ldots,\xi_k)$ on $J^kV$ that are
weighted polynomials of degree $m$ in~$(\xi_1,\ldots,\xi_m)$.
\endclaim

\noindent In coordinates and in multi-index notation, we can write
$$
P(z\,;\,\xi_1,\ldots,\xi_k)=
\sum_{\scriptstyle\alpha_1,\ldots,\alpha_k\in\bN^r\atop
\scriptstyle|\alpha_1|+2|\alpha_2|+\cdots+k|\alpha_k|=m}
a_{\alpha_1\ldots\alpha_k}(z)\,\xi_1^{\alpha_1}\ldots\xi_k^{\alpha_k}
$$
where the $a_{\alpha_1\ldots\alpha_k}(z)$ are holomorphic functions
in $z=(z_1,\ldots,z_n)$ and $\xi_s^{\alpha_s}$ actually means
$$
\xi_s^{\alpha_s}=\xi_{s,1}^{\alpha_{s,1}}\ldots\,\xi_{s,r}^{\alpha_{s,r}}\quad
\hbox{for}~~
\xi_s=(\xi_{s,1},\ldots,\xi_{s,r})\in\bC^r,~~
\alpha_s=(\alpha_{s,1},\ldots,\alpha_{s,r})\in\bN^r,
$$
and $|\alpha_s|=\sum_{j=1}^r\alpha_{s,j}$. Such sections can be
interpreted as algebraic differential operators acting on holomorphic
curves $f:\bD(0,R)\to X$ tangent to $V$, by putting
$P(f):=u$ where
$$
u(t)=\sum_{\scriptstyle\alpha_1,\ldots,\alpha_k\in\bN^r\atop
\scriptstyle|\alpha_1|+2|\alpha_2|+\cdots+k|\alpha_k|=m}
a_{\alpha_1\ldots\alpha_k}(f(t))\;f'(t)^{\alpha_1}\ldots\,f^{(k)}(t)^{\alpha_k}.
\leqno(1.10)
$$
Here $f^{(s)}(t)^{\alpha_s}$ is actually to be expanded as
$$
f^{(s)}(t)^{\alpha_s}=f_1^{(s)}(t)^{\alpha_{s,1}}\ldots\,f_r^{(s)}(t)^{\alpha_{s,r}}
$$
with respect to the components $f_j^{(s)}$ defined in (1.4). We also
set $u=P(f\,;\,f',f'',\ldots,f^{(k)})$ when we want to make more explicit
the dependence of the expression in terms of the derivatives of~$f$.
We thus get a sheaf of graded algebras
$$
\bigoplus_{m\in\bN}\cO_X(E_{k,m}V^*).\leqno(1.11)
$$
Locally in coordinates, the algebra is isomorphic
to the weighted polynomial ring
$$
\cO_X\big[f_j^{(s)}\big]_{1\le j\le r,\,1\le s\le k},\quad
\deg f_j^{(s)}=s\leqno(1.12)
$$
over $\cO_X$. An immediate consequence of these definitions is~:

\claim 1.13. Proposition|The projectivized bundle $\pi_k:X_k(V)\to X$ can be
identified with 
$$
\Proj\Bigg(\bigoplus_{m\in\bN}\cO_X(E_{k,m}V^*)\Bigg)\to X,
\leqno{\rm(a)}
$$
and, if $\cO_{X_k(V)}(m)$ denote the associated tautological sheaves,
we have the direct image formula
$$
(\pi_k)_*\smash{\cO_{X_k(V)}(m)}=\cO_X(E_{k,m} V^*).
\leqno{\rm(b)}
$$
\endclaim

\claim 1.14. Remark|{\rm These
objects were denoted $X_k^\GG$ and $E_{k,m}^\GG V^*$ in our previous
paper [Dem97], as a reference to the work of Green-Griffiths [GG79],
but we will avoid here the superscript GG to
simplify the notation.}
\endclaim

\noindent
Thanks to the Fa\`a di Bruno formula, a change of coordinates $w=\psi(z)$ on
$X$ leads to a transformation rule
$$
(\psi\circ f)^{(k)}=\psi'\circ f\cdot f^{(k)}+Q_\psi(f',\ldots,f^{(k-1})
$$
where $Q_\psi$ is a polynomial of weighted degree $k$ in the lower order
derivatives. This shows that the transformation rule of the top derivative
is linear and, as a consequence, the partial degree in $f^{(k)}$ of
the polynomial $P(f\,;\,f',\ldots,f^{k)})$ is intrinsically defined.
By taking the corresponding filtration and factorizing the monomials
$(f^{(k)})^{\alpha_k}$ with polynomials in $f',f'',\ldots,f^{(k-1)}$,
we get graded pieces
$$
G^\bullet(E_{k,m}V^*)=\bigoplus_{\ell_k\in\bN}
E_{k-1,m-k\ell_k}V^*\otimes S^{\ell_k}V^*.
$$
By considering successively the partial degrees with respect to
$f^{(k)}$, $f^{(k-1)}$, $\ldots\,$, $f'',f'$ and merging inductively
the resulting filtrations, we get a multi-filtration
such that
$$
G^\bullet(E_{k,m}V^*)=\bigoplus_{\ell_1,\ldots,\ell_k\in\bN,\,
\ell_1+2\ell_2+\cdots+k\ell_k=m}S^{\ell_1}V^*\otimes S^{\ell_2}V^*\otimes\cdots
\otimes S^{\ell_k}V^*.\leqno(1.15)
$$

\plainsubsection 1.B. Logarithmic directed varieties|

We now turn ourselves to the logarithmic case. Let $(X,V,D)$ be a
non singular logarithmic variety, where $D=\sum\Delta_j$ is a simple
normal crossing divisor. Fix a point $x_0\in X$. By the assumption that
$D$ is transverse to $V$, we can then select holomorphic coordinates
$(z_1,\ldots,z_n)$ centered at $x_0$ such that
$V_{x_0}=\Span({\partial\over\partial z_j})_{1\le j\le r}$
and $\Delta_j=\{z_j=0\}$, $1\le j\le p$, are the components of $D$
that contain $x_0$ (here $p\le r$ and we can have $p=0$
if $x_0\notin|D|$). What we want is to introduce an algebra of
differential operators, defined locally near $x_0$ as the weighted
polynomial ring
$$
\cO_X\big[(\log f_j)^{(s)}_{1\le j\le p}\,,(f_j^{(s)})_{p+1\le j\le r}
\big]_{1\le s\le k},\quad \deg f_j^{(s)}=\deg(\log f_j)^{(s)}=s,\leqno(1.16)
$$
or equivalently
$$
\cO_X\big[(f_j^{-1}f_j^{(s)})_{1\le j\le p}\,,(f_j^{(s)})_{p+1\le j\le r}
\big]_{1\le s\le k},\quad \deg f_j^{(s)}=s,~\deg f_j^{-1}=0.\leqno(1.16')
$$
For this we notice that
$$
\eqalign{
(\log f_1)''&=(f_1^{-1}f_1')'=f_1^{-1}f_1''-(f_1^{-1}f_1')^2,\cr
\noalign{\vskip4pt}
(\log f_1)'''&=f_1^{-1}f_1'''-3(f_1^{-1}f_1')(f_1^{-1}f_1'')+2
(f_1^{-1}f_1')^3,\ldots\,.\cr}
$$
A similar argument easily shows that the above graded rings do not depend on
the particular choice of coordinates made, as soon as they satisfy
$\Delta_j=\{z_j=0\}$.

Now (as is well known in the absolute case $V=T_X$), we have a
corresponding logarithmic directed structure
$V\langle D\rangle$ and its dual $V^*\langle D\rangle$.
If the coordinates $(z_1,\ldots,z_n)$ are
chosen so that $V_{x_0}=\{dz_{r+1}=\ldots=dz_n=0\}$, then
the fiber $V\langle D\rangle_{x_0}$ is spanned by the derivations
$$
z_1{\partial\over\partial z_1},\ldots,z_p{\partial\over\partial z_p},~
{\partial\over\partial z_{p+1}},\ldots,{\partial\over\partial z_r}.
$$
The dual sheaf $\cO_X(V^*\langle D\rangle)$ is the
locally free sheaf generated by
$$
{dz_1\over z_1},\ldots,{dz_p\over z_p},~dz_{p+1},\ldots,dz_r
$$
[where the $1$-forms are considered in restriction to
$\cO_X(V\langle D\rangle)\subset\cO_X(V)\,$]. It follows from this
that $\cO_X(V\langle D\rangle)$ and
$\cO_X(V^*\langle D\rangle)$ are locally free sheaves of rank~$r$.
By taking $\det(V^*\langle D\rangle)$ and using the above generators,
we find
$$
\det(V^*\langle D\rangle)=\det(V^*)\otimes\cO_X(D)=K_V+ D
\leqno(1.17)
$$
in additive notation. Quite similarly to 1.13 and 1.15, we have~:

\claim 1.18. Proposition|Let $\bigoplus_{m\in\bN}\cO_X(E_{k,m}V^*
\langle D\rangle)$ be the graded algebra defined in coordinates by
$(1.16)$ or $(1.16')$. We define the logarithmic $k$-jet bundle to be
$$
X_k(V\langle D\rangle):=
\Proj\Bigg(\bigoplus_{m\in\bN}\cO_X(E_{k,m}V^*\langle D\rangle)\Bigg)\to X.
\leqno{\rm(a)}
$$
If $\cO_{X_k(V\langle D\rangle)}(m)$ denote the
associated tautological sheaves, we get the direct image formula
$$
(\pi_k)_*\smash{\cO_{X_k(V\langle D\rangle)}(m)}=\cO_X(E_{k,m} V^*
\langle D\rangle).\leqno{\rm(b)}
$$
Moreover, the multi-filtration by the partial degrees in the derivatives
$f_j^{(s)}$ has graded pieces
$$
G^\bullet\big(E_{k,m}V^*\langle D\rangle\big)=
\bigoplus_{\ell_1,\ldots,\ell_k\in\bN,\, \ell_1+2\ell_2+\cdots+k\ell_k=m}
S^{\ell_1}V^*\langle D\rangle\otimes
S^{\ell_2}V^*\langle D\rangle\otimes\cdots\otimes
S^{\ell_k}V^*\langle D\rangle.
\leqno{\rm(c)}
$$
\endclaim

\plainsubsection 1.C. Orbifold directed varieties|

We finally consider a non singular directed orbifold $(X,V,D)$,
where $D=\sum(1-{1\over\rho_j})\Delta_j$ is a simple normal crossing
divisor transverse to~$V$. Let $\lceil D\rceil=\sum\Delta_j$ be
the corresponding reduced divisor. By \S$\,$1.B, we have associated
logarithmic sheaves $\cO_X(E_{k,m}V^*\langle\lceil D\rceil\rangle)$.
We want to introduce a graded subalgebra
$$
\bigoplus_{m\in\bN}\cO_X(E_{k,m}V^*\langle D\rangle)~~\subset~~
\bigoplus_{m\in\bN}\cO_X(E_{k,m}V^*\langle\lceil D\rceil\rangle)
\leqno(1.19)
$$
in such a way that for every germ $P\in \cO_X(E_{k,m}V^*\langle D\rangle)$
and every germ of orbifold curve $f:(\bC,0)\to(X,V,D)$ the germ of
meromorphic function $P(f)(t)$ is bounded at $t=0$ (hence holomorphic).
Assume that $\Delta_1=\{z_1=0\}$ and that $f$ has multiplicity $q\ge \rho_1>1$
along~$\Delta_1$ at $t=0$. Then $f_1^{(s)}$ still vanishes
at order${}\ge(q-s)_+$, thus $(f_1)^{-\beta}f_1^{(s)}$ is bounded as soon
as $\beta q\le(q-s)_+$, i.e.\ $\beta\le(1-{s\over q})_+$. Thus,
it is sufficient to ask that $\beta\le(1-{s\over \rho_1})_+$. At a point
$x_0\in |\Delta_1|\cap\ldots\cap|\Delta_p|$, a sufficient condition
for a monomial of the form
$$
f_1^{-\beta_1}\ldots\,f_p^{-\beta_p}\prod_{s=1}^k\prod_{j=1}^r
(f_j^{(s)})^{\alpha_{s,j}},
\quad
\alpha_s=(\alpha_{s,j})\in\bN^r,~\beta_1,\ldots,\beta_p\in\bN
\leqno(1.20)
$$
to be bounded is to require that the multiplicities of poles satisfy
$$
\beta_j\le\sum_{s=1}^k\alpha_{s,j}\Big(1-{s\over \rho_j}\Big)_+,\quad
1\le j\le p.
\leqno(1.20')
$$
\claim 1.21. Definition|The subalgebra
$\bigoplus_{m\in\bN}\cO_X(E_{k,m}V^*\langle D\rangle)$ is taken to be the
graded ring generated by monomials
$(1.20)$ of degree $\sum s|\alpha_s|=m$, satisfying the pole multiplicity
conditions $(1.20')$. These conditions do not depend on the choice of
coordinates, hence we get a globally and intrinsically defined sheaf
of algebras on~$X$. 
\endclaim

\plainproof. We only have to prove the last assertion. Consider a change of
variables $w=\psi(z)$ such that $\Delta_j$ can still be expressed as
$\Delta_j=\{w_j=0\}$. Then, for $j=1,\ldots,p$, we can write
$w_j=z_ju_j(z)$ with an invertible holomorphic factor~$u_j$. We need to check
that the monomials~(1.20) computed with $g=\psi\circ f$ are holomorphic
combinations of those associated with $f$. However, we have $g_j=f_ju_j(f)$,
hence $g_j^{(s)}=\sum_{0\le\ell\le s}{s\choose\ell}
f_j^{(\ell)}(u_j(f))^{(s-\ell)}$ by the Leibniz formula, and we see that
$$
g_1^{-\beta_1}\ldots\,g_p^{-\beta_p}\prod_{s=1}^k\prod_{j=1}^r
(g_j^{(s)})^{\alpha_{s,j}}
$$
expands as a linear combination of monomials
$$
f_1^{-\beta_1}\ldots\,f_p^{-\beta_p}\prod_{s=1}^k\prod_{j=1}^r
\prod_{m=1}^{\alpha_{s,j}}f_j^{(\ell_{s,j,m})},\quad \ell_{s,j,m}\le s,
$$
multiplied by holomorphic factors of the form
$$
\prod_{j=1}^pu_j(f)^{-\beta_j}\times
\prod_{s=1}^k\prod_{j=1}^r\prod_{m=1}^{\alpha_{s,j}}(u_j(f))^{(s-\ell_{j,s,m})}.
$$
However, we have
$$
\beta_j\le\sum_{s=1}^k\alpha_{s,j}\Big(1-{s\over \rho_j}\Big)_+
\le\sum_{s=1}^k\sum_{m=1}^{\alpha_{s,j}}
\Big(1-{\ell_{s,j,m}\over \rho_j}\Big)_+,\quad
$$
so the $f$-monomials satisfy again the required multiplicity conditions
for the poles~$f_j^{-\beta_j}$.\qed

\noindent
The above conditions $(1.20')$ suggest to introduce as in [CDR20]
a sequence of ``differentiated'' orbifold divisors
$$
 D^{(s)}=\sum_j\bigg(1-{s\over\rho_j}\bigg)_{\kern-3pt+}\Delta_j.
\leqno(1.22)
$$
We say that $D^{(s)}$ is the order $s$ orbifold divisor associated
to~$D\,$; its ramification numbers are $\rho_j^{(s)}=\max(\rho_j/s,1)$.
By definition, the logarithmic components ($\rho_j=\infty$) of
$D$ remain logarithmic in $D^{(s)}$, while all others eventually
disappear when $s$ is large.

Now, we introduce (in a purely formal way) a sheaf of rings
$\smash{\wt\cO}_X=\cO_X[z_j^\bullet]$ by adjoining
all positive real powers of coordinates $z_j$ such that
\hbox{$\Delta_j=\{z_j=0\}$} is locally a component of~$D$.
Locally over~$X$, this can be done by taking the universal cover $Y$ of
a punctured polydisk
$$
\bD^*(0,r):=\prod_{1\le j\le p}\bD^*(0,r_j)\times
\prod_{p+1\le j\le n}\bD(0,r_j)~~\subset~~
\bD(0,r):=\prod_{1\le j\le n}\bD(0,r_j)
$$
in the local coordinates $z_j$ on $X$. If $\gamma:Y\to\bD^*(0,r)
\hookrightarrow X$ is the covering map and $U\subset\bD(0,r)$ is an open
subset, we can then consider the functions
of~$\smash{\wt\cO}_X(U)$ as being defined on
$\gamma^{-1}(U\cap\bD^*(0,r))$. In case $X$ is projective,
one can even achieve such a construction ``globally'', at least on a
Zariski open set, by taking $Y$ to be the universal cover of a
complement $X\ssm(|D|\cup|A|)$, where
$A=\sum A_j$ is a very ample normal crossing divisor transverse to $D$, 
such that $\cO_X(\Delta_j)_{|X\ssm|A|}$ is trivial for every $j\,$; then
$\smash{\wt\cO}_X$ is well defined as a genuine sheaf on $X\ssm|A|$.

In this setting,
the subalgebra $\bigoplus_m\cO_X(E_{k,m}V^*\langle D\rangle)$ still has a
multi-filtration induced
by the one on $\bigoplus_m\cO_X(E_{k,m}V^*\langle\lceil D\rceil\rangle)$,
and by extending the structure sheaf $\cO_X$ into $\smash{\wt\cO}_X$, we get
an inclusion
$$
\wt\cO_X(G^\bullet E_{k,m}V^*\langle D\rangle)\subset
\bigoplus_{\ell_1+2\ell_2+\cdots+k\ell_k=m}
\wt\cO_X(S^{\ell_1}V^*\langle D^{(1)})\rangle\otimes \cdots
\otimes\wt\cO_X(S^{\ell_k}V^*\langle D^{(k)}\rangle),
\leqno(1.23)
$$
$\wt\cO_X(V^*\langle D^{(s)}\rangle)$ is the ``$s$-th orbifold
(dual) directed structure'', generated by the order
$s$ differentials
$$
z_j^{-(1-s/\rho_j)_+}d^{(s)}z_j,~~1\le j\le p,~~~d^{(s)}z_j,~~p+1\le j\le r.
\leqno(1.24)
$$
By construction, we have
$$
\det(\wt\cO_X(V^*\langle D^{(s)}\rangle))=\wt\cO_X(K_V+ D^{(s)}).
\leqno(1.25)
$$

\claim 1.26. Remark|{\rm When $\rho_j=a_j/b_j\in\bQ_+$, one can 
find a finite ramified Galois cover \hbox{$g:Y\to X$} from a smooth
projective variety $Y$ onto~$X$, such that the compositions
$(z_j\circ g)^{1/a_j}$ become single-valued functions $w_j$ on $Y$.
In this way, the pull-back
$\cO_Y(g^*V^*\langle D^{(s)}\rangle)$ is actually a
locally free $\cO_Y$-module. On can also introduce a sheaf
of algebras which we will denote by
$\bigoplus\cO_Y(E_{k,m}\widetilde V^*\langle D\rangle)$,
generated, according to the notation of \S$\,$1.B, by the elements
$g^*(z_j^{(1-s/\rho_j)_+}d^{(s)}z_j)$, $1\le j\le p$, and
$g^*(d^{(s)}z_j)$, $p+1\le j\le r$. Then, as already shown 
in [CDR20], there is indeed a multifiltration on 
$\cO_Y(E_{k,m}\widetilde V^*\langle D\rangle)$ whose 
graded pieces are
$$
\cO_Y(G^\bu E_{k,m}\widetilde V^*\langle D\rangle)=
\bigoplus_{\ell_1+2\ell_2+\cdots+k\ell_k=m}
\cO_Y(S^{\ell_1}\wt V^*\langle D^{(1)}\rangle)\otimes \cdots
\otimes\cO_Y(S^{\ell_k}\wt V^*\langle D^{(k)}\rangle).
\leqno(1.27)
$$
However, we will adopt here an alternative viewpoint that avoids the
introduction of finite or infinite covers, and suits better our
approach. Using the general philosophy of [Laz??], the idea is
to consider a ``jet orbifold directed
structure'' $X_k(V\langle D\rangle)$ as the underlying
``jet logarithmic directed structure''
$X_k(V\langle\lceil D\rceil\rangle)$, equipped additionally
with a submultiplicative sequence of ideal sheaves
$\cJ_m\langle D\rangle\subset \cO_{X_k(V\langle\lceil D\rceil\rangle)}$.
These are precisely defined as
the base loci ideals of the local sections defined by $(1.20)$ and $(1.20')$,
seen as sections of the logarithmic tautological sheaves
$\cO_{X_k(V\langle\lceil D\rceil\rangle)}(m)$. The corresponding analytic
viewpoint is to consider ad hoc singular
hermitian metrics on $\cO_{X_k(V\langle\lceil D\rceil\rangle)}(1)$ whose
singularities are asymptotically described by the limit of the formal
$m$-th root of $\cJ_m\langle D\rangle$, see \S$\,$3.B. It then becomes
possible to deal
without trouble with real coefficients $\rho_j\in{}]1,\infty]$, and
since we no longer have to worry about the existence of Galois covers,
the projectivity assumption on $X$ can be dropped as well.}
\endclaim

\plainsection{2. Preliminaries on holomorphic Morse inequalities}

\plainsubsection 2.A. Basic results|

We first recall the basic results concerning holomorphic Morse inequalities
for smooth hermitian line bundles, first proved in [Dem85].

\claim 2.1. Theorem|Let $X$
be a compact complex  manifolds, $E\to X$ a holomorphic vector bundle of
rank $r$, and $(L,h)$ a hermitian line bundle. We denote by
$\Theta_{L,h}={\ii\over 2\pi}\nabla_h^2=-{\ii\over2\pi}\ddbar\log h$ the
curvature form of $(L,h)$ and introduce the open subsets of $X$
$$
\plaincases{
X(L,h,q)=\big\{x\in X\,;\;\Theta_{L,h}(x)~\hbox{has signature $(n-q,q)$}\big\},
\cr
\noalign{\vskip5pt}
\displaystyle
X(L,h,S)=\bigcup_{q\in S} X(L,h,q),\quad \forall S\subset\{0,1,\ldots,n\}.\cr}
\leqno(*)
$$
Then, for all $q=0,1,\ldots,n$, the
dimensions $h^q(X,E\otimes L^m)$ of cohomology groups of the tensor powers 
$E\otimes L^m$ satisfy the following ``Strong Morse inequalities''
as $m\to +\infty\,:$
$$\sum_{0\le j\le q} (-1)^{q-j}h^j(X,E\otimes L^m) \le r {m^n\over n!}
\int_{X(L,h,\le q)}(-1)^q\Theta_{L,h}^n+o(m^n),
\leqno\SM(q):$$
with equality
$\chi(X,E\otimes L^m)= r{m^n\over n!}\int_X \Theta_{L,h}^n + o(m^n)$
for the Euler characteristic $(q=n)$.
\endclaim

\noindent
As a consequence, one gets upper and lower bounds for all cohomology
groups, and especially a very useful criterion for the existence of
sections of large multiples of $L$.
\vskip2mm

\claim 2.2. Corollary|Under the above hypotheses, we have
\vskip2pt
\plainitem{\rm(a)} Upper bound for $h^q$ $($Weak Morse inequalities$)\,:$
$$h^q(X,E\otimes L^m)\le r {m^n\over n!}\int_{X(L,h,q)} (-1)^q \Theta_{L,h}^n + o(m^n)~.$$
\vskip2pt
\plainitem{\rm(b)} Lower bound for $h^0\,:$
$$
h^0(X,E\otimes L^m)\ge h^0-h^1\ge
 r{m^n\over n!}\int_{X(L,h,\le 1)}\Theta_{L,h}^n -o(m^n)~.$$
Especially $L$ is big as soon as $\int_{X(L,h,\le 1)}\Theta_{L,h}^n>0$
for some hermitian metric $h$ on~$L$.
\vskip2pt
\plainitem{\rm(c)} Lower bound for $h^q\,:$
$$
h^q(X,E\otimes L^m)\ge h^q-h^{q-1}-h^{q+1}\ge
r{m^n\over n!}\int_{X(L,h,\{q,q\pm 1\})}
(-1)^q \Theta_{L,h}^n + o(m^n)~.$$
\endclaim

\plainproof. (a) is obtained by taking $\SM(q)+\SM(q\,{-}\,1)$, (b) is equivalent to
$-\SM(1)$ and (c) is equivalent to $-(\SM(q\,{+}\,1)+\SM(q\,{-}\,2))$.\qed

\noindent
The following simple lemma is the key to derive algebraic Morse
inequalities from their analytic form (cf.\ [Dem94], Theorem~12.3).

\claim 2.3.~Lemma|Let $\eta=\alpha-\beta$ be a difference of semipositive 
$(1,1)$-forms on an $n$-dimensional complex manifold~$X$, 
and let $\bOne_{\eta,\le q}$ be the characteristic function of the
open set where $\eta$ is non degenerate with a number of negative eigenvalues 
at most equal to~$q$.
Then
$$
(-1)^q\bOne_{\eta,\le q}~\eta^n\le \sum_{0\le j\le q}(-1)^{q-j}
{n\choose j}\alpha^{n-j}\wedge\beta^j,
$$
in particular
$$
\bOne_{\eta,\le 1}~\eta^n\ge \alpha^n-n\alpha^{n-1}\wedge \beta\qquad\hbox{for $q=1$.}
$$
\endclaim

\plainproof. Without loss of generality, we can assume $\alpha>0$ positive definite, so that
$\alpha$ can be taken as the base hermitian metric on~$X$. Let us denote by
$$
\lambda_1\ge\lambda_2\ge\ldots\ge\lambda_n\ge 0
$$
the eigenvalues of $\beta$ with respect to $\alpha$. The eigenvalues of $\eta=\alpha-\beta$
are then given by 
$$
1-\lambda_1\le\ldots\le 1-\lambda_q\le 1-\lambda_{q+1}\le\ldots\le 1-\lambda_n,
$$
hence the open set $\{\lambda_{q+1}<1\}$ coincides with the support of 
$\bOne_{\eta,\le q}$, except that it may also contain a part of 
the degeneration set $\eta^n=0$. On the other hand we have
$${n\choose j}\alpha^{n-j}\wedge\beta^j=\sigma_n^j(\lambda)\,\alpha^n,$$
where $\sigma_n^j(\lambda)$ is the $j$-th elementary symmetric function in the $\lambda_j$'s.
Thus, to prove the lemma, we only have to check that
$$\sum_{0\le j\le q}(-1)^{q-j}\sigma_n^j(\lambda)-
\bOne_{\{\lambda_{q+1}<1\}}(-1)^q\prod_{1\le j\le n}(1-\lambda_j)\ge 0.$$
This is easily done by induction on~$n$ (just split apart the parameter
$\lambda_n$ and write $\sigma_n^j(\lambda)=
\sigma_{n-1}^j(\lambda)+\sigma_{n-1}^{j-1}(\lambda)\,\lambda_n$).\qed

\claim 2.4.~Corollary|Assume that $\eta=\Theta_{L,h}$ can be expressed as
a difference $\eta=\alpha-\beta$ of smooth $(1,1)$-forms $\alpha,\beta\ge 0$.
Then we have
$$\sum_{0\le j\le q} (-1)^{q-j}h^j(X,E\otimes L^m) \le r {m^n\over n!}
\int_X\sum_{0\le j\le q}(-1)^{q-j}{n\choose j}\alpha^{n-j}\wedge\beta^j+o(m^n),
\leqno\SM(q):$$
and in particular, for $q=1$,
$$h^0(X,E\otimes L^m)\ge h^0-h^1\ge
r{m^n\over n!}\int_X\alpha^n-n\alpha^{n-1}\wedge\beta+o(m^n).
$$
\endclaim

\claim 2.5. Remark|{\rm These estimates are consequences of Theorem~2.1
and Lemma 2.3, by taking the integral over $X$. The estimate for $h^0$
was stated and studied by Trapani [Tra93]. In the special case
$\alpha=\Theta_{A,h_A}>0$, $\beta=\Theta_{B,h_B}>0$ where $A,B$ are ample
line bundles, a direct proof can be obtained by purely algebraic means,
via the Riemann-Roch formula. However, we will later have to use
Corollary 2.4 in case $\alpha$ and $\beta$ are not closed, a situation in
which no algebraic proof seems to exist.}
\endclaim

\plainsubsection 2.B. Singular holomorphic Morse inequalities|

The case of singular hermitian metrics has been considered
in Bonavero's PhD thesis {\rm [Bon93]} and will be important for~us.
We assume that $L$ is equipped with a singular hermitian metric
$h=h_\infty e^{-\varphi}$ with analytic singularities, i.e.,
$h_\infty$ is a smooth metric, and on an neighborhood
$V\ni x_0$ of an arbitrary point $x_0\in X$, the weight $\varphi$
is of the form $$
\varphi(z)=c\log\sum_{1\le j\le N}|g_j|^2+u(z)
\leqno(2.6)
$$
where $g_j\in\cO_X(V)$ and $u\in C^\infty(V)$. We then have
$\Theta_{L,h}=\alpha+{\ii\over 2\pi}\ddbar\varphi$ where
$\alpha=\Theta_{L,h_\infty}$ is a smooth closed $(1,1)$-form on~$X$.
In this situation, the multiplier ideal sheaves 
$$
\cI(h^m)=\cI(k\varphi)=\big\{f\in\cO_{X,x},\;\;\exists V\ni x,~
\int_V|f(z)|^2e^{-m\varphi(z)}d\lambda(z)<+\infty\big\}\leqno(2.7)
$$
play an important role. We define the singularity set of $h$ by
$\Sing(h)=\Sing(\varphi)=\varphi^{-1}(-\infty)$ which,
by definition, is an analytic subset of $X$. The associated $q$-index sets are
$$
X(L,h,q)=\big\{x\in X\ssm\Sing(h)\,;\;
\Theta_{L,h}(x)~\hbox{has signature $(n-q,q)$}\big\}.
\leqno(2.8)
$$
We can then state:

\claim 2.9. Theorem {\rm([Bon93])}|Morse inequalities still hold in
the context of singular hermitian metric with analytic singularities,
provided the cohomology groups under consideration are twisted by
the appropriate multiplier ideal sheaves, i.e.\ replaced by
$H^q(X,E\otimes L^m\otimes\cI(h^m))$.
\endclaim

\claim 2.10. Remark|{\rm The assumption (2.6) guarantees that the measure
$\bOne_{X\ssm\Sing(h)}(\Theta_{L,h})^n$ is locally integrable on~$X$,
as is easily seen by using the Hironaka desingularization theorem and by
taking a log resolution $\mu:\wt X\to X$ such that $\mu^*(g_j)=(\gamma)\subset
\cO_{\smash{\wt X}}$ becomes a
principal ideal associated with a simple normal crossing divisor
$E=\div(\gamma)$. Then $\mu^*\Theta_{L,h}=
c[E]+\beta$ where $\beta$ is a smooth closed $(1,1)$-form on $\wt X$, hence
$$
\mu^*(\bOne_{X\ssm\Sing(h)}\Theta_{L,h}^n)=\beta^n~~\Rightarrow~~
\int_{X\ssm\Sing(h)}\Theta_{L,h}^n=\int_{\wt X}\beta^n.
$$
It should be observed that the multiplier ideal sheaves $\cI(h^m)$ and the
integral $\int_{X\ssm\Sing(h)}\Theta_{L,h}^n$ only depend on the equivalence
class of singularities of $h\,$: if we have two metrics with analytic
singularities $h_j=h_\infty e^{-\varphi_j}$, $j=1,2$, such that
$\psi=\varphi_2-\varphi_1$ is bounded, then, with the above notation,
we have $\mu^*\Theta_{L,h_j}=c[E]+\beta_j$ and
$\beta_2=\beta_1+{\ii\over 2\pi}\ddbar\psi$, therefore
$\int_{\wt X}\beta_2^n=\int_{\wt X}\beta_1^n$ by Stokes theorem. By using
Monge-Amp\`ere operators in the sense of Bedford-Taylor [BT76], it is in
fact enough to assume $u\in L^\infty_\loc(X)$ in (2.6), and
$\psi\in L^\infty(X)$ here. In general,
however, the Morse integrals $\int_{X(L,h_j,q)}(-1)^q\Theta_{L,h_j}^n$,
$j=1,2$, will~differ.}
\endclaim

\plainsubsection 2.C. Morse inequalities and semi-continuity|
Let $\cX\to S$ be a proper and flat morphism of reduced complex spaces,
and let $(X_t)_{t\in S}$ be the fibers. 
Given a sheaf $\cE$ over $\cX$ of locally free $\cO_\cX$-modules of rank $r$,
inducing on the fibres a family of sheaves $(E_t\to X_t)_{t\in S}$,
the following semicontinuity property holds ([CRAS]):

\claim 2.11. Proposition|For every $q\ge 0$,
the alternate sum 
$$
t\mapsto h^q (X_t,E_t)-h^{ q-1} (X_t,E_t)+. . .+(-1)^q h^0 (X_t,E_t)
$$
is upper semicontinuous with respect to the (analytic) Zariski topology
on~$S$.
\endclaim

Now, if  $\cL\to\cX$ is an invertible sheaf equipped with a smooth
hermitian metric $h$, and if $(h_t)$ are the fiberwise metrics on the
family $(L_t\to X_t)_{t\in S}$, we get
$$
\sum_{j=0}^q(-1)^{q-j}h^j(X_t,E_t\otimes L_t^{\otimes m})
\le
r{m^n\over n!}\int_{X(L_0,h_0,\le q)}
(-1)^q\Theta_{L_0,h_0}^n + \delta(t)m^n,
\leqno(2.12)
$$
where $\delta(t)\to 0$ as $t\to 0$.
In fact, the proof of holomorphic Morse inequalities shows that the
inequality holds uniformly on every relatively compact $S'\compact S$, with
$$
I(t)=\int_{X(L_t,h_t,\le q)}(-1)^q\Theta_{L_t,h_t}^n=
\int_X (-1)^q\bOne_{X(L_t,h_t,\le q)}\Theta_{L_t,h_t}^n
$$
in the right hand side, and $t\mapsto I(t)$ is clearly continuous with
respect to the ordinary topology. In other words, the Morse integral
computed on the central fibers
provides uniform upper bounds for cohomology groups of $E_t\otimes
L_t^{\otimes m}$ when $t$ is close to $0$ in ordinary topology
(and also, as a consequence, for $t$ in a complement
$S\ssm \bigcup S_m$ of at most countably many analytic strata
$S_m\subsetneq S$).

\claim 2.13. Remark|{\rm Similar results would hold when $h$ is a singular
hermitian metric with analytic singularities on $\cL\to\cX$, under
the restriction that the families of multiplier ideal sheaves
$(\cI(h_t^m))_{t\in S}$ ``never jump''.}
\endclaim

\plainsubsection 2.D. Case of filtered bundles|

Let $E\to X$ be a vector bundle over a variety, equipped with a filtration
(or multifiltration) $F^p(E)$, and let $G=\bigoplus G^p(E)\to X$
be the graded bundle associated to this filtration.

\claim 2.14. Lemma|In the above setting, one has for every $q\ge 0$
$$
\sum_{j=0}^q(-1)^{q-j}h^j(X,E)\le\sum_{j=0}^q(-1)^{q-j}h^j(X,G).
$$
\endclaim

\plainproof. One possible argument is to use the well known fact that
there is a family of filtered bundles $(E_t\to X)_{t\in \bC}$
(with the same graded pieces $G^p(E_t)=G^p(E)$), such
that $E_t\simeq E$ for all $t\neq 0$ and $E_0\simeq G$. The result is then
an immediate consequence of the semi-continuity result~2.11. A more
direct very elementary argument can be given as follows: by transitivity
of inequalities, it is sufficient to prove the result for simple filtrations;
then, by induction on the length of filtrations, it is sufficient to
prove the result for exact sequences $0\to S\to E\to Q\to 0$ of vector
bundles on $X$. Consider the associated (truncated) long exact sequence
in cohomology:
$$
\eqalign{  
0\to H^0(S)\to H^0(E)\to H^0(Q)&\build\to|\delta_1||\cdots\cr
&\build\to|\delta_{q-1}|| H^q(S)\to H^q(E)\to H^q(Q)\build\to|\delta_q||
\Im(\delta_q)\to 0.\cr}
$$
By the rank theorem of linear algebra,
$$
0\le\rank(\delta_q) = (-1)^q\sum_{j=0}^q(-1)^j(h^j(X,Q)- h^j(X,E)+ h^j(X,S)).
$$
The result follows, since here $h^j(X,G)=h^j(X,Q)+h^j(X,S)$.
\qed

\plainsubsection 2.E. Rees deformation construction (after Cadorel)|

In this short paragraph, we outline a nice algebraic interpretation by
Beno\^it Cadorel of certain semi-continuity arguments for cohomology
group dimensions that underline the analytic approach of [Dem11, Lemma~2.12
and Prop.~2.13] and [Dem12, Prop.~9.28] (we will anyway explain again
its essential points in \S$\,$3, since we have to deal here with a more
general situation). Recall after [Cad17, Prop.~4.2, Prop.~4.5], that
the Rees deformation construction allows one to construct natural 
deformations of Green-Griffiths jets spaces to weighted projectivized bundles.

Let $(X,V,D)$ be a non singular directed orbifold, and let
$g:Y\to(X,D)$ be an adapted Galois cover, as briefly
described in remark~1.26, see also [CDR18, \S$\,$2.1] for more details. We
then get a Green-Griffiths jet bundle of graded algebras
$E_{k,\bullet }\wt V^\star\langle D \rangle\to Y$ which admits a
multifiltration of associated graded algebra
$$
G^\bu E_{k,\bu}\widetilde V^*\langle D\rangle=\bigoplus_{m\in\bN}
\bigoplus_{\ell_1+2\ell_2+\cdots+k\ell_k=m}
S^{\ell_1}\wt V^*\langle D^{(1)}\rangle\otimes \cdots
\otimes S^{\ell_k}\wt V^*\langle D^{(k)}\rangle.
$$
where the tilde means taking  pull-backs by $g^*$.
Applying the Proj functor, one gets a weighted projective bundle:
$$
\bP_{(1,\cdots,k)}\left(
\wt V^*\langle D^{(1)})\rangle\oplus \cdots \oplus
\wt V^*\langle D^{(k)}\rangle \right) =
\Proj\left(G^\bu E_{k,\bu}\widetilde V^*\langle D\rangle\right)
\build\to|\rho_k|| Y,
$$
Then, following mutadis mutandus the arguments of Cadorel, one constructs 
a family $Y\build\leftarrow|p_k||\cY_k\to \bC$ parametrized by $\bC$, with a canonical line bundle $\cO_{\cY_k}(1)$ such that:

\plainitem{$\bu$}
  the central fiber $\cY_{k,0}$ is $\bP_{(1,\cdots,k)}
  \left(
\wt V^*\langle D^{(1)})\rangle\oplus \cdots \oplus
\wt V^*\langle D^{(k)}\rangle \right) $
and the restriction of $\cO_{\cY_k}(1)$ coincide with the canonical
line bundle of this weighted projective bundle.
Hence $(\pi_k)_*\cO_{\cY_{k,0}}(m)= G^\bu E_{k,m}\wt V^*\langle D\rangle$.

\plainitem{$\bu$} the other fibers $\cY_{k,t}$ are isomorphic to the
singular quotient $J^k(Y,\wt V,D)/\bC^*$ for the natural
$\bC^*$-action by homotheties, where $J^k(Y,\wt V,D)$ is the
affine algebraic bundle associated with the sheaf of algebras,
and $(\pi_k)_*\cO_{\cY_{k,t}}(m)\simeq E_{k,m}\wt V^*\langle D\rangle$.
\medskip

\noindent
Applying the semicontinuity result of [Dem95], and working with
holomorphic inequalities, we obtain a control about dimensions of
cohomology spaces of $E_{k,m}\wt V^*\langle D\rangle$ in terms
of dimensions of cohomology spaces of the a priori simpler
graded pieces $G^\bu E_{k,m}\wt V^*\langle D\rangle$.
This reduces the study of higher order jet differentials to
sections of the tautological sheaves on the weighted
projective space associated with a direct sum combination of
symmetric differentials. In particular, we have

\claim 2.15. Lemma|For every $q\in\bN$
$$
\sum_{j=0}^q(-1)^{q-j}h^j(Y,E_{k,m}\wt V^*\langle D\rangle)
\ge
\sum_{j=0}^q(-1)^{q-j}h^j(Y,G^\bu E_{k,m}\wt V^*\langle D\rangle).
$$
Especially, for $q=1$, we get
$$
\eqalign{
h^0(Y,E_{k,m}\wt V^*\langle D\rangle)
&\ge
h^0(Y,E_{k,m}\wt V^*\langle D\rangle)-
h^1(Y,E_{k,m}\wt V^*\langle D\rangle)\cr
&\ge
h^0(Y,G^\bu E_{k,m}\wt V^*\langle D\rangle)-
h^1(Y,G^\bu E_{k,m}\wt V^*\langle D\rangle).\cr}
$$
\endclaim

\plainsection{3. Construction of jet metrics and orbifold jet metrics}

\plainsubsection 3.A. Jet metrics and curvature tensor of jet bundles|

Let $(X,V)$ be a non singular directed variety and $h$ a hermitian metric
on $V$. We assume that $h$ is smooth at this point (but will later relax
a little bit this assumption and allow certain singularities).
Near any given point $z_0\in X$, we can choose local coordinates
$z=(z_1,\ldots,z_n)$ centered at $z_0$ and a local holomorphic coordinate frame
$(e_\lambda(z))_{1\le\lambda\le r}$ of $V$ on an open set $U\ni z_0$, 
such that
$$
\langle e_\lambda(z),e_\mu(z)\rangle_{h(z)} =\delta_{\lambda\mu}+
\sum_{1\le i,j\le n,\,1\le\lambda,\mu\le r}c_{ij\lambda\mu}z_i\overline z_j+
O(|z|^3)\leqno(3.1)
$$
for suitable complex coefficients $(c_{ij\lambda\mu})$. It is a standard fact
that such a normalized coordinate system always exists, and that the 
Chern curvature tensor ${\ii\over 2\pi}\nabla^2_{V,h}$ of $(V,h)$ at $z_0$ 
is given by
$$
\Theta_{V,h}(z_0)=-{\ii\over 2\pi}
\sum_{i,j,\lambda,\mu}
c_{ij\lambda\mu}\,dz_i\wedge d\overline z_j\otimes e_\lambda^*\otimes e_\mu.
\leqno(3.2)
$$
Therefore, $({\ii\over 2\pi}\,c_{ij\lambda\mu})$ are the coefficients of
$-\Theta_{V,h}$.  Up to taking the transposed tensor with respect
to $\lambda,\mu$, these
coefficients are also the components of the curvature tensor
$\Theta_{V^*,h^*}=-{}^t\Theta_{V,h}$ of the dual bundle $(V^*,h^*)$.
By (1.5), the connection $\nabla=\nabla_h$ yields a $C^\infty$
isomorphism $J_kV\to V^{\oplus k}$. Let us fix an integer $b\in\bN^*$ that 
is a multiple of $\lcm(1,2,\ldots,k)$, and positive numbers
$1=\varepsilon_1\gg\varepsilon_2\gg\cdots\gg \varepsilon_k>0$.
Following [Dem11], we define a global weighted Finsler metric
on $J^kV$ by putting for any $k$-jet $f\in J^kV_z$
$$
\Psi_{h,b,\varepsilon}(f):=\Bigg(
\sum_{1\le s\le k}\varepsilon_s^{2b}\Vert\nabla^s f(0)
\Vert_{h(z)}^{2b/s}\Bigg)^{1/b},
\leqno(3.3)
$$
where $\Vert~~\Vert_{h(z)}$ is the hermitian metric $h$ of $V$ evaluated
on the fiber $V_z$, $z=f(0)$. The function $\Psi_{h,b,\varepsilon}$ satisfies
the fundamental homogeneity property 
$$
\Psi_{h,b,\varepsilon}(\alpha\cdot f)=|\alpha|^2\,\Psi_{h,b,\varepsilon}(f)
\leqno(3.4)
$$
with respect to the $\bC^*$ action on $J^kV$, in other words, it induces
a hermitian metric on the dual $L_k^*$ of the tautological $\bQ$-line bundle
$L_k=\cO_{X_k(V)}(1)$ over $X_k(V)$. The curvature of $L_k$ is given by
$$
\pi_k^*\Theta_{L_k,\Psi^*_{h,b,\varepsilon}}={\ii\over 2\pi}\ddbar
\log\Psi_{h,b,\varepsilon}
\leqno(3.5)
$$
Our next goal is to compute precisely the curvature and to apply
holomorphic Morse inequalities to $L\to X_k(V)$ with the above metric.
This might look a priori like an untractable problem, since the definition of
$\Psi_{h,b,\varepsilon}$ is a rather complicated one, involving the hermitian
metric in an intricate manner. However, the ``miracle''
is that the asymptotic behavior of $\Psi_{h,b,\varepsilon}$ as
$\varepsilon_s/\varepsilon_{s-1}\to 0$ is in some sense uniquely defined,
and ``splits'' according to the natural multifiltration on jet differentials
(as already hinted in \S$\,$2.E). This leads to a computable asymptotic
formula, which is moreover simple enough to produce useful results.

\claim 3.6. Lemma|Let us consider the global $C^\infty$ bundle isomorphism
$J^kV\to V^{\oplus k}$ associated with an arbitrary global $C^\infty$ connection
$\nabla$ on $V\to X$, and let us introduce the rescaling transformation 
$$\rho_{\nabla,\varepsilon}(\xi_1,\xi_2,\ldots,\xi_k)=
(\varepsilon_1^1\xi_1,\varepsilon_2^2\xi_2,\ldots,
\varepsilon_k^k\xi_k)\quad
\hbox{on fibers $J^kV_z$, $z\in X$}.
$$
Such a rescaling commutes with the $\bC^*$-action. Moreover, if $p$ is 
a multiple of $\lcm(1,2,\ldots,k)$ and the ratios
$\varepsilon_s/\varepsilon_{s-1}$ tend to~$0$ for all $s=2,\ldots,k$, the
rescaled Finsler metric
\hbox{$\Psi_{h,b,\varepsilon}\circ\rho_{\nabla,\varepsilon}^{-1}
(\xi_1,\ldots,\xi_k)$} converges towards the limit
$$
\bigg(\sum_{1\le s\le k}\Vert \xi_s\Vert^{2b/s}_h\bigg)^{1/b}
$$
on every compact subset of $V^{\oplus k}\ssm\{0\}$,
uniformly in $C^\infty$ topology, and the limit is independent
of the connection~$\nabla$. The error is measured by a multiplicative factor
$1\pm O(\max_{2\le s\le k}(\varepsilon_s/\varepsilon_{s-1})^s)$.
\endclaim

\plainproof. Let us pick another $C^\infty$ connection $\wt\nabla=
\nabla+\Gamma$ where $\Gamma\in C^\infty(U,T^*_X\otimes
\Hom(V,V))$. Then $\wt\nabla^2f=\nabla^2f+\Gamma(f)(f')\cdot f'$, and
inductively we get
$$
\wt\nabla^sf=\nabla^sf+P_s(f\,;\,\nabla^1f,\ldots,\nabla^{s-1}f)
$$
where $P(z\,;\,\xi_1,\ldots,\xi_{s-1})$ is a polynomial with $C^\infty$
coefficients in $z\in U$, which is of weighted homogeneous degree
$s$ in $(\xi_1,\ldots,\xi_{s-1})$. In other words, the corresponding 
isomorphisms  $J^kV\simeq V^{\oplus k}$ correspond to each other
by a $\bC^*$-homogeneous transformation $(\xi_1,\ldots,\xi_k)\mapsto
(\wt\xi_1,\ldots,\wt\xi_k)$ such that
$$
\wt\xi_s=\xi_s+P_s(z\,;\,\xi_1,\ldots,\xi_{s-1}).
$$
Let us introduce the corresponding rescaled components
$$
(\xi_{1,\varepsilon},\ldots,\xi_{k,\varepsilon})=
(\varepsilon_1^1\xi_1,\ldots,\varepsilon_k^k\xi_k),\qquad
(\wt\xi_{1,\varepsilon},\ldots,\wt\xi_{k,\varepsilon})=
(\varepsilon_1^1\wt\xi_1,\ldots,\varepsilon_k^k\wt\xi_k).
$$
Then
$$
\eqalign{
\wt\xi_{s,\varepsilon}
&=\xi_{s,\varepsilon}+
\varepsilon_s^s\,P_s(x\,;\,\varepsilon_1^{-1}\xi_{1,\varepsilon},\ldots,
\varepsilon_{s-1}^{-(s-1)}\xi_{s-1,\varepsilon})\cr
&=\xi_{s,\varepsilon}+O(\varepsilon_s/\varepsilon_{s-1})^s\,
O(\Vert\xi_{1,\varepsilon}\Vert+\cdots+\Vert\xi_{s-1,\varepsilon}
\Vert^{1/(s-1)})^s\cr}
$$
and it is easily seen, as a simple consequence of the mean value inequality
$|\Vert x\Vert^\gamma-\Vert y\Vert^\gamma|\le\gamma\sup_{z\in[x,y]}
\Vert z\Vert^{\gamma-1}\Vert x-y\Vert$, that
the ``error term'' in the difference
$\Vert\wt\xi_{s,\varepsilon}\Vert^{2b/s}-\Vert\xi_{s,\varepsilon}\Vert^{2b/s}$
is bounded by
$$
(\varepsilon_s/\varepsilon_{s-1})^s\,
\big(\Vert\xi_{1,\varepsilon}\Vert+\cdots+
\Vert\xi_{s-1,\varepsilon}\Vert^{1/(s-1)}+
\Vert\xi_{s,\varepsilon}\Vert^{1/s}\big)^{2b}.
$$
When $b/s$ is an integer, similar bounds hold for all
derivatives $D_{z,\xi}^\beta(\Vert\wt\xi_{s,\varepsilon}\Vert^{2b/s}-
\Vert\xi_{s,\varepsilon}\Vert^{2b/s})$ and the lemma follows.\qed

Now, we fix a point $z_0\in X$, a local holomorphic frame 
$(e_\lambda(z))_{1\le\lambda\le r}$ satisfying (3.1) on a neighborhood $U$ 
of~$z_0$, and the {\it holomorphic} connection $\nabla$ on $V_{|U}$ such that
$\nabla e_\lambda=0$. Since the uniform estimates of Lemma~3.6 also apply
locally (provided they are applied on a relatively compact open
subset $U'\compact U$), we can use the corresponding holomorphic
trivialization $J^kV_{|U}\simeq V_{|U}^{\oplus k}\simeq U\times(\bC^r)^{\oplus k}$
to make our calculations. We do this in terms of the rescaled components 
$\xi_s=\varepsilon_s^s\nabla^sf(0)$. Then, uniformly on compact subsets
of $J^kV_{|U}\ssm\{0\}$, we have
$$
\Psi_{h,b,\varepsilon}\circ\rho_{\nabla,\varepsilon}^{-1}(z\,;\,\xi_1,\ldots,\xi_k)
=\bigg(\sum_{1\le s\le k}\Vert\xi_s\Vert^{2b/s}_{h(z)}\bigg)^{1/b}
+O(\max((\varepsilon_s/\varepsilon_{s-1})^{1/b}),
$$
and the error term remains of the same magnitude when we take
any derivative $D_{z,\xi}^\beta$. By (3.1) we find
$$
\Vert \xi_s\Vert_{h(z)}^2=
\sum_\lambda|\xi_{s,\lambda}|^2+
\sum_{i,j,\lambda,\mu}c_{ij\lambda\mu}\,z_i\overline z_j
\,\xi_{s,\lambda}\overline \xi_{s,\mu}+O(|z|^3|\xi|^2).
$$
The question is thus reduced to evaluating the curvature of the weighted
Finsler metric on $V^{\oplus k}$ defined by
$$
\eqalign{
\Psi(z\,;\,\xi_1,\ldots,\xi_k)
&=\bigg(\sum_{1\le s\le k}\Vert\xi_s\Vert^{2b/s}_{h(z)}\bigg)^{1/b}\cr
&=\bigg(\sum_{1\le s\le k}\Big(\sum_\lambda|\xi_{s,\lambda}|^2+
\sum_{i,j,\lambda,\mu}c_{ij\lambda\mu}\,z_i\overline z_j\,
\xi_{s,\lambda}\overline\xi_{s,\mu}\Big)^{b/s}\bigg)^{1/b}+O(|z|^3).\cr}
$$
We set $|\xi_s|^2=\sum_\lambda|\xi_{s,\lambda}|^2$. A straightforward 
calculation yields the Taylor expansion
$$
\eqalign{
&\log\Psi(z\,;\,\xi_1,\ldots,\xi_k)\cr
&~~{}={1\over b}\log\sum_{1\le s\le k}|\xi_s|^{2b/s}+
\sum_{1\le s\le k}{1\over s}\,{|\xi_s|^{2b/s}\over \sum_t|\xi_t|^{2b/t}}
\sum_{i,j,\lambda,\mu}c_{ij\lambda\mu}z_i\overline z_j
{\xi_{s,\lambda}\overline\xi_{s,\mu}\over|\xi_s|^2}+O(|z|^3).\cr}
$$
By (3.5), the curvature form of $L_k=\cO_{X_k(V)}(1)$ 
is given at the central point $z_0$ by the formula
$$
\Theta_{L_k,\Psi^*_{h,b,\varepsilon}}(z_0,[\xi])\simeq
\omega_{r,k,b}(\xi)+{\ii\over 2\pi}
\sum_{1\le s\le k}{1\over s}\,{|\xi_s|^{2b/s}\over \sum_t|\xi_t|^{2b/t}}
\sum_{i,j,\lambda,\mu}c_{ij\lambda\mu}
{\xi_{s,\lambda}\overline\xi_{s,\mu}\over|\xi_s|^2}\,dz_i\wedge d\overline z_j
\leqno(3.7)
$$
where $[\xi]=[\xi_1,\ldots,\xi_k ]\in\bP(1^{[r]},2^{[r]},\ldots,k^{[r]})$ and
$\omega_{r,k,b}(\xi)={\ii\over 2\pi}\ddbar({1\over b}\log\sum_{1\le s\le k}
|\xi_s|^{2b/s})$. The fibers $\bP(1^{[r]},2^{[r]},\ldots,k^{[r]})$ of
$X_k(V)\to X$ can be represented as a quotient of the
``weighted ellipsoid'' $\sum_{s=1}^k|\xi_s|^{2b/s}=1$ by the $\bS^1$-action
induced by the weighted $\bC^*$-action. This suggests to make use of
polar coordinates and to set
$$
\leqalignno{
  &x_s=|\xi_s|^{2b/s},\quad x=(x_1,\ldots,x_k)\in\bR^k,&(3.8)\cr
  &u_s={\xi_s\over |\xi_s|}\in \bS^{2r-1}\subset\bC^r,\quad
  u=(u_1,\ldots,u_k)\in(\bS^{2r-1})^k,&(3.8')\cr }
$$
so that
$$
\sum_{s=1}^kx_s=1\quad\hbox{and}\quad \xi_s=x_s^{s/2b}u_s.\kern102pt
\leqno(3.8'')
$$
The Morse integrals will then have to be computed for
$(x,u)\in\bDelta^{k-1}\times(\bS^{2r-1})^k$, where
$\bDelta^{k-1}\subset\bR^k$ is the $(k-1)$-dimensional simplex.

\claim 3.9. Proposition| With respect to the rescaled components
$\xi_s=\varepsilon_s^s\nabla^sf(0)$ at $z=f(0)\in X$ and the above
choice of coordinates $(3.8^*)$, the curvature of the tautological
sheaf $L_k=\cO_{X_k(V)}(1)$ admits an approximate expression
$$
\leqno\displaystyle{\rm(a)}\quad
\Theta_{L_k,\Psi^*_{h,b,\varepsilon}}(z,[\xi])=
\omega_{r,k,b}(\xi)+g_{V,k}(z,x,u)+\hbox{\rm(error terms)},
$$
where $(x,u)\in\bDelta^{k-1}\times (\bS^{2r-1})^k$, $\xi_s=x_s^{s/2b}u_s\in\bC^r$,
$$
\leqno\displaystyle{\rm (b)}\quad
\omega_{r,k,b}(\xi)={\ii\over 2\pi}\ddbar\bigg(
{1\over b}\sum_{1\le s\le k}|\xi_s|^{2b/s}\bigg)
$$
is a Fubini-Study type K\"ahler metric on
$\bP(1^{[r]},2^{[r]},\ldots,k^{[r]})$, associated with the canonical
$\bC^*$ action on $J^kV$ of weight $a=(1^{[r]},2^{[r]},\ldots,k^{[r]})$, and
$$
\leqno\displaystyle{\rm (c)}\quad
g_{V,k}(z,x,u)={\ii\over 2\pi}\sum_{1\le s\le k}{x_s\over s}
\sum_{i,j,\lambda,\mu}c_{ij\lambda\mu}(z)\,
u_{s,\lambda}\overline u_{s,\mu}\,dz_i\wedge d\overline z_j.
$$
Here $({\ii\over 2\pi}\,c_{ij\lambda\mu})$ are the coefficients
of $-\Theta_{V,h}$, and the error terms admit an upper bound
$$
\leqno\displaystyle{\rm (d)}\quad
\hbox{\rm(error terms)}\le
O\Big(\max_{2\le s\le k}(\varepsilon_s/\varepsilon_{s-1})^s\Big)\quad
\hbox{uniformly on the compact variety~$X_k(V)$}.
$$
\endclaim

\plainproof. The error terms on $\Theta_{L_k}$ come from the differentiation of
the error terms on the Finsler metric, found in Lemma 3.6. They can indeed
be differentiated if $b$ is a multiple of $\lcm(1,2,\ldots,k)$, since
$2b/s$ is then an even integer.\qed

\noindent
For the calculation of Morse integrals, it is useful to find the expression
of the volume form $\omega_{r,k,b}^{kr-1}$ on
$\bP(1^{[r]},2^{[r]},\ldots,k^{[r]}) = (\bDelta^{k-1}\times (\bS^{2r-1})^k)/\bS^1$
in terms of the coordinates $(x,u)$. We~refer to [Dem11, Prop.~1.13]
for the proof.

\claim 3.10. Proposition|
\plainitem{\rm(a)} The volume form $\omega_{r,k,b}^{kr-1}$ is the quotient of the
measure ${1\over k!^r}\nu_{k,r}\otimes\mu$ on $\bDelta^{k-1}\times(\bS^{2r-1})^k$,
where  
$$
d\nu_{k,r}(x)=(kr-1)!{(x_1\ldots\,x_k) ^{r-1}\over (r-1)!\,{}^k}
dx_1\wedge\ldots\wedge dx_{k-1},\quad
d\mu(u)=d\mu_1(u_1)\ldots d\mu_k(u_k)
$$
are probability measures on $\bDelta^{k-1}$ and $(\bS^{2r-1})^k$ respectively
$(\mu$ being the rotation inva\-riant one$)$.
\vskip2pt
\plainitem{\rm(b)} We have the equality~
$\displaystyle
\int_{\bP(1^{[r]},2^{[r]},\ldots,k^{[r]})}\omega_{r,k,b}^{kr-1}={1\over k!^r}$~
$($independent of~$b)$.
\endclaim


\plainsubsection 3.B. Logarithmic and orbifold jet metrics|

Consider now an arbifold directed structure $(X,V,D)$, where
$V\subset T_X$ is a subbundle, $r=\rank(V)$, and
$D=\sum(1-{1\over \rho_j})\Delta_j$ is
a normal crossing divisor that is assumed to intersect $V$ transversally
everywhere. One then performs very similar calculations to what we did in
\S$\,$3.A, but with adapted Finsler metrics.
Fix a point $z_0$ at which $p$ components $\Delta_j$ meet, and use coordinates
$(z_1,\ldots,z_n)$ such that $V_{z_0}$ is spanned by
$({\partial\over\partial z_1},\ldots,{\partial\over\partial z_r})$
and $\Delta_j$ is defined by $z_j=0$, $1\le j\le p\le r$.
In the logarithmic case $\rho_j=\infty$, the logarithmic dual bundle
$\cO(V^*\langle D\rangle)$ is spanned by
$$
{dz_1\over z_1},\ldots,{dz_p\over z_p},~dz_{p+1},\ldots,dz_n.
$$
The logarithmic jet differentials are just polynomials in
$$
{d^sz_1\over z_1},\ldots,{d^sz_p\over z_p},~d^sz_{p+1},\ldots,d^sz_n,\quad
1\le s\le k,
$$
and the corresponding $(\varepsilon_1,\ldots,\varepsilon_k)$-rescaled
Finsler metric is
$$
\Bigg(\sum_{s=1}^k\varepsilon_s^{2b}\bigg(\sum_{j=1}^p
|f_j|^{-2}|f_j^{(s)}|^2+\sum_{j=p+1}^r|f_j^{(s)}|^2\bigg)^{2b/s}\Bigg)^{1/b}.
\leqno(3.11)
$$
Alternatively, we could replace $|f_j|^{-2}|f_j^{(s)}|^2$ by
$|(\log f_j)^{(s)}|^2$ which has the same leading term and differs by
a weighted degree $s$ polynomial in the $f_j^{-1}f_j^{(\ell)}$,
$\ell<s\,$; an argument very similar to the one used in the
proof of Lemma 3.6 then shows that the difference is negligible
when $\varepsilon_1\gg \varepsilon_2\gg \cdots\gg\varepsilon_k$.
However (3.11) is just the case of the model metric, in fact we get
$r$-tuples $\xi_s=(\xi_{s,j})_{1\le j\le r}$ of components produced
by the trivialization of the logarithmic bundle
$\cO(V\langle D\rangle)$, such that
$$
\xi_{s,j}=f_j^{-1}f_j^{(s)}\quad\hbox{for $1\le s\le p$ and}\quad
\xi_{s,j}=f_j^{(s)}\quad\hbox{for $p+1\le s\le r$}.\leqno(3.12)
$$
In general, we are led
to consider Finsler metrics of the form
$$
\Bigg(\sum_{s=1}^k\varepsilon_s^{2b}\Vert\xi_s\Vert_{h(z)}^{2b/s}\Bigg)^{1/b},
\quad\xi_s=(\xi_{s,j})_{1\le j\le r},
\leqno(3.13)
$$
where $h(z)$ is a variable hermitian metric on the logarithmic bundle
$V\langle D\rangle$.
In the orbifold case, the appropriate ``model'' Finsler metric is
$$
\Bigg(\sum_{s=1}^k\varepsilon_s^{2b}\bigg(\sum_{j=1}^p
|f_j|^{-2(1-s/\rho_j)_+}|f_j^{(s)}|^2+\sum_{j=p+1}^r|f_j^{(s)}|^2\bigg)^{2b/s}\Bigg)^{1/b}.
\leqno(3.14)
$$
As a consequence of Remark~2.10, we would get a metric with equivalent
singularities on the dual $L_k^*$ of the tautological sheaf
$L_k=\cO_{X_k(V\langle D\rangle)}(1)$ by replacing
$\sum_{j=p+1}^r|f_j^{(s)}|^2$ with $\sum_{j=1}^r|f_j^{(s)}|^2$ (or by
any smooth hermitian norm $h$ on $V$), since the extra terms
$\sum_{j=1}^p|f_j^{(s)}|^2$ are anyway controlled by the ``orbifold part''
of the summation. Of course, we need to find a suitable Finsler metric
that is globally defined on $X$. This can be done by taking smooth
metrics $h_{V,s}$ on $V$ and $h_j$ on $\cO_X(\Delta_j)$ respectively, as
well as smooth connections $\nabla$ and $\nabla_j$. One can then
consider the globally defined metric
$$
\Bigg(\sum_{s=1}^k\varepsilon_s^{2b}\bigg(\Vert\nabla^{(s)}f\Vert_{h_{V,s}}^2
+\sum_j\Vert\sigma_j(f)\Vert_{h_j}^{-2(1-s/\rho_j)_+}
\Vert\nabla_j^{(s)}(\sigma_j\circ f)\Vert_{h_j}^2\bigg)^{2b/s}\Bigg)^{1/b}
\leqno(3.15)
$$
where $D=\sum(1-{1\over\rho_j})\Delta_j$ and
$\sigma_j\in H^0(X,\cO_X(\Delta_j))$ are the tautological sections; here,
we want the flexibility of not necessarily taking the same hermitian metrics
on $V$ to evaluate the various norms $\Vert\nabla^{(s)}f\Vert_{h_{V,s}}$.
We obtain Finsler metrics with equivalent singularities by just changing the
$h_{V,s}$ and $h_j$ (and keeping $\nabla$, $\nabla_j$ unchanged). If
we also change the connections, then an argument very similar to the
one used in the proof of Lemma~3.6 shows that the ratio of the
corresponding metrics is
$1\pm O(\max(\varepsilon_s/\varepsilon_{s-1}))$, and therefore
arbitrary close to $1$ whenever
$\varepsilon_1\gg\varepsilon_2\gg \cdots\gg\varepsilon_k$; in~any
case, we get metrics with equivalent singularities. Fix $z_0\in X$ and
use coordinates $(z_1,\ldots,z_n)$ as described at the beginning of
\S$\,$3.B, so that $\sigma_j(z)=z_j$, $1\le j\le p$, in a suitable
trivialization of $\cO_X(\Delta_j)$. Let $f$ be a $k$-jet of curve
such that $f(0)=z\in X\ssm|D|$ is in a sufficiently small
neighborhood of $z_0$. By employing the trivial connections
associated with the above coordinates, the derivative $f^{(s)}$
is described by components
$$
\xi_{s,j}=f_j^{(s)},~~1\le j\le r,\quad
\xi^{\log}_{s,j}=f_j^{-1}f_j^{(s)},\quad
\xi^\orb_{s,j}=f_j^{-(1-s/\rho_j)_+}f_j^{(s)},\quad 1\le j\le p,
$$
and $\xi^\orb_{s,j}=\xi^{\log}_{s,j}=\xi_{s,j}$ for $p+1\le j\le r$.
Here $\xi^\orb_{s,j}$ are to be thought of as the components of $f^{(s)}$
in the ``virtual'' vector bundle $V\langle D^{(s)}\rangle$, and the fact
that the argument of these complex numbers is not uniquely defined is
irrelevant, because the only thing we need to compute the norms
is~$|\xi^\orb_{s,j}|$. Accordingly, for
$v\in V_z$, $v\simeq(v_j)_{1\le j\le r}\in\bC^r$, we
put
$$
v^{\log}_j=z_j^{-1}v_j=\sigma_j(z)^{-1}\nabla_j\sigma_j(v)\quad\hbox{and}\quad
v^\orb_j=z_j^{-(1-s/\rho_j)_+}v_j,~~1\le j\le p,
$$
and define the orbifold hermitian norm on $V\langle D^{(s)}\rangle$
associated with $h_{V,s}$ and $h_j$ by
$$
\leqalignno{
\Vert v^\orb\Vert_{\wt h_s}{\kern-3pt}^2&=\Vert v\Vert_{h_{V,s}}^2
+\sum_{j=1}^p\Vert\sigma_j(z)\Vert_{h_j}^{-2(1-s/\rho_j)_+)}
\Vert \nabla_j\sigma_j(v)\Vert_{h_j}^2&(3.16)\cr
&=\Vert v\Vert_{h_{V,s}}^2+\sum_{j=1}^p
\Vert\sigma_j(z)\Vert_{h_j}^{2(1-(1-s/\rho_j)_+)}|v_j^{\log}|^2
&(3.16')\cr
&=\Vert v\Vert_{h_{V,s}}^2+
\sum_{j=1}^p\Vert v_j^\orb\Vert_{h_j^{1-(1-s/\rho_j)_+}}^2.
&(3.16'')
\cr}
$$
With this notation, the orbifold Finsler metric (3.15) on $k$-jets
is reduced to an expression
$$
\Vert\xi^\orb\Vert_{\Psi_{h,b,\varepsilon}}^{\,2\phantom{\big|}}=
\Bigg(\sum_{s=1}^k\varepsilon_s^{2b}\Vert\xi_s^\orb
\Vert_{\wt h_s}^{2b/s}\Bigg)^{1/b},
\quad\xi_s^\orb=(\xi_{s,j}^\orb)_{1\le j\le r}\,,~~
\xi^\orb=(\xi_s^\orb)_{1\le s\le k}\,,
\leqno(3.17)
$$
formally identical to what we had in the compact or logarithmic cases. If
$v$ is a local holomorphic section of $\cO_X(V)$, formula (3.16) shows that
the norm
$\Vert v^\orb\Vert_{\wt h_s}$ can take infinite values when $z\in|D|$,
while, by $(3.16')$, the norm is always bounded (but slightly degenerate along
$|D|$) if $v$ is a section of the logarithmic sheaf
$\cO_X(V\langle\lceil D\rceil\rangle)$; we think intuitively of the
orbifold total space $V\langle D^{(s)}\rangle$ as the subspace of $V$
in which the tubular neighborhoods of the zero section are
defined by $\Vert v^\orb\Vert_{\wt h_s}<\varepsilon$ for $\varepsilon>0$.

\claim 3.18. Remark|{\rm When $\rho_j\in\bQ$, we can take an adapted
Galois cover $g:Y\to X$ such that $(z_j\circ g)^{1-(1-s/\rho_j)_+}$
is univalent on $Y$ for all components $\Delta_j$ involved, and
we then get a well defined locally free
sheaf $\cO_Y\big(g^*V\langle D^{(s)})$ such that
$$
g^*\big(\cO_X(V\langle\lceil D\rceil\rangle)\big)\subset
\cO_Y\big(g^*V\langle D^{(s)}\rangle\big)\subset 
g^*\big(\cO_X(V)\big).
$$
However, as already stressed in Remark 1.26, this viewpoint is
not needed in our analytic approach.}
\endclaim

\plainsubsection 3.C. Orbifold tautological sheaves and their curvature|

In this context, we define the orbifold tautological sheaves
$$
\cO_{X_k(V\langle D\rangle)}(m):=
\cO_{X_k(V\langle\lceil D\rceil\rangle)}(m)\otimes
\cI((\Psi_{k,b,\varepsilon}^*)^m)
\leqno(3.19)
$$
to be the logarithmic tautological sheaves
$\cO_{X_k(V\langle\lceil D\rceil\rangle)}(m)$ 
twisted by the multiplier ideal sheaves associated
with the dual metric $\Psi_{k,b,\varepsilon}^*$ (cf.\ (3.17)),
when these are viewed
as singular hermitian metrics over the logarithmic $k$-jet bundle
$X_k(V\langle\lceil D\rceil\rangle)$. In accordance
with this viewpoint, we simply define the orbifold $k$-jet bundle to be
$X_k(V\langle D\rangle)=X_k(V\langle\lceil D\rceil\rangle)$.
The calculation of the curvature tensor is formally the same as in
the case $D=0$, and we obtain~:

\claim 3.20. Proposition| With respect to the $($rescaled$\,)$
orbifold $k$-jet components
$$
\xi_{s,\lambda}=\varepsilon_s^sf_\lambda^{(1-(1-\rho_\lambda/s)_+)}
f_\lambda^{(s)}(0),~~1\le\lambda\le p,\quad\hbox{and}\quad
\xi_{s,\lambda}=\varepsilon_s^sf_\lambda^{(s)}(0),~~
p+1\le\lambda\le r,
$$
and of the dual metric $\Psi^*_{h,b,\varepsilon}$, the curvature form of the
tautological sheaf $L_k=\cO_{X_k(V\langle D\rangle)}(1)$ admits
at any point $(z,[\xi])\in X_k(V\langle D\rangle)$
an approximate expression
$$
\leqno\displaystyle{\rm(a)}\quad
\Theta_{L_k,\Psi^*_{h,b,\varepsilon}}(z,[\xi])\simeq
\omega_{r,k,b}(\xi)+g_{V,D,k}(z,x,u),
$$
where $x_s=|\xi_s|^{2b/s}$, $u_s={\xi_s\over|\xi_s|}\in \bS^{2r-1}$ are polar
coordinates associated with $\xi_s=(\xi_{s,\lambda})_{1\le\lambda\le k}$
in~$\bC^r$, $x=(x_1,\ldots,x_k)\in\bDelta^{k-1}$,
$[\xi]=[\xi_1,\ldots,\xi_k]\in\bP(1^{[r]},2^{[r]},\ldots,k^{[r]})$ and
$$
\leqno\displaystyle{\rm(b)}\quad
g_{V,D,k}(z,x,u)={\ii\over 2\pi}
\sum_{1\le s\le k}{x_s\over s}\sum_{i,j,\lambda,\mu}c^{(s)}_{ij\lambda\mu}(z)\,
u_{s,\lambda}\overline u_{s,\mu}\,dz_i\wedge d\overline z_j.
$$
Here $({\ii\over 2\pi}\,c^{(s)}_{ij\lambda\mu})$ are the coefficients of
the curvature tensor $-\Theta_{V\langle D^{(s)}\rangle,\wt h_s}$, and the
error terms  are
$O(\max_{2\le s\le k}(\varepsilon_s/\varepsilon_{s-1})^s)$, uniformly on
the projectivized orbifold variety~$X_k(V\langle D\rangle)$.
\endclaim

\noindent
Notice, as is clear from the expressions $(3.16'')$, (3.17) and the fact that
$v_j=z_jv^\orb_j$, that our orbifold Finsler metrics always have
fiberwise positive curvature, equal to $\omega_{k,r,b}(\xi)$, along
the fibers of $X_k(V\langle D\rangle)\to X$ (even after taking into
account the so-called error terms, because fiberwise, the functions
under consideration are just sums of even powers $|\wt\xi_s^\orb|^{2b/s}$
in suitable $k$-jet components, and are therefore plurisubharmonic.)

\plainsection{4. Existence theorems for jet differentials}

\plainsubsection 4.A. Expression of the Morse integral|

Thanks to the uniform approximation provided by Proposition 3.20,
we can (and will) neglect the $O(\varepsilon_s/\varepsilon_{s-1})$ error
terms in our calculations. Since $\omega_{r,k,b}$ is positive definite on
the fibers of $X_k(V\langle D\rangle)\to X$ (at least outside of
the axes $\xi_s=0$), the index of the $(1,1)$ curvature form
$\Theta_{L_k,\Psi^*_{h,b,\varepsilon}}(z,[\xi])$ is equal to the index
of the $(1,1)$-form $g_{V,D,k}(z,x,u)$. By the binomial formula,
the $q$-index integral of $(L_k ,\Psi_{h,b,\varepsilon}^*)$ on
$X_k(V\langle D\rangle)$ is therefore equal to
$$
\leqalignno{
&\int_{X_k(V\langle D\rangle)(L_k,q)}
\Theta_{L_k,\Psi^*_{h,b,\varepsilon}}^{n+kr-1}\cr
&\qquad{}={(n+kr-1)!\over n!(kr-1)!}
\int_{z\in X}\int_{\xi\in \bP(1^{[r]},\ldots,k^{[r]})}\omega_{r,k,b}^{kr-1}(\xi)
\wedge\bOne_{g_{V,D,k},q}(z,x,u)\,g_{V,D,k}(z,x,u)^n
&(4.1)\cr}
$$
where $\bOne_{g_{V,D,k},q}(z,x,u)$ is the characteristic function
of the open set of points where $g_{V,D,k}(z,x,u)$ has signature
$(n-q,q)$ in terms of the $dz_j$'s. Notice that since
$g_{V,D,k}(z,x,u)^n$ is~a determinant, the product
$\bOne_{g_{V,D,k},q}(z,x,u)\,g_{V,D,k}(z,x,u)^n$ gives rise to
a continuous function on~$X_k(V\langle D\rangle)$. By Formula 3.10~(b),
we get
$$
\leqalignno{
&\int_{X_k(V\langle D\rangle)(L_k,q)}
\Theta_{L_k,\Psi^*_{h,b,\varepsilon}}^{n+kr-1}
={(n+kr-1)!\over n!\,k!^r(kr-1)!}~~\times\cr
&\qquad\int_{z\in X}
\int_{(x,u)\in\bDelta^{k-1}\times(\bS^{2r-1})^k}\bOne_{g_{V,D,k},q}(z,x,u)\,
g_{V,D,k}(z,x,u)^n\,d\nu_{k,r}(x)\,d\mu(u).&(4.2)\cr}
$$

\plainsubsection 4.B. Probabilistic estimate of cohomology groups|

We assume here that we are either in the ``compact'' case $(D=0)$,
or in the logarithmic case $(\rho_j=\infty)$. Then the curvature
coefficients $\smash{c^{(s)}_{ij\lambda\mu}}=c_{ij\lambda\mu}$ do not
depend on $s$
and are those of the dual bundle $V^*$ (resp.\ $V^*\langle D\rangle$).
In this situation, formula 3.20~(b) for $g_{V,D,k}(z,x,u)$ can be
thought of as
a ``Monte Carlo'' evaluation of the curvature tensor, obtained by
averaging the curvature
at random points $u_s\in \bS^{2r-1}$ with certain positive weights $x_s/s\,$; 
we then think of the \hbox{$k$-jet}
$f$ as some sort of random variable such that the derivatives 
$\nabla^kf(0)$ (resp.\ logarithmic derivatives) are uniformly
distributed in all directions. Let us compute the expected value of
$(x,u)\mapsto g_{V,D,k}(z,x,u)$ with respect to the probability measure
$d\nu_{k,r}(x)\,d\mu(u)$. Since 
$$\int_{\bS^{2r-1}}u_{s,\lambda}\overline u_{s,\mu}d\mu(u_s)={1\over r}
\delta_{\lambda\mu}\quad
\hbox{and}\quad\int_{\bDelta^{k-1}}x_s\,d\nu_{k,r}(x)={1\over k},
$$
we~find the expected value
$$
\bE(g_{V,D,k}(z,\bu,\bu))={1\over kr}
\sum_{1\le s\le k}{1\over s}\cdot{\ii\over 2\pi}\sum_{i,j,\lambda}
c_{ij\lambda\lambda}(z)\,dz_i\wedge d\overline z_j.
$$
In other words, we get the normalized trace of the curvature, i.e.
$$
\bE(g_{V,D,k}(z,\bu,\bu))={1\over kr}\Big(1+{1\over 2}+\cdots+{1\over k}\Big)
\Theta_{\det(V^*\langle D\rangle),\det h^*},
\leqno(4.3)
$$
where $\Theta_{\det(V^*\langle D\rangle),\det h^*}$ is the
$(1,1)$-curvature form of $\det(V^*\langle D\rangle)$ with the
metric induced by~$h$. It is natural to guess that 
$g_{V,D,k}(z,x,u)$ behaves asymptotically as its expected value
$\bE(g_{V,D,k}(z,\bu,\bu))$ when $k$ tends to infinity. If we replace brutally 
$g_{V,D,k}$ by its expected value in~(4.2), we get the integral
$$
{(n+kr-1)!\over n!\,k!^r(kr-1)!}{1\over (kr)^n}
\Big(1+{1\over 2}+\cdots+{1\over k}\Big)^n\int_X\bOne_{\eta,q}\eta^n,
$$
where $\eta:=\Theta_{\det(V^*\langle D\rangle),\det h^*}$ and
$\bOne_{\eta,q}$ is the
characteristic function of its $q$-index set in~$X$. The leading constant is
equivalent to $(\log k)^n/n!\,k!^r$ modulo 
a multiplicative factor \hbox{$1+O(1/\log k)$}. By working out a more
precise analysis of the deviation, the following result has been
proved in [Dem11] in the compact case; the more general logarithmic case
can be treated without any change, so we state the result in this situation
by just transposing the results of [Dem11].

\claim 4.4. Probabilistic estimate|Let $(X,V,D)$ be a non singular
logarithmic directed variety.
Fix smooth hermitian metrics $\omega$ on $T_X$, $h$ on $V\langle D\rangle$,
and write
$\omega={\ii\over 2\pi} \sum\omega_{ij}dz_i\wedge d\overline z_j$ on~$X$. 
Denote by $\Theta_{V\langle D\rangle,h}=-{\ii\over 2\pi}\sum
c_{ij\lambda\mu}dz_i\wedge d\overline z_j\otimes e_\lambda^*\otimes
e_\mu$ the curvature tensor of $V\langle D\rangle$ with respect to
an $h$-orthonormal frame $(e_\lambda)$, and put
$$
\eta(z):=\Theta_{\det(V^*\langle D\rangle),\det h^*}=
{\ii\over 2\pi}\sum_{1\le i,j\le n}\eta_{ij}
dz_i\wedge d\overline z_j,\qquad
\eta_{ij}:=\sum_{1\le\lambda\le r}c_{ij\lambda\lambda}.
$$
Finally consider the $k$-jet line bundle
$L_k=\smash{\cO_{X_k(V\langle D\rangle)}(1)}\to
X_k(V\langle D\rangle)$ equipped with the induced metric
$\Psi^*_{h,b,\varepsilon}$
$($as defined above, with $1=\varepsilon_1\gg\varepsilon_2\gg\ldots\gg
\varepsilon_k>0)$. When $k$ tends 
to infinity, the integral of the top power of the curvature of $L_k$ on its
$q$-index set $X_k(V\langle D\rangle)(L_k,q)$ is given by
$$
\int_{X_k(V\langle D\rangle)(L_k,q)}
\Theta_{L_k,\Psi^*_{h,b,\varepsilon}}^{n+kr-1}=
{(\log k)^n\over n!\,k!^r}\bigg(
\int_X\bOne_{\eta,q}\eta^n+O((\log k)^{-1})\bigg)
$$
for all $q=0,1,\ldots,n$, and the error term $O((\log k)^{-1})$ can be 
bounded explicitly in terms of $\Theta_{V\langle D\rangle}$, $\eta$
and $\omega$. Moreover, the  left hand side is identically zero for $q>n$.
\endclaim

The final statement follows from the observation that the curvature of
$L_k$ is positive along the fibers of $X_k(V\langle D\rangle)\to X$, by the 
plurisubharmonicity of the weight (this is true even 
when the error terms are taken into account, since they
depend only on the base); therefore the $q$-index sets are empty for
$q>n$. It will be useful to extend the above estimates to the 
case of sections of
$$
L_{F,k}=\cO_{X_k(V\langle D\rangle)}(1)\otimes
\pi_k^*\cO\Big(-{1\over kr}\Big(1+{1\over 2}+\cdots+{1\over k}\Big)F\Big)
\leqno(4.5)
$$
where $F\in\Pic_\bQ(X)$ is an arbitrary $\bQ$-line bundle on~$X$ and 
$\pi_k:X_k(V\langle D\rangle)\to X$ is the natural projection. We assume here
that $F$ is also equipped with a smooth hermitian metric $h_F$. In formulas
(4.2--4.4), the curvature $\Theta_{L_{F,k}}$ of $L_{F,k}$ takes 
the form $\Theta_{L_{F,k}}=\omega_{r,k,b}(\xi)+g_{V,D,F,k}(z,x,u)$ where
$$
g_{V,D,F,k}(z,x,u)=g_{V,D,k}(z,x,u)-
{1\over kr}\Big(1+{1\over 2}+\cdots+{1\over k}\Big)\Theta_{F,h_F}(z),
\leqno(4.6)
$$
and by the same calculations its normalized expected value is
$$
\eta_F(z):={1\over{1\over kr}(1+{1\over 2}+\cdots+{1\over k})}
\bE(g_{V,D,F,k}(z,\bu,\bu))=
\Theta_{\det V^*\langle D\rangle,\det h^*}(z)-\Theta_{F,h_F}(z).
\leqno(4.7)
$$
Then the variance estimate for $g_{V,D,F,k}$ is the same as the
variance estimate for $g_{V,D,k}$, and the recentered
$L^p$ bounds are still valid, since our forms are just shifted
by adding the constant smooth term $\Theta_{F,h_F}(z)$. The probabilistic
estimate 4.4 is therefore still true in exactly the same form for $L_{F,k}$,
provided we use $g_{V,D,F,k}$ and $\eta_F$ instead of $g_{V,D,k}$
and $\eta$. An application of holomorphic Morse inequalities gives the 
desired cohomology estimates for 
$$
\eqalign{
h^q\Big(X,E_{k,m}V^*\langle D\rangle&{}\otimes
\cO\Big(-{m\over kr}\Big(1+{1\over 2}+\cdots+{1\over k}\Big)F\Big)\Big)\cr
&{}=h^q(X_k(V\langle D\rangle),\cO_{X_k(V\langle D\rangle)}(m)\otimes
\pi_k^*\cO\Big(-{m\over kr}\Big(1+{1\over 2}+\cdots+{1\over k}\Big)F\Big)\Big),
\cr}
$$
provided $m$ is sufficiently divisible to give a multiple of $F$ which
is a $\bZ$-line bundle.

\claim 4.8. Theorem| Let $(X,V\langle D\rangle)$ be a non singular
logarithmic directed variety, \hbox{$F\to X$} a
\hbox{$\bQ$-line} bundle, $(V\langle D\rangle,h)$ and $(F,h_F)$ smooth
hermitian structure on $V\langle D\rangle$ 
and $F$ respectively. We define
$$
\eqalign{
L_{F,k}&=\cO_{X_k(V\langle D\rangle)}(1)\otimes
\pi_k^*\cO\Big(-{1\over kr}\Big(1+{1\over 2}+\cdots+{1\over k}\Big)F\Big),\cr
\eta_F&=\Theta_{\det V^*\langle D\rangle,\det h^*}-\Theta_{F,h_F}
=\Theta_{\det V^*\langle D\rangle\otimes F^{-1},\det h^*}.\cr}
$$
Then for all $q\ge 0$ and all $m\gg k\gg 1$ such that 
m is sufficiently divisible, we have
$$\leqalignno{\kern20pt
h^q(X_k(V\langle D\rangle),\cO(L_{F,k}^{\otimes m}))
&\le {m^{n+kr-1}\over (n+kr-1)!}{(\log k)^n\over n!\,k!^r}\bigg(
\int_{X(\eta_F,q)}(-1)^q\eta_F^n+O((\log k)^{-1})\bigg),&\hbox{\rm(a)}\cr
h^0(X_k(V\langle D\rangle),\cO(L_{F,k}^{\otimes m}))
&\ge {m^{n+kr-1}\over (n+kr-1)!}
{(\log k)^n\over n!\,k!^r}\bigg(
\int_{X(\eta_F,\le 1)}\eta_F^n-O((\log k)^{-1})\bigg),&\hbox{\rm(b)}\cr
\cr
\chi(X_k(V\langle D\rangle),\cO(L_{F,k}^{\otimes m}))&={m^{n+kr-1}\over (n+kr-1)!}
{(\log k)^n\over n!\,k!^r}\big(
c_1(V^*\langle D\rangle\otimes F)^n+O((\log k)^{-1})\big).
&\hbox{\rm(c)}\cr
\cr}
$$
\vskip-4pt
\endclaim

Green and Griffiths [GrGr80] already checked the Riemann-Roch
calculation (4.8$\,$c) in the special case $D=0$,
$V=T_X^*$ and $F=\cO_X$. Their proof is much simpler since it relies only
on Chern class calculations, but it cannot provide any information on
the individual cohomology groups, except in very special cases where
vanishing theorems can be applied; in fact in dimension 2, the
Euler characteristic satisfies $\chi=h^0-h^1+h^2\le h^0+h^2$, hence
it is enough to get the vanishing of the top cohomology group $H^2$
to infer $h^0\ge\chi\,$; this works for surfaces by means of a well-known
vanishing theorem of Bogomolov which implies in general
$$H^n\bigg(X,E_{k,m} T_X^*\otimes\cO\Big(-{m\over kr}
\Big(1+{1\over 2}+\cdots+{1\over k}\Big)F\Big)\Big)\bigg)=0
$$
as soon as $K_X\otimes F^{-1}$ is big and $m\gg 1$.

In fact, thanks to Bonavero's singular holomorphic Morse inequalities 
(Theorem 2.9, cf.\ [Bon93]), everything works almost unchanged in the
case where the metric $h$ on $V$ is taken to a product
$h=h_\infty e^\varphi$ of a smooth
metric $h_\infty$ by the exponential of a quasi-plurisubharmonic
weight~$\varphi$ with analytic singularities (so that
$\det(h^*)=\det(h_\infty^*)e^{-r\varphi}$). Then $\eta$ is a
$(1,1)$-current with logarithmic poles, and we just have to twist
our cohomology groups by the appropriate multiplier ideal
sheaves $\cI_{k,m}$ associated with the weight
${1\over k}(1+{1\over 2}+\cdots+{1\over k})m\,\varphi$, since this
is the multiple of $\det V^*$ that occurs in the calculation, up to
the factor ${1\over r}\times r\varphi$. The corresponding Morse
integrals need only
be evaluated in the complement of the poles, i.e., on
$X(\eta,q)\ssm S$ where $S=\Sing(\varphi)$. Since
$$
(\pi_k)_*\big(\cO(L_{F,k}^{\otimes m})\otimes\cI_{k,m}\big)\subset
E_{k,m} V^*\otimes
\cO\Big(-{m\over kr}\Big(1+{1\over 2}+\cdots+{1\over k}\Big)F\Big)\Big)
$$
we still get a lower bound for the $H^0$ of the latter sheaf (or for the $H^0$
of the un-twisted line bundle $\cO(L_k^{\otimes m})$ on $\smash{X_k(V)}$).
If we assume that $K_V\otimes F^{-1}$ is big, these considerations
also allow us to obtain a strong estimate in terms of the volume, by
using an approximate Zariski decomposition on a suitable blow-up of~$X$.

\claim 4.9. Corollary|
If $F$ is an arbitrary $\bQ$-line bundle over~$X$, one has
$$
\eqalign{
h^0\bigg(&X_k(V),\cO_{X_k(V)}(m)\otimes\pi_k^*\cO
\Big(-{m\over kr}\Big(1+{1\over 2}+\cdots+{1\over k}\Big)F\Big)\bigg)\cr
&\ge {m^{n+kr-1}\over (n+kr-1)!}
{(\log k)^n\over n!\,k!^r}\Big(
\Vol(K_V\otimes F^{-1})-O((\log k)^{-1})\Big)-o(m^{n+kr-1}),\cr}
$$
when $m\gg k\gg 1$, in particular there are many sections of the
$k$-jet differentials of degree $m$ twisted by the appropriate
power of $F$ if $K_V\otimes F^{-1}$ is big.
\endclaim

\plainproof. The volume is computed here as usual, i.e.\ after performing a
suitable modifi\-cation $\mu:\smash{\wt X}\to X$ which converts $K_V$ into 
an invertible sheaf. There is of course nothing to prove if
$K_V\otimes F^{-1}$ is not big, so we can assume $\Vol(K_V\otimes F^{-1})>0$.
Let us fix smooth hermitian metrics $h_0$ on $T_X$ and
$h_F$ on $F$. They induce a metric $\mu^*(\det h_0^{-1}\otimes h_F^{-1})$
on $\mu^*(K_V\otimes F^{-1})$ which, by our definition of $K_V$, is
a smooth metric. By the result of Fujita [Fuj94] on
approximate Zariski decomposition, for every $\delta>0$, one can find
a modification $\mu_\delta:\smash{\wt X_\delta}\to X$ dominating
$\mu$ such that
$$
\mu_\delta^*(K_V\otimes F^{-1}) =\cO_{\wt X_\delta}(A+E)
$$
where $A$ and $E$ are $\bQ$-divisors, $A$ ample and $E$ effective,
with 
$$\Vol(A)=A^n\ge \Vol(K_V\otimes F^{-1})-\delta.$$
If we take a smooth metric $h_A$ with positive definite curvature form
$\Theta_{A,h_A}$, then we get a singular hermitian metric $h_Ah_E$ on
$\mu_\delta^*(K_V\otimes F)$ with poles along $E$, i.e.\ the quotient
$h_Ah_E/\mu^*(\det h_0^{-1}\otimes h_F)$ is of the form $e^{-\varphi}$ where
$\varphi$ is quasi-psh with log poles $\log|\sigma_E|^2$ 
(mod $C^\infty(\smash{\wt X_\delta}))$ precisely given
by the divisor~$E$. We then only need to take the singular metric $h$
on $T_X$ defined by
$$
h=h_0e^{{1\over r}(\mu_\delta)^*\varphi}
$$
(the choice of the factor ${1\over r}$ is there to correct adequately 
the metric on $\det V$). By construction $h$ induces an 
admissible metric on $V$ and the resulting 
curvature current $\eta_F=\Theta_{K_V,\det h^*}-\Theta_{F,h_F}$ is such that
$$
\mu_\delta^*\eta_F = \Theta_{A,h_A} +[E],\qquad
\hbox{$[E]={}$current of integration on $E$.}
$$
Then the $0$-index Morse integral in the complement of the poles 
is given by
$$
\int_{X(\eta,0)\ssm S}\eta_F^n=\int_{\wt X_\delta}\Theta_{A,h_A}^n=A^n\ge
\Vol(K_V\otimes F^{-1})-\delta
$$
and Corollary 4.9 follows from the fact that $\delta$ can be taken arbitrary 
small.\qed

\claim 4.10. Remark|{\rm Since the probability estimate requires
$k$ to be very large, and since all non logarithmic components disappear
from $D^{(s)}$ when $s$ is large, the above lower bound does not work
in the general orbifold case. In that case, one can only hope to get
an interesting result when $k$ is fixed and not too large. This is what
we will do in \S$\,$6.}
\endclaim

\plainsection{5. Positivity concepts for vector bundles
and Chern inequalities}
  
\plainsubsection 5.A. Griffiths, Nakano and strong (semi-)positivity|

Let $E\to X$ be a holomorphic vector bundle equipped with a hermitian metric.
Then $E$ possesses a uniquely defined Chern connection $\nabla_h$ compatible
with $h$ and such that $\nabla_h^ {0,1}=\dbar$. The curvature tensor of $(E,h)$
is defined to be
$$
\Theta_{E,h}:={\ii\over 2\pi}
\ii\ddbar\nabla_h^2\in C^\infty(X,\Lambda^{1,1}T^*_X\otimes\Hom(E,E)).
\leqno(5.1)
$$
One can then associate bijectively to $\Theta_{E,h}$ a hermitian form
$\wt\Theta_{E,h}$ on $TX\otimes E$, such that
$$
\wt\Theta_{E,h}(\xi\otimes u,\xi\otimes u)=
\langle\Theta_{E,h}(\xi,\xi)\cdot u,u\rangle_h.\leqno(5.2)
$$
and can be written
$$
\Theta_{E,h}=
{\ii\over 2\pi}
\sum_{i,j,\lambda,\mu}
c_{ij\lambda\mu}\,dz_i\wedge d\overline z_j\otimes e_\lambda^*\otimes e_\mu
$$
Let $(z_1,\ldots,z_n)$ be a holomorphic coordinate system  and
$(e_\lambda)_{1\le\lambda\le r}$ a smooth frame of $E$. If $(e_\lambda)$ is
chosen to be orthonormal, then we can write
$$
\leqalignno{
\Theta_{E,h}&={\ii\over 2\pi}
\sum_{i,j,\lambda,\mu}
c_{ij\lambda\mu}\,dz_i\wedge d\overline z_j\otimes e_\lambda^*\otimes e_\mu,
&(5.3)\cr
\wt\Theta_{E,h}(\xi\otimes u,\xi\otimes u)&={1\over 2\pi}
\sum_{i,j,\lambda,\mu}
c_{ij\lambda\mu}\,\xi_i\overline\xi_j\,u_\lambda\overline u_\mu,
&(5.3')\cr}
$$
and more generally $\wt\Theta_{E,h}(\tau,\tau)={1\over 2\pi}
\sum_{i,j,\lambda,\mu}c_{ij\lambda\mu}\,\tau_{i\lambda}\overline
\tau_{j\mu}$ for every tensor $\tau\in T_X\otimes E$. We now consider
three concepts of (semi-)positivity, the first two being very classical.

\claim 5.4. Definition|Let $\theta$ be a hermitian form on a tensor
product $T\otimes E$ of complex vector spaces. We say that
\plainitem{\rm (a)} $\theta$ is Griffiths semi-positive if
$\theta(\xi\otimes u,\xi\otimes u)\ge 0$ for every
$\xi\in T$ and every $v\in E;$
\vskip2pt
\plainitem{\rm (b)} $\theta$ is Nakano semi-positive if
$\theta(\tau,\tau)\ge 0$ for every
$\tau\in T\otimes E\,;$
\vskip2pt
\plainitem{\rm (c)} $\theta$ is strongly semi-positive if
there exist a finite collection of linear forms $\alpha_j\in T^*$,
$\psi_j\in E^*$ such that $\theta=\sum_j |\alpha_j\otimes\psi_j|^2$, i.e.
$$
\theta(\tau,\tau)=\sum_j|(\alpha_j\otimes\psi_j)\cdot \tau|^2,\quad
\forall\tau\in T\otimes E.
$$
Semi-negativity concepts are introduced in a similar way.
\vskip2pt
\plainitem{\rm(d)} We say that the hermitian bundle $(E,h)$ is Griffiths
semi-positive,
resp.\ Nakano semi-positive, resp.\ strongly semi-positive, if
$\wt\Theta_{E,h}(x)
\in \Herm(T_{X,x}\otimes E_x)$ satisfies the corresponding property for
every point $x\in X$.
\vskip2pt
\plainitem{\rm(e)}
$($Strict$)$ Griffiths positivity means that
$\wt\Theta_{E,h}(\xi\otimes u,\xi\otimes u)>0$ for every non zero vectors
$\xi\in T_{X,x}$, $v\in E_x$.
\vskip2pt
\plainitem{\rm(f)} $($Strict$)$ strong positivity means that at
every point $x\in X$ we can decompose $\wt\Theta_{E,h}$ as
$\wt\Theta_{E,h}=\sum_j|\alpha_j\otimes\psi_j|^2$ where
$\Span(\alpha_j\otimes\psi_j)=T^*_{X,x}\otimes E^*_x$.
\vskip2pt
\endclaim

\noindent
We will denote respectively by $\ge_G$, $\ge_N$, $\ge_S$
(and $>_G$, $>_N$, $>_S$) the Griffiths, Nakano, strong 
(semi-)positivity relations. It is obvious that
$$
\theta\ge_S 0~~\Rightarrow~~\theta\ge_N 0~~\Rightarrow~~\theta\ge_G 0,
$$
and one can show that the reverse implications do not hold when
$\dim T>1$ and $\dim E>1$. The following result from [Dem80]
will be useful.

\claim 5.5. Proposition|Let $\theta\in\Herm(T\otimes E)$, where
$(E,h)$ is a hermitian vector space. We define $\Tr_E(\theta)\in\Herm(T)$
to be the hermitian form such that
$$
\Tr_E(\theta)(\xi,\xi)=\sum_{1\le\lambda\le r}
\theta(\xi\otimes e_\lambda,\xi\otimes e_\lambda)
$$
where $(e_\lambda)_{1\le\lambda\le r}$ is an arbitrary orthonormal basis of $E$.
Then
$$
\theta\ge_G 0~~\Longrightarrow~~
\theta+\Tr_E(\theta)\otimes h\ge_S 0.
$$
As a consequence, if $(E,h)$ is a Griffiths $($semi-$)$positive vector bundle,
then the tensor product
$(E\otimes \det E,h\otimes\det(h))$ is strongly $($semi-$)$positive.
\endclaim

\plainproof. Since [Dem80] is written in French and
perhaps not so easy to find, we repeat here briefly the arguments.
They are based on a Fourier inversion formula for discrete Fourier transforms.

\claim 5.6. Lemma|Let $q$ be an integer $\ge 3$, and
$x_\alpha,~y_\beta,~1\le\alpha,\beta\le r$, be complex numbers.
Let $\chi$ describe
the set $U^r_q$ of $r$-tuples of $q$-th roots of unity and put
$$\wh x(\chi) = \sum_{1\le\alpha\le r} x_\alpha\ol \chi_\alpha,~~~~
\wh y(\chi) = \sum_{1\le\beta\le r} y_\beta\ol\chi_\beta,~~~~ \chi\in U^r_q.$$
Then for every pair $(\lambda,\mu),~1\le\lambda,\mu\le r$, the following 
identity holds:
$$
q^{-r} \sum_{\chi\in U^r_q}
\wh x(\chi)\,\ol{\wh y(\chi)}\,\chi_\lambda\ol \chi_\mu
=\plaincases{x_\lambda\ol y_\mu&if~~$\lambda\ne\mu,$\cr
\noalign{\vskip6pt}  
\sum_{1\le\alpha\le r} x_\alpha\ol y_\alpha&if~~$\lambda=\mu.$\cr}$$
\endclaim

\noindent
In fact, the coefficient of $x_\alpha\ol y_\beta$ in the summation
$q^{-r} \sum_{\chi\in U^r_q}\wh x(\chi)\,\ol{\wh y(\chi)}\,
\chi_\lambda\ol \chi_\mu$ is given by
$$
q^{-r} \sum_{\chi\in U^r_q}\chi_\alpha \ol \chi_\beta\ol
\chi_\lambda \chi_\mu,
$$
so it is equal to $1$ when the pairs $\{ \alpha,\mu\}$ and
$\{\beta,\lambda\}$ coincide, and is equal to $0$ otherwise.
The identity stated in Lemma~5.6 follows immediately.\qed

\noindent
Now, let $(t_j)_{1\le j\le n}$ be a basis of
$T$, $(e_\lambda)_{1\le \lambda\le r}$ an orthonormal basis of $E$ and
$\xi=\sum_j\xi_jt_j\in T$, $w=\sum_{j,\lambda}w_{j\lambda}\,t_j\otimes 
e_\lambda \in T\otimes E$. The coefficients $c_{jk\lambda\mu}$ of $\theta$
with respect to the basis $t_j\otimes e_\lambda$ satisfy the symmetry
relation $\ol c_{jk\lambda\mu}=c_{kj\mu\lambda}$, and we have the formulas 
$$\eqalign{\theta(w,w)
&=\sum_{j,k,\lambda,\mu}c_{jk\lambda\mu}w_{j\lambda}\ol w_{k\mu},\quad
\Tr_E\theta(\xi,\xi)=\sum_{j,k,\lambda}c_{jk\lambda\lambda}\xi_j\ol\xi_k,\cr
(\theta+\Tr_E\theta\otimes h)(w,w)&=\sum_{j,k,\lambda,\mu}c_{jk\lambda\mu} 
w_{j\lambda} \ol w_{k\mu}+ c_{jk\lambda\lambda}w_{j\mu}\ol w_{k\mu}.\cr}$$
For every $\chi\in U^r_q$, let us put
$$
\wh w_j(\chi)=\sum_\alpha w_{j\alpha}\ol\chi_\alpha,\quad
\wh w(\chi)=\sum_j\wh w_j(\chi)\,t_j\in T\,,\quad
\wh e_\chi=\sum_\lambda\chi_\lambda e_\lambda\in E.$$
Lemma 5.6 implies
$$\eqalign{
q^{-r} \sum_{\chi\in U^r_q}\theta(\wh w(\chi)\otimes\wh e_\chi,
\wh w(\chi)\otimes\wh e_\chi)
&=q^{-r}\sum_{\chi\in U^r_q}~\sum_{j,k,\lambda,\mu}c_{jk\lambda\mu}\,
\wh w_j(\chi)\ol{\wh w_k(\chi)}\,\chi_\lambda \ol \chi_\mu\cr
&= \sum_{j,k,\lambda\ne \mu} c_{jk\lambda\mu} w_{j\lambda} \ol w_{k\mu} +
\sum_{j,k,\lambda,\mu} c_{jk\lambda\lambda} w_{j\mu} \ol w_{k\mu}.\cr}$$
The Griffiths positivity assumption $\theta_G\ge 0$ shows that
$\xi\mapsto q^{-r}\,\theta(\xi\otimes\wh e_\chi,\xi\otimes\wh e_\chi)$
is a semi-positive hermitian form on~$T$, hence there are linear forms
$\ell_{\chi,j}\in T^*$ such that $q^{-r}\,\theta(\xi\otimes\wh e_\chi,
\xi\otimes\wh e_\chi)=\sum_j|\ell_{\chi,j}(\xi)|^2$ for all~$\xi\in T$.
Similarly, there are $\ell'_{\lambda,j}\in T^*$ such that
$$
\sum_{j,k} c_{jk\lambda\lambda}\,\xi_j\overline\xi_k
=\sum_j|\ell'_{\lambda,j}(\xi)|^2,\quad
\hbox{for all $\lambda=1,\ldots,r$}.
$$
Our final Fourier identity can be rewritten
$$
\eqalign{
(\theta + \Tr_E \theta \otimes h)(w,w)
&=\sum_{j,k,\lambda,\mu}c_{jk\lambda\mu}w_{j\lambda} \ol w_{k\mu}+
\sum_{j,k,\lambda,\mu}c_{jk\lambda\lambda}w_{j\mu} \ol w_{k\mu}\cr
&=q^{-r} \sum_{\chi\in U^r_q}\theta(\wh w(\chi)\otimes\wh e_\chi,
\wh w(\chi)\otimes\wh e_\chi)+\sum_{j,k,\lambda} c_{jk\lambda\lambda}\,
w_{j\lambda}\overline w_{j\lambda}\cr
&=\sum_{\chi\in U^r_q}\sum_j|\ell_{\chi,j}(\wh w(\chi))|^2+\sum_{j,\lambda}
|\ell'_{\lambda,j}(w_{{\scriptscriptstyle\bullet},\lambda})|^2\cr
&=\sum_{\chi\in U^r_q}\sum_j|\ell_{\chi,j}\otimes\chi^*(w)|^2+
\sum_{j,\lambda}|\ell'_{\lambda,j}\otimes e_\lambda^*(w)|^2\cr}
$$
where $\chi^*=\langle\bu,\chi\rangle\in E^*$, thus
$\theta + \Tr_E \theta \otimes h\ge_S 0$.\qed

\claim 5.7. Corollary|Let $r=\dim E$ and $\Theta\in\Herm(T\otimes E)$.
\vskip2pt
\plainitem{\rm(a)} If $\theta\ge_G 0$, then\kern27pt
$-\Tr_E\theta\otimes h~\le_S~\theta~\le_S~r\,\Tr_E\theta\otimes h$.
\vskip2pt
\plainitem{\rm(b)} If $\theta\le_G 0$, then~
$-r\,\Tr_E(-\theta)\otimes h~\le_S~\theta~\le_S~\Tr_E(-\theta)\otimes h$.
\vskip2pt
\plainitem{\rm(c)} If $\pm\theta\le_G\tau\otimes h$ where $\tau\in\Herm(T)$ is
semi-positive, then
$$
-(2r+1)\,\tau\otimes h~\le_S~\theta~\le_S~(2r+1)\,\tau\otimes h.
$$
\endclaim

\plainproof. (a) It is easy to chech that $\theta'=\Tr_E\theta\otimes h-\theta$
satisfies $\theta'\ge_G 0$ and that we have
$\Tr_E\theta'=(r-1)\Tr_E\theta$. Lemma~5.6 implies
$$
\theta'+\Tr_E\theta'\otimes h=r\,\Tr_E\theta\otimes h-\theta\ge_S 0.
$$
(b) follows from (a), after replacing $\theta$ with $-\theta$.
\smallskip
\noindent
(c) also follows from Lemma~5.6 by taking $\theta'=\tau\otimes h+\theta$
(resp.\ $\theta'=\tau\otimes h-\theta$),
since $\Tr_E\theta\le r\,\tau$ and we have e.g.
$$
0\le_S\theta'+\Tr_E\theta'\otimes h=
\theta+\Tr_E\theta\otimes h+(r+1)\tau\otimes h
\le_S\theta+(2r+1)\tau\otimes h.
\eqno\square
$$

\plainsubsection 5.B. Chern form inequalities|

In view of the estimates developed in section~6, we will have
to evaluate integrals involving powers of curvature tensors, and
the following basic inequalities will be useful.

\claim 5.8. Lemma|Let $\ell_j\in(\bC^r)^*$, $1\le j\le p$, be non zero
complex linear forms on $\bC^r$, where $(\bC^r)^*\simeq\bC^r$ is
equipped with its standard hermitian form, and let $\mu$ the rotation
invariant probability measure on $\bS^{2r-1}\subset\bC^r$. Then
$$
I(\ell_1,\ldots,\ell_p)=\int_{\bS^{2r-1}}|\ell_1(u)|^2\ldots\,|\ell_p(u)|^2\,
d\mu(u)
$$
satisfies the following inequalities$\;:$
\vskip2pt
\plainitem{\rm(a)} $\displaystyle I(\ell_1,\ldots,\ell_p)\le
{p!\,(r-1)!\over (p+r-1)!}\,\prod_{j=1}^p|\ell_j|^2,$\vskip2pt
and the equality occurs if and only if the
$\ell_j$ are proportional$\,;$
\vskip2pt
\plainitem{\rm(b)} $\displaystyle I(\ell_1,\ldots,\ell_p)\ge
{(r-1)!\over (p+r-1)!}\,\prod_{j=1}^p|\ell_j|^2,$\vskip2pt
and the equality occurs if and only if $p\le r$ and the
$\ell_j$ are pairwise orthogonal.
\endclaim

\plainproof. Denote by $d\lambda$ the Lebesgue measure on Euclidean space
and by $d\sigma$ the area measure of the sphere. One can easily check that
the projection
$$
\bS^{2r-1}\to \bB^{2r-2},\quad
u=(u_1,\ldots,u_r)\mapsto v=(u_1,\ldots,u_{r-1}),
$$
yields $d\sigma(u)=d\theta\wedge d\lambda(v)$ where $u_r=|u_r|\,e^{i\theta}$
$[\,$just check that the wedge products of both sides with ${1\over 2}d|u|^2$
are equal to $d\lambda(u)$, and use the fact that $d\theta={1\over 2i}(du_r/u_r
-d\overline u_r/\overline u_r)$], thus, in terms of polar
coordinates $v=t\,u'$, $u'\in\bS^{2r-1}$, we have
$d\sigma(u)=d\theta\wedge t^{2r-3}\,dt\wedge d\sigma'(u')$,
and going back to the invariant probability measures $\mu$ on $\bS^{2r-1}$ and
$\mu'$ on $\bS^{2r-3}$, we get
$|u_r|^2=1-|v|^2=1-t^2$ and an equality
$$
d\mu(u)={2r-2\over 2\pi}\,d\theta\wedge t^{2r-3}\,dt
\wedge d\mu'(u').
\leqno(5.9)
$$
If $\ell_1,\ldots,\ell_p$ are independent of $u_r$, (5.9) and
the Fubini theorem imply by homogeneity
$$\leqalignno{
\qquad&\int_{\bS^{2r-1}}|\ell_1(u')|^2\ldots\,|\ell_p(u')|^2\,d\mu(u)
={r-1\over p+r-1}
\int_{\bS^{2r-3}}|\ell_1(u')|^2\ldots\,|\ell_p(u')|^2\,d\mu'(u'),
&(5.10)\cr
\noalign{\vskip8pt}
&\int_{\bS^{2r-1}}|\ell_1(u')|^2\ldots\,|\ell_{p-1}(u')|^2\,|u_r|^2\,d\mu(u)
=\cr
& \kern101pt{r-1\over (p+r-2)(p+r-1)}
\int_{\bS^{2r-3}}|\ell_1(u')|^2\ldots\,|\ell_{p-1}(u')|^2\,d\mu'(u')
&(5.10')\cr}
$$
(for instance, in case $(5.10')$, we have to integrate $t^{2p-2}(1-t^2)\times
t^{2r-3}\,dt$). The formulas
$$
\int_{\bS^{2r-1}}|u_1|^{2p}\,d\mu(u)={p!\,(r-1)!\over (p+r-1)!},\quad
\int_{\bS^{2r-1}}|u_1|^2\ldots\,|u_p|^2\,d\mu(u)={(r-1)!\over (p+r-1)!}\quad
(p\le r),
$$
are then obtained by induction on $r$ and $p$.\medskip

\noindent
(a) For any $\ell\in(\bC^r)^*$, we can find orthonormal coordinates
on $\bC^r$ such that $\ell(u)=|\ell|\,u_1$ in the new
coordinates. Hence
$$
\int_{\bS^{2r-1}}|\ell(u)|^{2p}\,d\mu(u)=m_{r,p}\,|\ell|^{2p}\quad
\hbox{where}~~m_{r,p}=\int_{\bS^{2r-1}}|u_1|^{2p}\,d\mu(u)=
{p!\,(r-1)!\over (p+r-1)!}.    
$$
It follows from H\"older's inequality that
$$
I(\ell_1,\ldots,\ell_p)\le \prod_{j=1}^p
\bigg(\int_{\bS^{2r-1}}|\ell_j|^{2p}\,d\mu(u)\bigg)^{1/p}=m_{r,p}
\prod_{j=1}^p|\ell_j|^2,
$$
and that the equality occurs if and only if all $\ell_j$ are proportional.
\medskip

\noindent (b) We prove the inequality
$$
I(\ell_1,\ldots,\ell_p)\ge {(r-1)!\over (p+r-1)!}\prod_{j=1}^p|\ell_j|^2
$$
by induction on $p$, the result being clear for $p=0$ or $p=1$. If we
choose an orthonormal basis $(e_1,\ldots,e_r)\in\bC^r$ such that
$\ell_j(e_r)\ne 0$ for all $j$ and replace $\ell_j$ by
$(\ell_j(e_r))^{-1}\ell_j$, we can assume $\ell_j(e_r)=1$. We then
write $u=u'+u_re_r$ with $u'\in e_r^\perp\simeq\bC^{r-1}$ and
$$
\ell_j(u)=\ell_j'(u')+u_r,\quad 1\le j\le p,\quad\ell'_j:=\ell_{j|e_r^\perp}.
$$
Let $s_k(\ell'_\bu(u'))$ be the elementary symmetric functions
in $\ell'_j(u')$, $1\le j\le p$, with $s_0:=1$. We have
$$
I(\ell_1,\ldots,\ell_p)
=\int_{\bS^{2r-1}}\prod_{j=1}^p|\ell'_j(u')+u_r|^2\,d\mu(u)
=\int_{\bS^{2r-1}}\Bigg|\sum_{k=0}^ps_k(\ell'_\bu(u'))\,u_r^{p-k}\Bigg|^2d\mu(u).
$$
We make a change of variable $u_r\mapsto u_r\,e^{i\theta}$ and take the average
over $\theta\in[0,2\pi]$. Parseval's formula gives
$$
I(\ell_1,\ldots,\ell_p)
=\int_{\bS^{2r-1}}\sum_{k=0}^p\big|s_k(\ell'_\bu(u'))\big|^2\,|u_r|^{2(p-k)}
d\mu(u),
$$
and since
$$
(2r-2)\int_0^1t^{2k}(1-t^2)^{p-k}\,t^{2r-3}dt=
{(r-1)\,(k+r-2)!\,(p-k)!\over(p+r-1)!},
$$
formula (5.9) implies
$$
I(\ell_1,\ldots,\ell_p)
=\int_{\bS^{2r-3}}\sum_{k=0}^p{(r-1)\,(k+r-2)!\,(p-k)!\over(p+r-1)!}\,
\big|s_k(\ell'_\bu(u'))\big|^2\,d\mu'(u').
$$
As $|\ell_j|^2=1+|\ell'_j|^2$, our inequality (5.8~(b)) is equivalent to
$$
\int_{\bS^{2r-3}}\sum_{k=0}^p{(k+r-2)!\,(p-k)!\over(r-2)!}\,
\big|s_k(\ell'_\bu(u'))\big|^2\,d\mu'(u')\ge
\prod_{j=1}^p(1+|\ell'_j|^2)
\leqno(5.11)
$$
for all linear forms $\ell'_j\in(\bC^{r-1})^*$. We actually prove (5.11) by
induction on $p$ (observing that the inequality is a trivial equality
for $p=0,1$). Assume that (5.11) (and hence (5.8~(b))) is known
for any $(p-1)$-tuple of linear forms $(\ell'_1,\ldots,\ell'_{p-1})$.
As (5.8~(b)) is invariant under the action
of $U(r)$, it is sufficient to consider the case when $\ell_p(u)=u_r$,
i.e.\ $\ell'_p=0$. The induction hypothesis tells us that
$$
\int_{\bS^{2r-3}}\sum_{k=0}^{p-1}{(k+r-2)!\,(p-1-k)!\over(r-2)!}\,
\big|s_k(\ell'_\bu(u'))\big|^2\,d\mu'(u')\ge
\prod_{j=1}^{p-1}(1+|\ell'_j|^2).
$$
However, when we add the factor $\ell_p$, the elementary symmetric
functions $s_k(\ell'_\bu(u'))$ are left unchanged for $k\le p-1$, while
$s_p(\ell'_\bu(u'))=0$ and $1+|\ell'_p|^2=1$. Therefore (5.11) holds true
for $p$, since $(p-k)!\ge (p-1-k)!$ for all $k=0,1,\ldots,p-1$.
We have proved the inequality at order $p$ whenever
$\ell_p=\alpha_p\langle\bu,e_r\rangle$ and $\ell_j(e_r)\ne 0$ for $j\le p-1$.
Since those $(\ell_1,\ldots,\ell_p)$ are dense in the space $((\bC^r)^*)^p$
of $p$-tuples of linear forms, the proof of the lower bound is
complete.\medskip

\noindent (b, equality case) We argue by induction on $r$. For $r=1$, we have
in fact $\ell_j(u)=\alpha_ju_1$, $\alpha_j\in\bC^*$, and
$I(\ell_1,\ldots,\ell_r)=\prod|\ell_j|^2$,
thus the coefficient ${1\over (p+r-1)!}={1\over p!}$ is reached if
and only if $p\le 1$. Now, assume $r\ge 2$ and the equality case solved
for dimension $r-1$. By rescaling and reordering the $\ell_j$,
we can always assume that $\ell_j(e_r)\ne 0$ (and hence
$\ell_j(e_r)=1$) for $q+1\le j\le p$, while $\ell_j(e_r)=0$ for
$1\le j\le q$ (we can possibly have $q=0$ here). Then we write
$\ell_j(u)=\ell'_j(u')$ for $1\le j\le q$ and
$\ell_j(u)=\ell'_j(u')+u_r$ for $q+1\le j\le p$. Therefore, if
$s_k(\ell'(u'))$ denotes the $k$-th elementary symmetric function in
$(\ell'_j(u')_{q+1\le j\le p}$, we find
$$
\eqalign{
I(\ell_1,\ldots,\ell_p)
&=\int_{\bS^{2r-1}}\prod_{j=1}^q|\ell'_j(u')|^2\,
\prod_{j=q+1}^p|\ell'_j(u')+u_r|^2\,d\mu(u)\cr
&=\int_{\bS^{2r-1}}\prod_{j=1}^q|\ell'_j(u')|^2
\bigg|\sum_{k=0}^{p-q}s_k(\ell'(u'))\,u_r^{p-q-k}\bigg|^2\,d\mu(u)\cr
&=\int_{\bS^{2r-1}}\prod_{j=1}^q|\ell'_j(u')|^2
\sum_{k=0}^{p-q}\big|s_k(\ell'(u'))\big|^2\,|u_r|^{2(p-q-k)}\,d\mu(u)\cr
&=\int_{\bS^{2r-3}}\prod_{j=1}^q|\ell'_j(u')|^2
\sum_{k=0}^{p-q}{(r-1)\,(k+r-2)!\,(p-q-k)!\over(p-q+r-1)!}\,
\big|s_k(\ell'(u'))\big|^2\,d\mu'(u')\cr
&\ge{(r-1)!\over(p+r-1)!}\,\prod_{j=1}^q|\ell'_j|^2
\prod_{j=q+1}^p(1+|\ell'_j|^2)\cr}
$$
by what we have just proved. In an equivalent way, we get
$$
\eqalign{
\int_{\bS^{2r-3}}&\prod_{j=1}^q|\ell'_j(u')|^2
\sum_{k=0}^{p-q}{(k+r-2)!\,(p-q-k)!\,(p+r-1)!\over(r-2)!\,(p-q+r-1)!}\,
\big|s_k(\ell'(u'))\big|^2\,d\mu'(u')\cr
&\ge\prod_{j=1}^q|\ell'_j|^2
\prod_{j=q+1}^p(1+|\ell'_j|^2)\cr}
$$
for all $0\le q\le p-1$ and all choices of the forms $\ell'_j\in(\bC^{r-1})^*$.
In general, we can rotate coordinates in such a way that $\ell_p(u)=u_r$
and $\ell'_p=0$, and we see that the above inequality holds
when $p$ is replaced by $p-1$, as soon as $q\le p-2$. Then the
corresponding coefficients $k=0$ for $p$, $p-1$ are
$$
{(p-q)!\,(p+r-1)!\over(p-q+r-1)!}>{(p-1-q)!\,(p-1+r-1)!\over(p-1-q+r-1)!},
$$
and since $s_0=1$, we infer that the inequality is strict. The only
possibility for the equality case is $q=p-1$, but then
$$
I(\ell_1,\ldots,\ell_p)=
\int_{\bS^{2r-1}}\prod_{j=1}^{p-1}|\ell'_j(u')|^2\,|u_r|^2\,d\mu(u)=
{r-1\over p+r-1}\int_{\bS^{2r-3}}\prod_{j=1}^{p-1}|\ell'_j(u')|^2\,d\mu'(u'),
$$
and we see that we must have equality in the case $(r-1,p-1)$.
By induction, we conclude that $p-1\le r-1$ and that the
$\ell_j(u)=\ell'_j(u')$ are orthogonal for $j\le p-1$, as desired.\qed

\claim 5.12. Remark|{\rm When $r=2$, our inequality (5.11) is equivalent
to the ``elementary'' inequality
$$
\prod_{j=1}^p(1+|a_j|^2)\le \sum_{k=0}^p k!\,(p-k)!\,|s_k|^2,\leqno(*)
$$
relating a polynomial $X^p-s_1X^{p-1}+\cdots+(-1)^ps_p$ and its complex
roots~$a_j$ (just consider~$\ell'_j(u')=a_ju_1$ and $\ell_j(u)=a_ju_1+u_2$
on $\bC^2$ to get this). It should be observed 
that $(*)$ is not optimal symptotically when $p\to+\infty\,$; in fact,
Landau's inequa\-lity [Land05] gives
$\prod\max(1,|a_j|)\le(\sum|s_k|^2)^{1/2}$, from which one
can easily derive that $\prod(1+|a_j|^2)\le 2^p\sum|s_k|^2$, which improves
$(*)$ as soon as $p\ge 7$ (observe that $2^7=128$ and
$k!(7-k)!\ge 3!\,4!=144$).
Our discussion of the equality case shows that inequality (5.8~(b))
is never sharp when $p>r$. It would be interesting, but probably
challenging, if~not impossible, to compute the optimal constant for 
all pairs $(r,p)$, $p>r$, since this is an optimization problem 
involving the distribution of a large number of points in projective 
space.}
\endclaim

\noindent We finally state one of the main consequences of these estimates
concerning the Chern curvature form of a hermitian holomorphic vector bundle.

\claim 5.13. Proposition|Let $T$, $E$ be complex vector spaces of
respective dimensions $\dim T=n$,\break $\dim E=r$. Assume that $E$ is
equipped with a hermitian structure $h$, and denote by $\mu$ the
unitary invariant probability measure $\mu$ on the unit
sphere bundle $S(E)=\{u\in E\,;|u|_h=1\}$ of~$E$.
\vskip2pt
\plainitem{\rm(a)} If $\ell_1,\ldots,\ell_k\in E^*$ and
$\theta_1,\ldots,\theta_{p-k}\ge_S 0$
are strongly semi-positive hermitian tensors in\break
$\Herm(T\otimes E)\simeq\Lambda_\bR^{1,1}T^*\otimes_\bR\Herm(E,E)$, then
$$
\eqalign{
\int_{u\in S(E)}|\ell_1(u)|^2\ldots|\ell_k(u)|^2\,
&\langle\theta_1(u),u\rangle_h\wedge\ldots
\wedge\langle\theta_{p-k}(u),u\rangle_h\,d\mu(u)\cr
&\plaincases{\displaystyle
\ge {(r-1)!\over (p+r-1)!}\,\bigg(\prod_{j=1}^k|\ell_j|^2\bigg)
\Tr_h\theta_1\wedge\ldots\wedge\Tr_h\theta_{p-k},\cr
\noalign{\vskip5pt}
\displaystyle
\le {p!\,(r-1)!\over (p+r-1)!}\,\bigg(\prod_{j=1}^k|\ell_j|^2\bigg)
\Tr_h\theta_1\wedge\ldots\wedge\Tr_h\theta_{p-k},\cr}\cr}
$$
as pointwise strong inequalities of $(p-k,p-k)$-forms.
\vskip2pt
\plainitem{\rm(b)} If $\theta\ge_G 0$ in $\Lambda_\bR^{1,1}T^*\otimes_\bR\Herm(E,E)$
and $\ell_j\in E^*$, then
$$
\int_{u\in S(E)}|\ell_1(u)|^2\ldots|\ell_k(u)|^2\,
\langle\theta(u),u\rangle_h^{p-k}\,d\mu(u)\le
{p!\,(r-1)!\over (p+r-1)!}\bigg(\prod_{j=1}^k|\ell_j|^2\bigg)(\Tr_h\theta)^{p-k}
$$
as a pointwise weak inequality of $(p-k,p-k)$-forms.
\vskip2pt
\noindent
In particular, the above inequalities apply when
$(E,h)$ is a hermitian holomorphic vector bundle
of rank $r$ on a complex $n$-dimensional manifold $X$, and one takes
$\theta_j=\Theta_{E,h}$ to be the curvature tensor of~$E$, so that
$\Tr_h\theta_j=c_1(E,h)$ is the first Chern form of $(E,h)$.
\endclaim

\plainproof. (a) The assumption $\theta_q\ge_S 0$ means that at every point
$x\in X$ we can write $\theta$ as
$$
\theta_q=\sum_{1\le j\le N_q}|\beta_{qj}\otimes\ell_{qj}|^2
\simeq \sum_{1\le j\le N_q}\ii\beta_{qj}\wedge\overline\beta_{qj}
\otimes\ell_{qj}\otimes\ell_{qj}^*,\quad\beta_{qj}\in T^*,~~\ell_{qj}\in E^*
$$
as an element of $\Lambda_\bR^{1,1}T^*\otimes_\bR\Herm(E,E)$, hence
$$
\langle\theta_q(u),u\rangle_h=\sum_{1\le j\le N_q}
\ii\beta_{qj}\wedge\overline\beta_{qj}\,|\ell_{qj}(u)|^2.
$$
Without loss of generality, we can assume $|\ell_{qj}|_{h^*}=1$.
Then
$$
\eqalign{
&|\ell_1(u)|^2\ldots|\ell_k(u)|^2\,\langle\theta_1(u),u\rangle_h\wedge
\ldots\wedge\langle\theta_{p-k}(u),u\rangle_h\cr
&\kern2mm{}=\sum_{j_1,\ldots,j_{p-k}}
\ii\beta_{1j_1}\wedge\overline\beta_{1j_1}\wedge\ldots\wedge
\ii\beta_{p-k\,j_{p-k}}\wedge\overline\beta_{p-k\,j_{p-k}}\,
\prod_{1\le s\le k}|\ell_s(u)|^2\prod_{1\le s\le p-k}|\ell_{sj_s}(u)|^2,\cr}
$$
and since $|\ell_{qj}|_{h^*}=1$, Lemma~5.8~(b) implies
$$
\eqalign{
\int_{u\in S(E)}&|\ell_1(u)|^2\ldots|\ell_k(u)|^2\,
\langle\theta_1(u),u\rangle_h\wedge\ldots\wedge
\langle\theta_{p-k}(u),u\rangle_h\,d\mu(u)\cr
&\ge {(r-1)!\over(p+r-1)!}~
\sum_{j_1,\ldots,j_{p-k}}\ii\beta_{1j_1}\wedge\overline\beta_{1j_1}
\wedge\ldots\wedge\ii\beta_{p-k\,j_{p-k}}\wedge\overline\beta_{p-k\,j_{p-k}}\,
\prod_{1\le s\le k}|\ell_s|^2\cr
&={(r-1)!\over(p+r-1)!}\,\bigg(\prod_{1\le j\le k}|\ell_j|^2\bigg)
\Tr_h\theta_1\wedge\ldots\wedge\Tr_h\theta_p,\cr}
$$
where $\ge$ is in the sense of the strong positivity of $(p,p)$-forms.
The upper bound is obtained by the same argument, via 5.8~(a).
\medskip
\noindent
(b) By the definition of weak positivity of forms, it is enough to
show the inequality in restriction to every $(p-k)$-dimensional subspace
$T'\subset T$. Without loss of generality, we can assume that $\dim T=p-k$
(and then take $T'=T$), that $|\ell_j|=1$, and also that
$\theta>_G0$ (otherwise take a positive definite form
$\eta\in\Lambda_\bR^{1,1}T^*$,
replace $\theta$ with $\theta_\varepsilon=\theta+\varepsilon\,\eta\otimes h$,
and let $\varepsilon$ tend to $0$). For any $u\in S(E)$, let
$$
0\le\lambda_1(u)\le\cdots\le\lambda_{p-k}(u)
$$
be the eigenvalues of the hermitian form
$q_u(\bu)=\langle \theta(u),u\rangle$ on $T$ with respect to
$$
\omega=\Tr_h\theta=\sum_{j=1}^r\langle \theta(e_j),e_j\rangle\in
\Herm(T),\quad\omega>0,
$$
$(e_j)$ being any orthonormal frame of $E$. We have to show that
$$
\int_{u\in S(E)}
|\ell_1(u)|^2\ldots|\ell_k(u)|^2\,\lambda_1(u)\cdots\lambda_{p-k}(u)\,d\mu(u)\le
{p!\,(r-1)!\over(p+r-1)!}.
$$
However, the inequality between geometric and arithmetic means implies
$$
\lambda_1(u)\cdots\lambda_p(u)\le
\bigg({1\over p-k}\sum_{j=1}^{p-k}\lambda_j(u)\bigg)^p,
$$
thus, putting $Q(u)={1\over p-k}\langle\Tr_\omega\theta(u),u\rangle$,
$Q\in\Herm(E)$, it is enough to prove that
$$
\int_{u\in S(E)}|\ell_1(u)|^2\ldots|\ell_k(u)|^2\,Q(u)^{p-k}\,d\mu(u)
\le {p!\,(r-1)!\over(p+r-1)!}.\leqno(5.14)
$$
Our assumption $\theta>_G0$ implies
$Q(u)=\sum_{1\le j\le r} c_j|\ell'_{qj}(u)|^2$
for some $c_j>0$ and some orthonormal basis $(\ell'_{qj})_{1\le j\le r}$ of
$E^*$, and
$$
\sum_{j=1}^r c_j=\Tr_hQ={1\over p-k}\Tr_h(\Tr_\omega\theta)
={1\over p-k}\Tr_\omega(\Tr_h\theta)={1\over p-k}\Tr_\omega(\omega)=1.
$$
Inequality (5.14) is a consequence of Lemma~5.8~(a), by
Newton's multinomial expansion.\qed

\claim 5.15. Remark|{\rm For $p=1$, the inequalities of Proposition~5.13 are
identities, and no semi-positivity assumption is needed in that case.
This can be seen directly from the fact that we have
$$
\int_{u\in S(E)}Q(u)\,d\mu(u)={1\over r}\,\Tr Q
$$
for every hermitian quadratic form $Q$ on $E$. However,
when $p\ge 2$, inequality 5.13~(a) does not hold under the assumption that
$E\ge_G0$ (or even that $E$ is dual Nakano semi-positive, i.e.\  $E^*$ 
Nakano semi-negative).
Let us take for instance $E=T_{\bP^n}\otimes\cO(-1)$. It is well known that
$E$ is isomorphic to the tautological quotient vector bundle 
$\bC^{n+1}/\cO(-1)$ over $\bP^n$, and that its curvature tensor form for the
Fubini-Study metric is given by
$$
\Theta_{E,h}(\xi\otimes u,\xi\otimes u)=
|\langle\xi,u\rangle|^2\ge 0
$$
(where $v$ is identified which a tangent vector via the choice of a unit
element $e\in\cO(-1)$). Then $\det E=\cO(1)$ and thus
$c_1(E,h)=\omega_\FS>0$, although $\langle\Theta_{E,h}(u),u\rangle_h^p=0$
for all $p\ge 2$, as one can easily check.}
\endclaim

\plainsection{6. On the curvature of orbifold tangent bundles}

\plainsubsection 6.A. Evaluation of the orbifold curvature tensor|

The main qualitative result is summarized in the following statement.

\claim 6.1.~Proposition|Let $X$ be a projective variety, $A$ an ample line
bundle, and $(X,V,D)$ an orbifold directed structure where
$D=\sum_{1\le j\le N}(1-{1\over\rho_j})\Delta_j$ is a 
normal crossing divisor transverse to~$V$ in~$X$. Let $d_j$ be the
infimum of numbers
$\lambda\in\bR_+$ such that $\lambda A-\Delta_j$ is~nef, and $\gamma_V$
be the infimum of numbers $\gamma\ge 0$ such that $\gamma\,
\Theta_{A,h_A}\otimes\Id_V-\Theta_{V,h_V}\ge_G0$
for suitable smooth hermitian metrics $h_V$ on $V$. Then for every
$\gamma>\gamma_{V,D}:=\max(\max_j(d_j/\rho_j),\gamma_V)$,
the~orbifold vector bundle $V\langle D\rangle$
possesses a hermitian metric $h_{V\langle D\rangle,\gamma,\varepsilon}$ such that
\vskip2pt
\plainitem{\rm(a)} $h_{V\langle D\rangle,\gamma,\varepsilon}$ is smooth on $X\ssm|D|,$
\vskip2pt
\plainitem{\rm(b)} $h_{V\langle D\rangle,\gamma,\varepsilon}$ has the appropriate
orbifold singularities along $D,$
\vskip2pt
\plainitem{\rm(c)} we have $\gamma\,\Theta_{A,h_A}\otimes\Id-
\Theta_{V\langle D\rangle,h_{V\langle D\rangle,\gamma,\varepsilon}}\ge_G0$ 
on $X\ssm|D|$.
\vskip2pt  
\endclaim

\plainproof. Let $h_A$ be a metric on $A$ such that $\Theta_{A,h_A}>0$,
written locally as $h_A=e^{-\psi}$, and take
\hbox{$\gamma>\max(\max_j(d_j/\rho_j),\gamma_V)$}.
Consider the tautological sections
$\sigma_j\in H^0(X,\cO_X(\Delta_j))$ defining $\Delta_j=\sigma_j^{-1}(0)$,
and let $h_V$, $h_j$ be smooth hermitian metrics on
$V$ and $\cO_X(\Delta_j)$ such that
$$
\leqalignno{
&\gamma\,\Theta_{A,h_A}\otimes\Id_V-\Theta_{V,h_V}>_G0,&(6.2_0)\cr
&\gamma\,\Theta_{A,h_A}-{1\over\rho_j}\Theta_{\cO_X(\Delta_j),h_j}>0,\quad
\forall j=1,\ldots,N,
&(6.2_j)\cr}
$$
as is possible by our choice of the constants $d_j$ and $\gamma$.
Finally, denote by $\nabla_j$ the associated Chern connection on
$\cO_X(\Delta_j)$. If we write $h_j=e^{-\varphi_j}$ in some local
trivialization, then $\nabla_j\sigma_j=\nabla_j^{1,0}\sigma_j=
\partial\sigma_j-\sigma_j\partial\varphi_j$. 
Take $\omega_A=\Theta_{A,h_A}$ as the K\"ahler metric on $X$. We have
$$
\ii\ddbar|\sigma_j|_{h_j}^{2/\rho_j}=
{1\over\rho_j^2}\,|\sigma_j|_{h_j}^{-2+2/\rho_j}\,
i\langle\nabla_j\sigma_j,\nabla_j\sigma_j\rangle_{h_j}
-{1\over\rho_j}|\sigma_j|_{h_j}^{2/\rho_j}\,\ii\ddbar\varphi_j,
$$
hence there exists $\delta>0$ small such that the metric
$h_{A,\delta}=h_A\exp(-\delta\sum_j|\sigma_j|_{h_j}^{2/\rho_j})$
of weight $\psi_\delta=\psi+\delta\sum_j|\sigma_j|_{h_j}^{2/\rho_j}$ satisfies
$$
\ii\ddbar\psi_\delta(\xi,\xi)=
|\xi|^2_{\omega_A}+\delta\ii\ddbar
\sum_j\kern-0.5pt|\sigma_j|_{h_j}^{2/\rho_j}(\xi,\xi)\ge
(1-C\delta)|\xi|^2_{\omega_A}+\delta\sum_j{1\over\rho_j^2}|
\sigma_j|_{h_j}^{-2+2/\rho_j}\,|\nabla_j\sigma_j(\xi)|_{h_j}^2.
$$
We can consider $\omega_{A,\delta}=\Theta_{A,h_{A,\delta}}=\ii\ddbar\psi_\delta$
as an orbifold K\"ahler metric, that is ``smooth'' from the point of view
of the orbifold structure. Let us explain the more precise meaning of this
``orbifold smoothness'' assumption. In fact, there exists a ramified cover
$g_Y:Y\to X$ such that $g^*\sigma_j=w_j^{m_j}$ for some
local coordinate $w_j$ on $Y$, with arbitrary high multiplicity $m_j\in\bN^*$
along $\smash{g_Y^{-1}}(\Delta_j)=\{w_j=0\}$.
Then $g_Y^*h_{A,\delta}=g_Y^*h_A\exp(-\delta\sum_j
\smash{|w_j|^{2m_j/\rho_j})}$
can be taken in any regularity class $C^p$, $p\in \bN^*$, by taking
$m_j\ge p\,\rho_j$. Therefore, by pulling-back our calculations to $Y$, we
would actually get forms of high regularity on~$Y$. Of course, if we
compute an integral over $X$, pulling-back forms to $Y$ multiplies the
integral by the degree of $g_Y$, and it suffices to divide by that degree
to recover the integral over~$X$. For $\delta>0$ sufficiently small, our
positivity conditions $(6.2_j)$ can
be turned into the stronger form
$$
\leqalignno{
&\gamma\,\ii\ddbar\psi_\delta(\xi,\xi)\,|u|^2-\widetilde\Theta_{V,h_V}(\xi\otimes u)\ge 
c\bigg(|\xi|_{\omega_A}^2+\sum_j|\sigma_j|_{h_j}^{-2+2/\rho_j}\,|\nabla_j\sigma_j(\xi)|_{h_j}^2
\bigg),&(6.3_0)\cr
&\gamma\,\ii\ddbar\psi_\delta(\xi,\xi)-{1\over\rho_j}\,\ii\ddbar\varphi_j(\xi,\xi)\ge 
c\bigg(|\xi|_{\omega_A}^2+\sum_j|\sigma_j|_{h_j}^{-2+2/\rho_j}\,|\nabla_j\sigma_j(\xi)|_{h_j}^2
\bigg),&(6.3_j)\cr}
$$
for some constant $c>0$ and all $\xi\in T_X$, $u\in V$ (observe that 
the right hand side can in fact be seen as a positive definite hermitian 
form with respect to the orbifold coordinates, we just exploit the fact that
$A$ remains ample when viewed as a line bundle on the orbifold structure).
We~are going to estimate the curvature of the orbifold metric
$h_{V\langle D\rangle,\varepsilon}$ on $V\langle D\rangle$ defined by
$$
\Vert u\Vert_{h_{V\langle D\rangle,\varepsilon}}^2=|u|_{h_V}^2+\sum_j\varepsilon_j\,
|\sigma_j|_{h_j}^{-2(1-1/\rho_j)}\,|\nabla_j\sigma_j(u)|_{h_j}^2,
\quad \varepsilon_j\ll 1.
\leqno(6.4)
$$
Again, this metric can be seen as orbifold smooth (in the sense that
the metric $g_Y^*h_{V\langle D\rangle,\varepsilon}$ on $g_Y^*(V\langle D\rangle)$
may be taken of arbitrary high regularity; in case $\rho_j=\infty$, it is
actually a smooth metric on the logarithmic bundle). Since
$$
\ii\ddbar\Vert u\Vert_{h_{V\langle D\rangle},\varepsilon}^2=
\ii\langle\nabla u,\nabla u\rangle_{h_{V\langle D\rangle},\varepsilon}-2\pi\,
\langle\Theta_{V\langle D\rangle,h_{V\langle D\rangle,\varepsilon}}(u),u
\rangle_{h_{V\langle D\rangle,\varepsilon}}
$$
where $\nabla u=du+\Gamma(dz)\cdot u$ is the Chern connection of
$(V\langle D\rangle,h_{V\langle D\rangle,\varepsilon})$, what we need to prove is
that on the total space of $V$ over $X\ssm|D|$, the $(1,1)$-form
$$
V\ni (z,u)\mapsto
\ii\ddbar\Vert u\Vert_{h_{V\langle D\rangle,\varepsilon}}^2+\gamma\,
\ii\ddbar\psi_\delta\,\Vert u\Vert_{h_{V\langle D\rangle,\varepsilon}},
$$
is non negative. For this, we calculate the associated hermitian
quadratic form on $T_V$
$$
Q_{V\langle D\rangle,\gamma,\varepsilon}(z,u)(\xi,\eta),\quad
(\xi,\eta)\in T_{V,(z,u)},\quad
\xi=\sum_{\ell=1}^n\xi_\ell{\partial\over\partial z_\ell},\quad
\eta=\sum_{\lambda=1}^r\eta_\lambda{\partial\over\partial u_\lambda},
\leqno(6.5)
$$
and observe that the curvature tensor is obtained by taking the
restriction to the ``parallel'' directions $\nabla u=0$, that is,
by substituting $du=-\Gamma(dz)\cdot u$, i.e.\
$\eta=-\Gamma(\xi)\cdot u$. Let us fix an arbitrary point
$z_0\in X\ssm|D|$.
We take local holomorphic coordinates $(z_1,\ldots,z_n)$ centered at~$z_0$,
and let $(e_1,\ldots,e_r)$ be a local holomorphic frame of $V$ such that
$$
\langle e_\lambda,e_\mu\rangle_{h_V}=\delta_{\lambda\mu}+
\sum_{\ell,m,\lambda,\mu}c_{\ell m\lambda\mu}\,z_\ell\overline z_m+O(|z|^3),
$$
where the ${\ii\over 2\pi}c_{\ell m\lambda\mu}$ are the coefficients of
$-\Theta_{V,h_V}$.
Let us write $u=\sum_{\lambda=1}^ru_\lambda e_\lambda$ and denote by $\langle
u,v\rangle=\sum_{1\le\lambda\le r}u_\lambda\overline v_\lambda$ the standard
hermitian form, $|u|$ the associated norm. We~find
$$
\leqalignno{
&\Vert u\Vert_{h_{V\langle D\rangle,\varepsilon}}^2
=|u|^2+\sum_{\ell,m,\lambda,\mu}c_{\ell m\lambda\mu}\,z_\ell\overline z_m
u_\lambda\overline u_\mu+O(|z|^3)\cr  
&\kern50pt{}
+\sum_j\varepsilon_j\,\big(|\sigma_j|^2e^{-\varphi_j}\big)^{-1+1/\rho_j}\,
\big|\partial\sigma_j(u)-\sigma_j\partial\varphi_j(u)\big|^2e^{-\varphi_j},
&(6.4_0)\cr}
$$
since $\dbar\sigma_j=0$. In order to simplify the calculation, we set formally
$$
\left\{\plainmatrix{
&\tilde\sigma_j=\sigma_j^{1/\rho_j},\hfill
&\tilde\varepsilon_j=\rho_j^2\varepsilon_j,\hfill
&\tilde\varphi_j=\rho_j^{-1}\varphi_j,\hfill
&&\hbox{if $\rho_j<\infty$,}\hfill\cr
\noalign{\vskip6pt}
&\tilde\sigma_j=\log\sigma_j,\hfill
&\tilde\varepsilon_j=\varepsilon_j,\hfill
&\tilde\varphi_j=\varphi_j,\hfill
&&\hbox{if $\rho_j=\infty$.}\hfill\cr}\right.
\leqno(6.6)
$$
Respectively to the non logarithmic and logarithmic situations, we then
get the more tractable expression
$$
\leqalignno{\strut\kern25pt
&\Vert u\Vert_{h_{V\langle D\rangle,\varepsilon}}^2
\,{=}\,|u|^2+{}\kern-4pt
\sum_{\ell,m,\lambda,\mu}\kern-4pt c_{\ell m\lambda\mu}\,z_\ell\overline z_m
u_\lambda\overline u_\mu\,{+}\,O(|z|^3)
\,{+}\sum_j\tilde\varepsilon_j\,
\big|\partial\tilde\sigma_j(u)\,{-}\,\sigma_j\partial\tilde\varphi_j(u)\big|^2\,
e^{-\tilde\varphi_j},&(6.7)\cr
\noalign{\vskip6pt}
&\kern8.5pt\Vert u\Vert_{h_{V\langle D\rangle,\varepsilon}}^2
\,{=}\,|u|^2+{}\kern-4pt
\sum_{\ell,m,\lambda,\mu}\kern-4pt c_{\ell m\lambda\mu}\,z_\ell\overline z_m
u_\lambda\overline u_\mu\,{+}\,O(|z|^3)
\,{+}\sum_j\tilde\varepsilon_j\,
\big|\partial\tilde\sigma_j(u)-\partial\tilde\varphi_j(u)\big|^2.
&(6.7_\infty)\cr}
$$
More importantly, the poles have disappeared -- a fact reflecting the
orbifold smoothness of the metric.
In what follows, for the sake of simplicity, we remove the tildes in the
notation, and conduct the calculation only in the non logarithmic
situation $(\rho_j<\infty)$, since the logarithmic case
can be recovered by taking $\rho_j$ very large; this actually amounts to using
a ramified change of variable $\tilde z'_\ell=\smash{z_\ell^{1/\rho_\ell}}$
in suitable coordinates, allowing us in this way to take $\rho_j=1$ 
in~$(6.4_0)$. Also, our later calculations will be done by adding the orbifold
divsior components one by one. This essentially
reduces the situation to the case where  $D=(1-{1\over\rho})\Delta$ only
has one component, and the notation becomes much lighter. Therefore, we drop
the indices $j$ and the summations $\sum_j$, and consider the simple
situation where the metric is given by
$$
\leqalignno{
\Vert u\Vert_{h_{V\langle D\rangle,\varepsilon}}^2
&=|u|^2+\sum_{\ell,m,\lambda,\mu}c_{\ell m\lambda\mu}\,z_\ell\overline z_m
u_\lambda\overline u_\mu+O(|z|^3)+\varepsilon\,|\partial\sigma(u)
-\sigma\,\partial\varphi(u)|^2&(6.8)\cr
\noalign{\vskip5pt}
\langle\!\langle u,v\rangle\!\rangle_{h_{V\langle D\rangle,\varepsilon}}^2
&=\langle u,v\rangle^2+
\sum_{\ell,m,\lambda,\mu}c_{\ell m\lambda\mu}\,z_\ell\overline z_m
u_\lambda\overline v_\mu+O(|z|^3)&(6.8')\cr
&\kern41pt{}+\varepsilon\,\big(\partial\sigma(u)-\sigma\,\partial\varphi(u)\big)
\big(\,\overline{\partial\sigma(v)-\sigma\,\partial\varphi(v)}\,\big)
\,e^{-\varphi}.\cr}
$$
We also take a holomorphic trivialization of the line bundle $\cO_X(\Delta)$
so that the associated weight $\varphi$ satisfies
$\varphi(z)=\sum_{\ell,m}\alpha_{\ell m}\,z_\ell\overline z_m+O(|z|^3)$ near
$z_0=0$. Then
$$
\partial\varphi=\sum_{\ell,m}\alpha_{\ell m}\,\overline z_m dz_\ell+O(|z|^2),
\quad
\dbar\varphi=\sum_{\ell,m}\alpha_{\ell m}\,z_\ell\,d\overline z_m+O(|z|^2).
$$
At the point $z=z_0$, we have $\partial\varphi(z_0)=\partial\varphi(z_0)=0$,
$\nabla\sigma=\partial\sigma$, and our metric admits the expression
$$
\Vert u\Vert_{h_{V\langle D\rangle,\varepsilon}}^2=
|u|^2+\varepsilon\,|\partial\sigma(u)|^2,\quad
\langle\!\langle u,v\rangle\!\rangle_{h_{V\langle D\rangle,\varepsilon}}=
\langle u,v\rangle+\varepsilon\,\partial\sigma(u)\,
\overline{\partial\sigma(v)}.
\leqno(6.9)
$$
Let $u,v$ be arbitrary local holomorphic sections of $V$, and denote by
$\nabla_\xi$ the Chern covariant differentiation
of $(V\langle D\rangle,h_{V\langle D\rangle,\varepsilon})$ in the direction
$\xi\in T_X$.
By polarizing the quadratic form
$\Vert u\Vert_{h_{V\langle D\rangle,\varepsilon}}^2$ into a hermitian inner product
$\partial_\xi\langle\!\langle u,v
\rangle\!\rangle_{h_{V\langle D\rangle,\varepsilon}}$
and setting $\nabla_\xi u=
\nabla^{1,0}_\xi u=\partial_\xi u+\Gamma(\xi)\cdot u$,
a~differentiation of $(6.8')$ at $z=z_0$ yields
$$
\eqalign{
\partial_\xi\langle\!\langle u,v
\rangle\!\rangle_{h_{V\langle D\rangle,\varepsilon}}
={}&\langle\nabla_\xi u,v\rangle+\varepsilon\,
\partial\sigma(\nabla_\xi u)\,\overline{\partial\sigma(v)}
\cr
={}&
\langle\partial_\xi u\,,\,v\rangle+\varepsilon\;
\partial\sigma(\partial_\xi u)\;\overline{\partial\sigma(v)}
+\varepsilon\,\partial^2\sigma(\xi,u)\,\overline{\partial\sigma(v)}
-\varepsilon\,
\partial\sigma(u)\,\overline \sigma\,\ddbar\varphi(\xi,v),\cr}
$$
where $\partial^2\sigma(\xi,u):=\sum_\lambda\partial_\xi\big(\partial\sigma
(e_\lambda)\big)\,u_\lambda$ is viewed as an element of
$(T_X^*\otimes V^*)_{z_0}$ and $\ddbar\varphi$ as a hermitian form on $T_X$,
operating on $T_X\otimes\overline V\subset T_X\otimes \overline T_X$. In fact,
$u\mapsto\partial\sigma(u)$ and $(\xi,u)\mapsto\partial^2\sigma(\xi,u)$
can be intrinsically defined as $\nabla^{1,0}\sigma_{|V}$ and
$\smash{\nabla_{V^*\otimes \cO(\Delta)}^{1,0}}(\nabla^{1,0}\sigma_{|V})$ at~$z_0$,
and we will denote them by $\nabla\sigma$ and $\nabla^2\sigma$.
In this setting, a subtraction of the last two lines in our equalities
shows that the $(1,0)$-form $\Gamma$ of
the connection of $(V\langle D\rangle,h_{V\langle D\rangle})$
is given at~$z_0$ by the formula
$$
\langle \Gamma(\xi)\cdot u,v\rangle
+\varepsilon\,
\nabla\sigma(\Gamma(\xi)\cdot u)\,\overline{\nabla\sigma(v)}
=\varepsilon\,
\nabla^2\sigma(\xi,u)\,\overline{\nabla\sigma(v)}
-\varepsilon\,\nabla\sigma(u)\;\overline\sigma\,\ddbar\varphi(\xi,v).
\leqno(6.10)
$$
This equality if valid pointwise for any $u,v\in V_{z_0}$.
As a consequence
$$
\Gamma(\xi)\cdot u+\varepsilon\,
\nabla\sigma(\Gamma(\xi)\cdot u)\,(\nabla\sigma)^*
=\varepsilon\,
\nabla^2\sigma(\xi,u)\,(\nabla\sigma)^*
-\varepsilon\,\nabla\sigma(u)\;\overline\sigma\,
(\ddbar\varphi(\bu,\xi))^*
\leqno(6.11)
$$
where $\alpha^*\in V$ is the dual vector to a $1$-form
$\alpha\in V^*$, such that
$\langle\alpha^*,\bu\rangle_{h_V}=\overline\alpha$. 
The special choice $v=\Gamma(\xi)\cdot u$ yields a (non negative) real
value in the left hand side of~$(6.10)$, and by taking the real
part of the right hand side, we obtain
$$
\leqalignno{
|\Gamma(\xi)\cdot u|^2
&+\varepsilon\,
\big|\nabla\sigma(\Gamma(\xi)\cdot u)\big|^2\cr
={}&\varepsilon\,\Re\big(
\nabla^2\sigma(\xi,u)\,
\overline{\nabla\sigma(\Gamma(\xi)\cdot u)}\,\big)-\varepsilon\,
\Re\big(
\nabla\sigma(u)\,\overline\sigma\;\ddbar\varphi(\xi,\Gamma(\xi)
\cdot u)\big).
&(6.12_0)\cr
\noalign{\vskip5pt}
}
$$
Also, by applying $\nabla\sigma$ to (6.11), we obtain
$$
\eqalign{
\nabla\sigma(\Gamma(\xi)\cdot u)&+\varepsilon\,
\nabla\sigma(\Gamma(\xi)\cdot u)\,
\langle\nabla\sigma,\nabla\sigma\rangle\cr
&=\varepsilon\,\nabla^2\sigma(\xi,u)\,
\langle\nabla\sigma,\nabla\sigma\rangle
-\varepsilon\,\nabla\sigma(u)\;\overline\sigma\,
\langle\nabla\sigma,\ddbar\varphi(\bu,\xi)\rangle,
\cr}
$$
hence
$$
\leqalignno{
\nabla\sigma(\Gamma(\xi)\cdot u)&=
{\varepsilon\over 1+\varepsilon\,|\nabla\sigma|^2}\big(
\nabla^2\sigma(\xi,u)\,|\nabla\sigma|^2-\nabla\sigma(u)\;\overline\sigma\,
\langle\nabla\sigma,\ddbar\varphi(\bu,\xi)\rangle\big).&(6.12_1)\cr}
$$
As $2\pi\,\Theta_{A,h_A}=\ii\ddbar\psi_\delta$, we infer
by a brute force calculation from (6.8) that
$$
\leqalignno{
Q_{V\langle D\rangle,\gamma,\varepsilon}(z,u)(\xi,\eta)={}
&\ddbar\Vert u\Vert_{h_{V\langle D\rangle,\varepsilon}}^2\cdot(\xi,\eta)
+\gamma\,\ddbar\psi_\delta(\xi,\xi)\,
\Vert u\Vert_{h_{V\langle D\rangle,\varepsilon}}^2\cr
\noalign{\vskip6pt}
\kern40pt={}&\gamma\,\ddbar\psi_\delta(\xi,\xi)\,|u|^2
+\sum_{\ell,m,\lambda,\mu}c_{\ell m\lambda\mu}\,\xi_\ell\overline\xi_m\,
u_\lambda\overline u_\mu
&(6.13_1)\cr
&+\varepsilon\,\big(\gamma\,\ddbar\psi_\delta(\xi,\xi)-
\ddbar\varphi(\xi,\xi)\big)\,|\nabla\sigma(u)|^2\kern80pt
&(6.13_2)\cr
&+|\eta|^2+\varepsilon\,
\big|\nabla\sigma(\eta)+\nabla^2\sigma(\xi,u)\big|^2
&(6.13_3)\cr
&\kern6pt{}-2\varepsilon\,\Re\big(
\nabla\sigma(u)\,\overline\sigma\,\ddbar\varphi(\xi,\eta)\big)
&(6.13_4)\cr
&\kern6pt{}-2\varepsilon\,\Re\big(
\nabla\sigma(u)\,\ddbar\varphi(\xi,u)\,
\overline{\nabla\sigma(\xi)}\,\big)
&(6.13_5)\cr
&\kern6pt{}-2\varepsilon\,\Re\big(
\nabla\sigma(u)\,\overline\sigma\,\ddbar{\kern1pt}^2
\varphi(\xi,\xi,u)\big)
&(6.13_6)\cr
&\kern6pt{}+\varepsilon\,|\sigma|^2\,|\ddbar\varphi(u,\xi)|^2,
&(6.13_7)\cr}
$$
where we identify a $(1,1)$-form such as $\ddbar\varphi$ with a hermitian
form, and take $\eta=-\Gamma(\xi)\cdot u$.
The second term in $(6.13_2)$ is obtained by differentiating
$\varepsilon\,|\nabla\sigma(u)|^2$, while
$(6.13_3)$, $(6.13_4)$ and $(6.13_5)$ actually come from the differentiation of
the term $...\Re(...)$ in $(6.8)$. By our assumptions $(6.3_j)$, the first
two terms $(6.13_1)$, $(6.13_2)$
are positive in the sense of Griffiths, and such that
$$
\eqalign{
&(6.13_1)\ge c\,(|\xi|^2+|\nabla\sigma(\xi)|^2)\,|u|^2,\cr
\noalign{\vskip4pt}  
&(6.13_2)\ge c\,\varepsilon\,(|\xi|^2+|\nabla\sigma(\xi)|^2)\,
|\nabla\sigma(u)|^2,\quad c>0.\cr}
$$
(Here the term $|\nabla\sigma(\xi)|^2$ is significant, because
we will later replace $\sigma$ by $\sigma^{1/\rho}$ in the orbifold case,
and then  $\nabla\sigma^{1/\rho}(\xi)$ is unbounded with respect to $|\xi|)$.
The third term $(6.13_3)$ is semi-positive.
We claim that the terms $(6.13_{4,5,6,7})$ are negligible for
$\varepsilon\ll 1$, in the sense that
$Q_{V\langle D\rangle,\gamma,\varepsilon}(z,u)(\xi,\eta)$ is comprised
between $(1\pm\delta)\sum_{j=1,2,3}(6.13_j)$, with $\delta>0$ as small as
we want when $\varepsilon\le\varepsilon_0(\delta)$. In fact, since
$\ddbar\varphi$ is smooth, there exists $C>0$ such that
$$
\eqalign{
\big|(6.13_4)\big|
&\le C\,\varepsilon|\sigma|\,|\nabla\sigma(u)|\,|\xi|\,|\eta|\cr
&\le \varepsilon^{3/2}|\xi|^2\,|\nabla\sigma(u)|^2+
C^2\,\varepsilon^{1/2}\,|\sigma|^2\,|\eta|^2
\ll(6.13_2)+(6.13_3).\cr}
$$
Similarly
$$
\eqalign{
\big|(6.13_5)\big|
&\le C\,\varepsilon\,
|\xi|\,|u|\,|\nabla\sigma(\xi)|\,|\nabla\sigma(u)|\cr
&\le C\,\varepsilon^{3/2}\,|\xi|^2\,|\nabla\sigma(u)|^2
+C\,\varepsilon^{1/2}\,|\nabla\sigma(\xi)|^2\,|u|^2
\ll (6.13_1)+(6.13_2).\cr}
$$
The last two terms $(6.13_{6,7})$ are even easier, since
$$
\eqalign{
\big|(6.13_6)\big|
&\le C\,\varepsilon\,|\sigma|\,|\xi|^2\,|u|\,|\nabla\sigma(u)|
\le \varepsilon^{1/2}|\xi|^2\,|u|^2+
C^2\,\varepsilon^{3/2}\,|\sigma|^2\,|\xi|^2\,|\nabla\sigma(u)|^2\cr
&\ll(6.13_1)+(6.13_2),\cr
\noalign{\vskip4pt}
\big|(6.13_7)\big|
&\le C\,\varepsilon\,|\xi|^2|u|^2\ll (6.13_1).\cr}
$$
Finally, by replacing $\eta$ with $-\Gamma(\xi)\cdot u$
and using $(6.12_{0,1})$, we find
$$
\leqalignno{
(6.13_3)+(6.13_4)&=\,\big|\Gamma(\xi)\cdot u\big|^2\cr
&\quad{}+\varepsilon\,
\big|\nabla\sigma(\Gamma(\xi)\,{\cdot}\,u)
-\nabla^2\sigma(\xi,u)\big|^2
+2\varepsilon\,\Re\big(
\nabla\sigma(u)\,\overline\sigma\,
\ddbar\varphi(\xi,\Gamma(\xi)\cdot u)\big)\cr
&=(6.12_0)+\varepsilon\,\big|
\nabla^2\sigma(\xi,u)\big|^2-2\varepsilon\,\Re\big(
\,\overline{\nabla^2\sigma(\xi,u)}\,
\nabla\sigma(\Gamma(\xi)\cdot u)\big)\cr
&\quad{}+2\varepsilon\,\Re\big(
\nabla\sigma(u)\,\overline\sigma\,
\ddbar\varphi(\xi,\Gamma(\xi)\cdot u)\big)\cr
&=\varepsilon\,\big|\nabla^2\sigma(\xi,u)\big|^2
-\varepsilon\,\Re\big(\,
\overline{\nabla^2\sigma(\xi,u)}\,
\nabla\sigma(\Gamma(\xi)\cdot u)\big)\cr
&\quad{}+\varepsilon\,\Re\big(
\nabla\sigma(u)\,\overline\sigma\,
\ddbar\varphi(\xi,\Gamma(\xi)\cdot u)\big).\cr
&={\varepsilon\over 1+\varepsilon|\nabla\sigma|^2}\,
\big|\nabla^2\sigma(\xi,u)\big|^2&(6.14_1)\cr
&\quad{}+{\varepsilon\over 1+\varepsilon\,|\nabla\sigma|^2}\,\Re\Big(
\overline{\nabla^2\sigma(\xi,u)}\,
\varepsilon\,\nabla\sigma(u)\,\overline\sigma\,
\langle\nabla\sigma,\ddbar\varphi(\bu,\xi)\rangle\Big)&(6.14_2)\cr
&\quad{}+\varepsilon\,\Re\big(
\nabla\sigma(u)\,\overline\sigma\,
\ddbar\varphi(\xi,\Gamma(\xi)\cdot u)\big).&(6.14_3)\cr}
$$
The term $(6.14_3)$ equals ${1\over 2}(6.13_4)$, thus it is negligible, and
the term $(6.14_2)$ admits an obvious bound
$$
\eqalign{
(6.14_2)&\le{\varepsilon\over 1+\varepsilon|\nabla\sigma|^2}
\Big(\varepsilon^{1/2}\,|\nabla^2\sigma(\xi,u)|^2+\varepsilon^{3/2}\,
|\sigma|^2\,|\nabla\sigma|^2\,|\nabla\sigma(u)|^2\,|\xi|^2\Big)\cr
&\le \varepsilon^{1/2}(6.14_1)+
\varepsilon^{3/2}|\sigma|^2\,|\nabla\sigma(u)|^2\,|\xi|^2\ll
(6.14_1)+(6.13_2).
\cr}
$$
By collecting all non negligible terms $(6.13_{1,2})$
and $(6.14_1)$, we obtain a curvature form
$$
\eqalign{
Q_{V\langle D\rangle,\gamma,\varepsilon}(z)(\xi\otimes u)
\simeq{}&\gamma\,\ddbar\psi_\delta(\xi,\xi)\,|u|^2
+\sum_{\ell,m,\lambda,\mu}c_{\ell m\lambda\mu}\,\xi_\ell\overline\xi_m\,
u_\lambda\overline u_\mu\cr 
&+\varepsilon\,
\big(\gamma\,\ddbar\psi_\delta(\xi,\xi)-\ddbar\varphi(\xi,\xi)\big)
\,|\nabla\sigma(u)|^2
+{\varepsilon\over 1+\varepsilon|\nabla\sigma|^2}\,
\Big|\nabla^2\sigma(\xi,u)\Big|^2.\cr}
$$
At this point, we come back to
the orbifold situation, and thus replace $\sigma$ by $\sigma^{1/\rho}$,
$\varphi$ by $\rho^{-1}\varphi$ and $\varepsilon$ by
$\rho^2\,\varepsilon$. This gives the curvature estimate
$$
\leqalignno{\kern10pt
Q_{V\langle D\rangle,\gamma,\varepsilon}(z)(\xi\otimes u)\simeq{}
&\gamma\,\ddbar\psi_\delta(\xi,\xi)\,|u|^2
+\sum_{\ell,m,\lambda,\mu}c_{\ell m\lambda\mu}\,\xi_\ell\overline\xi_m\,
u_\lambda\overline u_\mu\cr
&+\varepsilon\,|\sigma|^{-2+2/\rho}\,
\big(\gamma\,\ddbar\psi_\delta(\xi,\xi)-\rho^{-1}\,\ddbar\varphi(\xi,\xi)\big)
\,|\nabla\sigma(u)|^2&(6.15)\cr
&+{\varepsilon\,|\sigma|^{-2+2/\rho}\over
1+\varepsilon\,|\sigma|^{-2+2/\rho}\,
|\nabla\sigma|^2}\,\big|\,\nabla^2\sigma(\xi,u)
-(1-1/\rho)\,\sigma^{-1}\,\nabla\sigma(\xi)\nabla\sigma(u)\big|^2,\cr}
$$
In the general situation $D=\sum_{1\le j\le N}(1-1/\rho_j)\Delta_j$ of a
multi-component orbifold divisor, we add the components $\Delta_j$ one
by one, and obtain inductively the following quantitative estimate,
which is a rephrasing of Theorem~0.9.

\claim 6.16. Corollary|With a choice of 
$\gamma>\gamma_{V,D}:=\max(\max(d_j/\rho_j),\gamma_V)\ge 0$
determined by the curvature assumptions of Proposition~$6.1$, and of
hermitian metrics on $A$, $V$, $\cO_X(D)$ as prescribed by conditions
$(6.3_j)$, the orbifold metric
\vskip5pt\noindent
{\rm(a)}~~$\displaystyle
|u|_{h_{V\langle D\rangle,\varepsilon}}^2:=
|u|^2_{h_V}+\sum_{1\le j\le N}\varepsilon_j\,|\sigma_j|^{-2+2/\rho_j}\,
|\nabla_j\sigma_j(u)|^2_{h_j}$
\vskip5pt\noindent
yields a curvature tensor 
$\theta_{V\langle D\rangle,\gamma,\varepsilon}:=
\gamma\,\Theta_{A,h_{A,\delta}}\otimes\Id-
\Theta_{V\langle D\rangle,h_{V\langle D\rangle,\varepsilon}}$
such that the associated quadratic form $Q_{V\langle D\rangle,\gamma,\varepsilon}$ on
$T_X\otimes V$ satisfies for
$\varepsilon_N\ll\varepsilon_{N-1}\ll\cdots\ll\varepsilon_1\ll 1$
the curvature estimate
$$
\leqalignno{\kern10pt Q_{V\langle D\rangle,\gamma,\varepsilon}
&(z)(\xi\otimes u)\simeq\gamma\,\ddbar\psi_\delta(\xi,\xi)\,|u|^2
+\sum_{\ell,m,\lambda,\mu}c_{\ell m\lambda\mu}\,\xi_\ell\overline\xi_m\,
u_\lambda\overline u_\mu&{\rm(b)}\cr
+\sum_j&\varepsilon_j\,|\sigma_j|^{-2+2/\rho_j}\,
\big(\gamma\,\ddbar\psi_\delta(\xi,\xi)-\rho_j^{-1}\,\ddbar\varphi_j(\xi,\xi)\big)
\,|\nabla_j\sigma_j(u)|^2\cr
+\sum_j&{\varepsilon_j\,|\sigma_j|^{-2+2/\rho_j}\over
1+\varepsilon_j\,|\sigma_j|^{-2+2/\rho_j}\,
|\nabla_j\sigma_j|^2}\,\big|\,\nabla_j^2\sigma_j(\xi,u)-(1-1/\rho_j)\,
\sigma_j^{-1}\,\nabla_j\sigma_j(\xi)\nabla_j\sigma_j(u)\big|^2,\cr}
$$
where
$$
\nabla_{A,h_{A,\delta}}^2=\ddbar\psi_\delta,\quad
\nabla_{\Delta_j,h_j}^2=\ddbar\varphi_j,\quad
\hbox{$(c_{\ell m\lambda\mu})={}$coefficients of $-2\pi\,\Theta_{V,h_V}$}.
$$
Here, the symbol $\simeq$ means that the ratio of the left and right hand
sides can be chosen in $[1-\alpha,1+\alpha]$ for any
$\alpha>0$ prescribed in advance.
\endclaim

\plainsubsection 6.B. Evaluation of some Chern form integrals and their limits|

Our aim is to apply Lemma 5.8 and Corollary 6.16 to compute Morse
integrals of the curvature tensor of a directed orbifold~$(X,V,D)$,
where $D=\sum_j(1-1/\rho_j)\Delta_j$ is transverse to~$V$. Let $A\in
\Pic(X)$ be an ample line bundle, and $d_j$, $\gamma_V$,
$\gamma>\gamma_{V,D}$ be defined as in~6.16.
We get hermitian metrics
$h_{V\langle D\rangle,\varepsilon}$ on $V\langle D\rangle$ and
corresponding curvature tensors
$\theta_{V\langle D\rangle,\gamma,\varepsilon}$ in $C^\infty(X\ssm|D|,
\Lambda^{1,1}T^*_X\otimes\Hom(V,V))$ that are ``orbifold smooth'', and
such that $\theta_{V\langle D\rangle,\gamma,\varepsilon}\ge_G 0$.
Given a smooth strongly positive $(n-p,n-p)$-form $\beta\ge_S 0$
on~$X$, we want to evaluate the integrals
$$
\leqalignno{
I_{p,\varepsilon}(\beta)&=\int_{S_\varepsilon(V\langle D\rangle)}
\langle\theta_{V\langle D\rangle,\gamma,\varepsilon}\cdot u,u\rangle^p
\wedge\beta\,d\mu_\varepsilon(u)&(6.17)\cr
&=\int_{z\in X}\int_{u\in S_\varepsilon(V\langle D\rangle)_z}
\langle\theta_{V\langle D\rangle,\gamma,\varepsilon}\cdot u,u\rangle^p
\wedge\beta(z)\,d\mu_\varepsilon(u),
&(6.17')\cr}
$$
where $S_\varepsilon(V\langle D\rangle)$ denotes the unit sphere bundle
of $V\langle D\rangle$ with respect to $h_\varepsilon$, and
$\mu_\varepsilon$ the unitary invariant probability measure on the sphere.
Proposition 5.13~(b) and the Fubini theorem imply the upper bound
$$
I_{p,\varepsilon}(\beta)\le {p!\,(r-1)!\over(p+r-1)!}
\int_X (\Tr\theta_{V\langle D\rangle,\gamma,\varepsilon})^p\wedge \beta.
\leqno(6.18)
$$
When $\beta$ is closed, the upper bound can be evaluated by a
cohomology class calculation, thanks to the following lemma.

\claim 6.19. Lemma|The $(1,1)$-form
$\Tr\theta_{V\langle D\rangle,\gamma,\varepsilon}\ge 0$ is closed and
belongs to the cohomology class\vskip8pt
\centerline{$\displaystyle
r\gamma\,c_1(A)-c_1(V)+\sum_j(1-1/\rho_j)c_1(\Delta_j).
$}
\endclaim

\plainproof. The trace can be seen as the curvature of
$$
\det(\cO_X(\gamma A)\otimes V\langle D\rangle^*)
=\cO_X(r\gamma\,A)\otimes \det(V\langle D\rangle^*)
=\cO_X(r\gamma\,A)\otimes \det(V^*)\otimes\cO_X(D)
$$
with the determinant metric. Since all metrics have equivalent behaviour
along $|D|$ (and can be seen as orbifold smooth), Stokes' theorem
shows that the cohomology class is independent of $\varepsilon$. Formally,
the result follows from (1.25). One can also consider the intersection product
$$
\{\Tr\theta_{V\langle D\rangle,\gamma,\varepsilon}\}\cdot \{\beta\}
=\int_X\Tr\theta_{V\langle D\rangle,\gamma,\varepsilon}\wedge\beta
=r\int_{u\in S_\varepsilon(V\langle D\rangle)}
\langle\theta_{V\langle D\rangle,\gamma,\varepsilon}\cdot u,u\rangle
\wedge\beta\,d\mu_\varepsilon(u)
$$
for all smooth closed
$(n-1,n-1)$-forms $\beta$ on $X$, and apply Corollary 6.16~(b)
to evaluate the limit as $\varepsilon\to 0$. This will be checked later as
the special case $p=1$ of (6.17).\qed

\noindent
We actually need even more
general estimates. The proof follows again from the Fubini theorem.

\claim 6.20. Proposition|Consider orbifold directed structures
$(X,V,D_s)$, $1\le s\le k$, with
$D_s=\sum_{1\le j\le N}(1-{1\over\rho_{s,j}})\Delta_j$.
We assume that the divisors $D_s$ are simple normal crossing divisors 
transverse to~$V$, sharing the same components~$\Delta_j$. Let $d_j$
be the infimum of numbers
$\lambda\in\bR_+$ such that $\lambda\,A-\Delta_j$ is~nef, and let
$\gamma_V$ be the infimum of numbers
$\gamma\ge 0$ such that 
$\theta_{V,\gamma}:=\gamma\,\Theta_{A,h_A}\otimes\Id_V-\Theta_{V,h_V}\ge_G0$
for suitable hermitian metrics $h_V$ on~$V$.
Take $p=(p_1,\ldots,p_k)\in\bN^k$ such that $p'=n-(p_1+\ldots+p_k)\ge 0$ and a
smooth, closed, strongly positive $(p',p')$ form $\beta\ge_S0$ on $X$.
Then for every
$$
\gamma_s>\gamma_{V,D_s}:=\max(\max_j(d_j/\rho_{s,j}),\gamma_V)
$$
there exist hermitian metrics
$h_{V\langle D_s\rangle,\varepsilon_s}$ on the orbifold vector bundles
$V\langle D_s\rangle$ such that
$$
\theta_{V\langle D_s\rangle,\gamma_s,\varepsilon_s}:=\gamma_s
\,\Theta_{A,h_A}\otimes\Id_V-\Theta_{V\langle D_s\rangle,
h_{V\langle D_s\rangle,\varepsilon_s}}>_G0,\quad
\varepsilon_s=(\varepsilon_{s1},\ldots,\varepsilon_{sN}),\quad 1\le s\le k
$$
in the sense of Griffiths. Moreover, the integrals
$$
I_{k,p,\varepsilon}(\beta)=
\int_{z\in X}\int_{(u_s)\,\in\,\Pi_s S(V\langle D_s\rangle)_z}~\bigwedge_{1\le s\le k}
\langle\theta_{V\langle D_s\rangle,\gamma_s,\varepsilon_s}(u_s),u_s\rangle^{p_s}
\wedge\beta(z)\,\prod_{1\le s\le k}d\mu_{\varepsilon_s}(u_s)
\leqno{\rm(a)}
$$
admit upper bounds
$$
I_{k,p,\varepsilon}(\beta)\le
\int_X\bigwedge_{1\le s\le k}{p_s!\,(r-1)!\over(p_s+r-1)!}
\Big(r\gamma_s\,\Theta_{A,h_{A,\delta}}-\Tr\Theta_{V,h_V}+
\sum_j(1-1/\rho_{s,j})\Theta_{\Delta_j,h_j}\Big)^{p_s}\!\wedge\beta.
\leqno{\rm(b)}
$$
When $\beta$ is closed, we get a purely cohomological upper bound
$$
I_{k,p,\varepsilon}(\beta)\le
\int_X\prod_{1\le s\le k}{p_s!\,(r-1)!\over(p_s+r-1)!}\,
\Big(r\gamma_s\,c_1(A)-c_1(V)+\sum_j(1-1/\rho_{s,j})c_1(\Delta_j)\Big)^{p_s}
\cdot\{\beta\}.
\leqno{\rm(c)}
$$
\endclaim

\claim 6.21. Complement|When $p_1=\ldots=p_k=1$, formulas {\rm 6.20~(b)} and
{\rm 6.20~(c)} are equalities.
\endclaim

\plainproof. This follows from Remark 5.15.\qed

\noindent
In general, getting a lower bound for $I_{p,\varepsilon}(\beta)$ and
$I_{k,p,\varepsilon}(\beta)$ is substantially harder. We start with
$I_{p,\varepsilon}(\beta)$ and content ourselves to evaluate the iterated limit
$$
\lim_{\varepsilon\to 0}I_{p,\varepsilon}(\beta):=
\lim_{\varepsilon_1\to 0}\lim_{\varepsilon_2\to 0}\ldots\lim_{\varepsilon_N\to 0}
I_{p,\varepsilon}(\beta),\quad
~~\varepsilon_N\ll\varepsilon_{N-1}\ll\cdots\ll\varepsilon_1\ll 1.
\leqno(6.22)
$$
For this, we consider the expression of the curvature form
in a neigborhood of an arbitrary point
$z_0\in\Delta_{j_1}\cap\ldots\cap\Delta_{j_m}$ (if $z_0\in X\ssm|\Delta|$,
we have $m=0$). We take trivializations of the line bundles
$\cO_X(\Delta_j)$ so that the hermitian metrics have weights
$e^{-\varphi_j}$ with $\varphi_j(z_0)=d\varphi_j(z_0)=0$, and
introduce the corresponding ``orbifold'' coordinates
$$
t_{j,\varepsilon}=\varepsilon_j^{1/2}\,\sigma_j(z)^{-(1-1/\rho_j)}\,|\nabla_j\sigma_j(z_0)|,
\quad j=j_1,\ldots,j_m,
\leqno(6.23)
$$
We complete these coordinates with $n-m$ variables $z_\ell$ that
define coordinates along $\Delta_{j_1}\cap\ldots\cap\Delta_{j_m}$. In
this way, we get a $n$-tuple $(t_{j,\varepsilon},z_\ell)$ of complex numbers that
provide local coordinates on the universal cover of $\Omega_{z_0}\ssm|D|$,
where $\Omega_{z_0}$ is a small neighborhood of~$z_0$. Viewed on $X$, the
coordinates $t_{j,\varepsilon}$ are multivalued near $z_0$, but we can
make a ``cut'' in~$X$ along $\Delta_j$ to exclude the negligible set 
of points where
$\sigma_j(z)\in\bR_-$, and take the argument in~$]-\pi,\pi[$, so that
$\Arg(t_{j,\varepsilon})\in{}]-(1-1/\rho_j)\pi,(1-1/\rho_j)\pi[\,$. If we integrate
over complex numbers $t_{j,\varepsilon}$ without such a restriction on the argument, 
the integral will have to be multiplied by the factor $(1-1/\rho_j)$ to get
the correct value. Since $|\sigma_j|$ is bounded, the range of the
absolute value $|t_{j,\varepsilon}|$ is an interval
$]O(\varepsilon_j^{1/2}),+\infty[\,$, thus $t_{j,\varepsilon}$ will cover asymptotically
an entire angular sector in $\bC$ as $\varepsilon_j\to 0$.
In the above coordinates, we~have
$$
{dt_{j,\varepsilon}\over t_{j,\varepsilon}}=-(1-1/\rho_j){d\sigma_j\over\sigma_j}=
-(1-1/\rho_j)\bigg({\nabla_j\sigma_j\over\sigma_j}+\partial\varphi_j\bigg)=
-(1-1/\rho_j){\nabla_j\sigma_j\over\sigma_j}+O(1),
\leqno(6.24)
$$
since $\nabla_j\sigma_j=d\sigma_j-\sigma_j\partial\varphi_j$ and the weight
$\varphi_j$ of the metric of $\cO_X(\Delta_j)$ is smooth.
Denote
$$
\leqalignno{
&\theta_{V,\gamma}=\gamma\,\Theta_{A,h_{A,\delta}}\otimes\Id_V-\Theta_{V,h_V},
&(6.25_1)\cr
\noalign{\vskip5pt}
&\theta_{V\langle D\rangle,\gamma,\varepsilon}:=\gamma\,\Theta_{A,h_{A,\delta}}\otimes\Id-
\Theta_{V\langle D\rangle,h_{V\langle D\rangle,\varepsilon}},&(6.25_2)\cr
&e_j^*={\nabla_j\sigma_j\over|\nabla_j\sigma_j|}\in S(V^*).&(6.25_3)\cr}
$$
By Corollary~6.16, we have
$$
\leqalignno{&\kern24pt
\langle\theta_{V\langle D\rangle,\gamma,\varepsilon}\cdot u,u\rangle
\simeq\gamma\,\Theta_{A,h_{A,\delta}}-\langle\Theta_{V,h_V}\cdot u,u\rangle
&(6.26)\cr
\noalign{\vskip8pt}  
&+\sum_j\varepsilon_j\,|\sigma_j|^{-2+2/\rho_j}\,
\big(\gamma\,\Theta_{A,h_{A,\delta}}-\rho_j^{-1}\,\Theta_{\Delta_j,h_j}\big)
\,|\nabla_j\sigma_j(u)|^2\cr
&+{1\over 2\pi}\sum_j{\varepsilon_j\,|\sigma_j|^{-2+2/\rho_j}\over
1+\varepsilon_j\,|\sigma_j|^{-2+2/\rho_j}\,
|\nabla_j\sigma_j|^2}\,\big|\,\nabla_j^2\sigma_j(\xi,u)-(1-1/\rho_j)\,
\sigma_j^{-1}\,\nabla_j\sigma_j(\xi)\nabla_j\sigma_j(u)\big|^2,\cr}
$$
therefore
$$
\leqalignno{\kern32pt
\langle\theta_{V\langle D\rangle,\gamma,\varepsilon}\cdot u,u\rangle  
&\simeq\langle\theta_{V,\gamma}\cdot u,u\rangle+\sum_j
\big(\gamma\,\Theta_{A,h_{A,\delta}}-\rho_j^{-1}\,\Theta_{\Delta_j,h_j}\big)
\,|t_{j,\varepsilon}|^2\,|e_j^*(u)|^2&(6.27_1)\cr
&\kern10pt{}+{i\over 2\pi}
\sum_j{|t_{j,\varepsilon}|^2\over 1+|t_{j,\varepsilon}|^2}\,\Big\langle
{dt_{j,\varepsilon}\over t_{j,\varepsilon}}\,e_j^*(u)+b_j(u)\,,\,
{dt_{j,\varepsilon}\over t_{j,\varepsilon}}\,e_j^*(u)+b_j(u)\Big\rangle_{h_j},&(6.27_2)\cr}
$$
where 
$$
b_j={1\over|\nabla_j|}\nabla^2\sigma_j\in 
C^\infty(\Omega_{z_0},\Lambda^{1,0}T^*_X\otimes V^*\otimes\cO_X(\Delta_j))
\leqno(6.27_3)
$$
is a smooth $(1,0)$-form near $z_0$.
The approximate equality $\simeq$ in formula $(6.27_{1,2})$ involves
the approximation $|\nabla_j\sigma_j(z)|/|\nabla_j\sigma_j(z_0)|\simeq 1$,
which holds in a sufficiently small neighborhood of~$z_0\,$;
if we apply the Fubini theorem and consider the fiber integral over
$z_0\in X$, there is actually no error coming from this approximation.
Now, we want to integrate the volume form 
$\langle\theta_{V\langle D\rangle,\gamma,\varepsilon}\cdot u,u\rangle^p
\wedge\beta\,d\mu_\varepsilon(u)$ along the fibers of
$S_\varepsilon(V\langle D\rangle)\to X$.
The sphere bundle $S_\varepsilon(V\langle D\rangle)$ is defined by
$|u|_{h_{V\langle D\rangle,\varepsilon}}^2=1$
where
$$
|u|_{h_{V\langle D\rangle,\varepsilon}}^2=
|u|^2+\sum_j\varepsilon_j|\sigma_j|^{-2+2/\rho_j}\,|\nabla_j\sigma_j(u)|^2
\simeq
|u|^2+\sum_j|t_{j,\varepsilon}|^2|e_j^*(u)|^2=1.\leqno(6.28)
$$
For the sake of simplicity, we first deal with the case where the
divisor $D=(1-1/\rho_j)\Delta_j$ has a single component. Along 
$\Delta_j$, we then get an orthogonal decomposition
$V=(V\cap T_{\Delta_j})\oplus\bC e_j$, and by (6.28) we can write
$$
u=u'_j+e_j^*(u)\,e_j\in S(V),\quad
|u|^2=|u'_j|^2+|e_j^*(u)|^2,\quad
u'_j\in V\cap T_{\Delta_j}.
\leqno(6.28^\perp)
$$
We reparametrize the  integration in $u\in S_\varepsilon(V\langle D\rangle)$
on the sphere $S(V)$ by introducing the change of variables 
$$
\eqalign{
&\tau=\tau_{j,\varepsilon}={|t_{j,\varepsilon}|^2\over 1+|t_{j,\varepsilon}|^2}
\in[0,1],\quad
1-\tau={1\over 1+|t_{j,\varepsilon}|^2},\quad
d\tau={d|t_{j,\varepsilon}|^2\over
(1+|t_{j,\varepsilon}|^2)^2},\cr
\noalign{\vskip5pt}  
&g_{j,\varepsilon}:S(V)\to S_\varepsilon(V\langle D\rangle),\quad
u\mapsto u_{j,\varepsilon}=u'_j+\sqrt{1-\tau}\,e_j^*(u)\,e_j
=u'_j+{e_j^*(u)\over(1+|t_{j,\varepsilon}|^2)^{1/2}}\,e_j,\cr}
$$
so that $u_{j,\varepsilon}$ satisfies
$|u_{j,\varepsilon}|^2_\varepsilon=|u|^2$ and
$$
e_j^*(u_{j,\varepsilon})=\sqrt{1-\tau}\,e_j^*(u)
={1\over (1+|t_{j,\varepsilon}|^2)^{1/2}}\,e_j^*(u),
\quad
|t_{j,\varepsilon}|^2\,|e_j^*(u_{j,\varepsilon})|^2=\tau\,|e_j^*(u)|^2.
$$
This gives $d\mu_\varepsilon(u_{j,\varepsilon})=d\mu(u)$, and
as a consequence (6.17) can be rewritten as
$$
I_{p,\varepsilon}(\beta)=\int_{S(V)}\langle\theta_{V\langle D\rangle,\gamma,\varepsilon}\cdot 
u_{j,\varepsilon},u_{j,\varepsilon}\rangle^p\wedge\beta(z)\,d\mu(u).
\leqno(6.29)
$$
Finally, a use of polar coordinates with 
$\alpha=\Arg(t_{j,\varepsilon})$ shows that
$$
{\ii dt_{j,\varepsilon}\wedge d\overline t_{j,\varepsilon}\over
(1+|t_{j,\varepsilon}|^2)^2}=
{2\,|t_{j,\varepsilon}|\,d|t_{j,\varepsilon}|\wedge d\alpha\over
(1+|t_{j,\varepsilon}|^2)^2}
=d\tau\wedge d\alpha.
$$
A substitution $u\mapsto u_{j,\varepsilon}$ in $(6.27_{1,2})$ yields
$$
\leqalignno{
&\langle\theta_{V\langle D\rangle,\gamma,\varepsilon}\cdot u_{j,\varepsilon},
u_{j,\varepsilon}\rangle  
\simeq\langle\theta_{V,\gamma}\cdot u_{j,\varepsilon},u_{j,\varepsilon}
\rangle+
\big(\gamma\,\Theta_{A,h_{A,\delta}}-\rho_j^{-1}\,\Theta_{\Delta_j,h_j}\big)
\,{|t_{j,\varepsilon}|^2\,|e_j^*(u)|^2\over 1+|t_{j,\varepsilon}|^2}\cr
&\kern25pt{}+{\ii\over 2\pi}{1\over 1+|t_{j,\varepsilon}|^2}\,\Big\langle
{dt_{j,\varepsilon}\,e_j^*(u)\over(1+|t_{j,\varepsilon}|^2)^{1/2}}+
t_{j,\varepsilon}\,b_j(u_{j,\varepsilon})\,,\,
{dt_{j,\varepsilon}\,e_j^*(u)\over(1+|t_{j,\varepsilon}|^2)^{1/2}}+
t_{j,\varepsilon}\,b_j(u_{j,\varepsilon})\Big\rangle_{h_j}.&(6.30)\cr}
$$
The last term is a $(1,1)$-form that is a square of a $(1,0)$-form 
(when $u$ is fixed), hence the expansion of the $p$-th power can involve
at most one such factor. Therefore we get
$$
\leqalignno{
&\langle\theta_{V\langle D\rangle,\gamma,\varepsilon}\cdot u_{j,\varepsilon},
u_{j,\varepsilon}\rangle^p
\simeq
\bigg(
\langle\theta_{V,\gamma}\cdot u_{j,\varepsilon},u_{j,\varepsilon}\rangle
+\big(\gamma\,\Theta_{A,h_{A,\delta}}-\rho_j^{-1}\,\Theta_{\Delta_j,h_j}\big)
\,{|t_{j,\varepsilon}|^2\,|e_j^*(u)|^2\over 1+|t_{j,\varepsilon}|^2}\bigg)^p\cr
&\kern20pt{}+p\,{\ii\over 2\pi}{1\over 1+|t_{j,\varepsilon}|^2}\,\Big\langle
{dt_{j,\varepsilon}\,e_j^*(u)\over(1+|t_{j,\varepsilon}|^2)^{1/2}}+
t_{j,\varepsilon}\,b_j(u_{j,\varepsilon})\,,\,
{dt_{j,\varepsilon}\,e_j^*(u)\over(1+|t_{j,\varepsilon}|^2)^{1/2}}+
t_{j,\varepsilon}\,b_j(u_{j,\varepsilon})\Big\rangle\wedge{}\cr
&\kern50pt\bigg(
\langle\theta_{V,\gamma}\cdot u_{j,\varepsilon},u_{j,\varepsilon}\rangle
+\big(\gamma\,\Theta_{A,h_{A,\delta}}-\rho_j^{-1}\,\Theta_{\Delta_j,h_j}\big)
\,{|t_{j,\varepsilon}|^2\,|e_j^*(u)|^2\over 1+|t_{j,\varepsilon}|^2}
\bigg)^{p-1}.&(6.30^p)\cr}
$$
The integrals involving $b_j(u_{j,\varepsilon})$ are of the form
$$
\int_{S(V)}{\overline t_{j,\varepsilon}\,dt_{j,\varepsilon}\wedge
\langle e_j^*(u),b_j(u_{j,\varepsilon})\rangle
\over(1+|t_{j,\varepsilon}|^2)^{3/2}}\wedge A_{j,\varepsilon}(u),\quad
\int_{S(V)}{|t_{j,\varepsilon}|^2
\langle b_j(u_{j,\varepsilon}),b_j(u_{j,\varepsilon})\rangle
\over 1+|t_{j,\varepsilon}|^2}\wedge A'_{j,\varepsilon}(u)
$$
where $A_{j,\varepsilon}(u)$, $A'_{j,\varepsilon}(u)$ are forms with
uniformly bounded coefficients in orbifold coordinates.
Since ${|t_{j,\varepsilon}|^2\over 1+|t_{j,\varepsilon}|^2}$ is
bounded by $1$ and converges to $0$ on $X\ssm\Delta_j$, Lebesgue's
dominated convergence theorem shows that the second integral converges to~$0$.
The second integral can be estimated by the Cauchy-Schwarz inequality.
We obtain an upper bound
$$
\bigg(\int_{S(V)}{|t_{j,\varepsilon}|^2
\langle b_j(u_{j,\varepsilon}),b_j(u_{j,\varepsilon})\rangle
\over 1+|t_{j,\varepsilon}|^2}\wedge A_{j,\varepsilon}(u)\bigg)^{1/2}
\bigg(\int_{S(V)}{\ii dt_{j,\varepsilon}\wedge d\overline t_{j,\varepsilon}\,
|e_j^*(u)|^2
\over(1+|t_{j,\varepsilon}|^2)^2}\wedge A_{j,\varepsilon}(u)\bigg)^{1/2}
$$
where the first factor converges to $0$ and the second one is bounded
by Fubini, since $\int_\bC \ii dt\wedge d\overline t/(1+|t|^2)^2<+\infty$.
Modulo negligible terms, and changing variables into our new parameters
$(\tau,\alpha)$,  we finally obtain
$$
\leqalignno{
&\langle\theta_{V\langle D\rangle,\gamma,\varepsilon}\cdot u_{j,\varepsilon},
u_{j,\varepsilon}\rangle^p
\simeq
\langle\theta_{V,\gamma}\cdot u_{j,\varepsilon},u_{j,\varepsilon}\rangle^p\cr
&\kern30pt{}+p\,{d\tau\wedge d\alpha\over 2\pi}\,|e_j^*(u)|^2\wedge\Big(
\langle\theta_{V,\gamma}\cdot u_{j,\varepsilon},u_{j,\varepsilon}\rangle
+\big(\gamma\,\Theta_{A,h_{A,\delta}}-\rho_j^{-1}\,\Theta_{\Delta_j,h_j}\big)
\,\tau\,|e_j^*(u)|^2\Big)^{p-1}.&(6.31)\cr}
$$
Therefore
$$
\leqalignno{
&\int_{S(V)}\langle\theta_{V\langle D\rangle,\gamma,\varepsilon}\cdot u_{j,\varepsilon},
u_{j,\varepsilon}\rangle^p\wedge\beta\,d\mu(u)
\simeq \int_{S(V)}
\langle\theta_{V,\gamma}\cdot u_{j,\varepsilon},u_{j,\varepsilon}\rangle^p
\wedge\beta\,d\mu(u)
\cr
&\kern6pt{}+\int_{S(V)}p\,{d\tau\wedge d\alpha\over 2\pi}\,|e_j^*(u)|^2\wedge\Big(
\langle\theta_{V,\gamma}\cdot u_{j,\varepsilon},u_{j,\varepsilon}\rangle
+\big(\gamma\,\Theta_{A,h_{A,\delta}}-\rho_j^{-1}\,\Theta_{\Delta_j,h_j}\big)
\,\tau\,|e_j^*(u)|^2\Big)^{p-1}\cr
\noalign{\vskip4pt}
&\kern120mm{}\wedge\beta\,d\mu(u).&(6.32)\cr}
$$
Since $u_{j,\varepsilon}\to u$ almost everywhere and boundedly, we have
$$
\lim_{\varepsilon\to 0}
\int_{S(V)}\langle\theta_{V,\gamma}\cdot u_{j,\varepsilon},u_{j,\varepsilon}\rangle^p
\wedge\beta\,d\mu(u)=
\int_{S(V)}
\langle\theta_{V,\gamma}\cdot u,u\rangle^p\wedge\beta\,d\mu(u).
$$
Here, we have to remember that $\tau=\tau_{j,\varepsilon}$ converges uniformly
to $0$ (even in the $C^\infty$ topology), on all compact subsets 
of $X\ssm\Delta_j$, hence the second integral in (6.32) 
asymptotically concentrates on $\Delta_j$ as $\varepsilon\to 0$.
Also, the angle $\alpha=\Arg(t_{j,\varepsilon})$ runs 
over the interval $]-(1-1/\rho_j)\pi,(1-1/\rho_j)\pi[$. In the easy case 
$p=1$, we get
$$
\eqalign{
\lim_{\varepsilon\to 0}
\int_{S(V)}\langle\theta_{V,\gamma,\varepsilon}&\cdot u_{j,\varepsilon},
u_{j,\varepsilon}\rangle\wedge\beta\,d\mu(u)\cr
&=\int_{S(V)}\langle\theta_{V,\gamma}\cdot u,u\rangle\wedge\beta\,d\mu(u)
+(1-1/\rho_j)\int_{S(V)_{|\Delta_j}}|e_j^*(u)|^2\beta\,d\mu(u)\cr
&=\int_{X}{1\over r}\,\Tr\theta_{V,\gamma}\wedge\beta
+(1-1/\rho_j)\int_{\Delta_j}{1\over r}\beta.\cr}
$$
If we assume $\beta$ closed, this is equal to the intersection product
$$
{1\over r}\big(\rho\gamma\,c_1(A)-c_1(V)+(1-1/\rho_j)
c_1(\Delta_j)\big)\cdot\beta
$$
and the final assertion of the proof of Lemma~6.19 is thus confirmed, 
adding the components $\Delta_j$ one by one (see below).
Now, in the general case $p\ge 1$, we will obtain a lower bound of
the second integral involving $d\tau\wedge d\alpha$ in (6.32) by using
a change of variable
$$
\eqalign{
&h_{j,\varepsilon}:S(V)\to S(V),\quad u\mapsto h_{j,\varepsilon}(u)
=\big((1-\tau)|u'|^2+|e_j^*(u)|^2\big)^{-1/2}\,
\big(\sqrt{1-\tau}\,u'+e_j^*(u)\,e_j\big)\cr}
$$
where $\tau=\tau_{j,\varepsilon}$.
Observe that the composition $g_{j,\varepsilon}\circ h_{j,\varepsilon}:
S(V)\to S(V)\to S_\varepsilon(V\langle D\rangle)$ is given by
$$
g_{j,\varepsilon}\circ h_{j,\varepsilon}(u)={\sqrt{1-\tau}\over
\big((1-\tau)|u'|^2+|e_j^*(u)|^2\big)^{1/2}}\,u.
$$
Since $(1-\tau)|u'|^2+|e_j^*(u)|^2\le|u|^2=1$,
it is easy to check that $d\mu(h_{j,\varepsilon}(u))\ge(1-\tau)^{r-1}\,d\mu(u)$
on the unit sphere, that $|e_j^*(h_{j,\varepsilon}(u))|\ge |e_j^*(u)|$, and
finally, that
$$
\langle\theta_{V,\gamma}\cdot g_{j,\varepsilon}(h_{j,\varepsilon}(u)),g_{j,\varepsilon}(h_{j,\varepsilon}(u))
\rangle=
{1-\tau\over (1-\tau)|u'|^2+|e_j^*(u)|^2}\,
\langle\theta_{V,\gamma}\cdot u,u\rangle
\ge(1-\tau)\,\langle\theta_{V,\gamma}\cdot u,u\rangle.
$$
Hence, by a change a variable $u\mapsto h_{j,\varepsilon}(u)$ we find
$$
\leqalignno{
&\int_{S(V)}p\,{d\tau\wedge d\alpha\over 2\pi}\,|e_j^*(u)|^2\wedge\Big(
\langle\theta_{V,\gamma}\cdot u_{j,\varepsilon},u_{j,\varepsilon}\rangle
+\big(\gamma\,\Theta_{A,h_{A,\delta}}-\rho_j^{-1}\,\Theta_{\Delta_j,h_j}\big)
\,\tau\,|e_j^*(u)|^2\Big)^{p-1}\cr
&\kern120mm{}\wedge\beta\,d\mu(u)\cr
&\ge\int_{S(V)}p\,{d\tau\wedge d\alpha\over 2\pi}\,|e_j^*(u)|^2\wedge\Big(
(1-\tau)\langle\theta_{V,\gamma}\cdot u,u\rangle
+\big(\gamma\,\Theta_{A,h_{A,\delta}}-\rho_j^{-1}\,\Theta_{\Delta_j,h_j}\big)
\,\tau\,|e_j^*(u)|^2\Big)^{p-1}\cr
\noalign{\vskip4pt}
&\kern100mm{}\wedge\beta\,(1-\tau)^{r-1}\,d\mu(u).&(6.33)\cr}
$$
Here, we have to remember that $\tau=\tau_{j,\varepsilon}$ converges uniformly
to $0$ (even in the $C^\infty$ topology), on all compact subsets 
of $X\ssm\Delta_j$.  Therefore, the 
last integral concentrates over the divisor~$\Delta_j$. If we apply
the binomial formula with an index $q'=q-1$, we see that the
limit as $\varepsilon\to 0$ is equal to
$$
\leqalignno{
&p(1-1/\rho_j)\int_{S(V)_{|\Delta_j}}\sum_{q=1}^{p}{p-1\choose q-1}
\langle\theta_{V,\gamma}\cdot u,u\rangle^{p-q}\wedge
\big(\gamma\,\Theta_{A,h_{A,\delta}}-\rho_j^{-1}\,\Theta_{\Delta_j,h_j}\big)^{q-1}\,
|e_j^*(u)|^{2q}\cr
&\kern50mm\bigg(\int_0^1(1-\tau)^{p-q+r-1}\,\tau^{q-1}\,d\tau\bigg)\wedge
\beta\,d\mu(u).&(6.34)\cr}
$$
We have
$$
\int_0^1(1-\tau)^{p-q+r-1}\,\tau^{q-1}\,d\tau={
(p-q+r-1)!\,(q-1)!\over (p+r-1)!}
\leqno(6.35)
$$ 
and the combination of (6.29) and $(6.32 - 6.35)$ implies
$$
\leqalignno{
&\lim_{\varepsilon\to 0}I_{p,\varepsilon}(\beta)\ge
\int_{S(V)}\langle\theta_{V,\gamma}\cdot u,u\rangle^p\wedge\beta\,d\mu(u)
+p(1-1/\rho_j)\,\sum_{q=1}^{p}{(p-1)!\,(p-q+r-1)!\over
(p-q)!\,(p+r-1)!}\times{}\cr
&\kern60pt{}\int_{S(V)_{|\Delta_j}}
\langle\theta_{V,\gamma}\cdot u,u\rangle^{p-1-q}\wedge
\big(\gamma\,\Theta_{A,h_{A,\delta}}-\rho_j^{-1}\,\Theta_{\Delta_j,h_j}\big)^q\,
|e_j^*(u)|^{2q+2}\wedge\beta\,d\mu(u).&(6.36)\cr}
$$
Inductively, formula (6.36) requires the investigation of more general 
integrals
$$
I_{p,p',Y,\varepsilon}=
\int_{S_\varepsilon(V\langle D\rangle)_{|Y}}
\langle\theta_{V,\gamma,\varepsilon}\cdot u,u\rangle^{p-p'}\wedge 
\prod_{1\le j\le p'}|\ell_j(u)|^2\,\beta\,d\mu_\varepsilon(u)\leqno(6.37)
$$
where $Y$ is a subvariety of $X$ (which we assume to be transverse
to the $\Delta_j$'s, and $\ell_j\in C^\infty(Y,V^*)$  with $|\ell_j|=1$,
and $\beta\ge_S 0$ is a smooth form of suitable bidegree on $Y$.
Not much is changed in the calculation, except that the change of variable
$u\mapsto g_{j,\varepsilon}\circ h_{j,\varepsilon}(u)$ applied to 
$\prod_{1\le j\le p'}|\ell_j(u)|^2$ introduces an extra factor $(1-\tau)^{p'}$
in the lower bound, entirely compensated by the corresponding factor 
$(1-\tau)^{p-p'-q}$ appearing in 
$\langle\theta_{V,\gamma,\varepsilon}\cdot u,u\rangle^{p-p'}$.
The binomial formula yields a coefficient ${p-p'-1\choose q-1}$ instead of
${p-1\choose q-1}$. We thus obtain 
$$
\leqalignno{
&\lim_{\varepsilon\to 0}I_{p,p',Y,\varepsilon}(\beta)\ge
\int_{S(V)_{|Y}}\langle\theta_{V,\gamma}\cdot u,u\rangle^{p-p'}\wedge
\prod_{1\le j\le p'}|\ell_j(u)|^2\,\beta\,d\mu(u)\cr
&\kern100pt{}+(1-1/\rho_j)\,\sum_{q=1}^{p-p'}{(p-p')!\,(p-q+r-1)!\over
(p-p'-q)!\,(p+r-1)!}\times{}&(6.38)\cr
&\int_{S(V)_{|Y\cap\Delta_j}}\kern-4pt
\langle\theta_{V,\gamma}\cdot u,u\rangle^{p-p'-q}\wedge
\big(\gamma\,\Theta_{A,h_{A,\delta}}-\rho_j^{-1}\,
\Theta_{\Delta_j,h_j}\big)^{q-1}\!
\wedge |e_j^*(u)|^{2q}\prod_{1\le j\le p'}|\ell_j(u)|^2\,
\beta\,d\mu(u).\cr}
$$
When $D$ contains several components, we apply induction on $N$ and put
$$
\leqalignno{
&|u|_{h_{V\langle D\rangle,\varepsilon}}^2=
|u|_{h_{V\langle D\rangle,\varepsilon'}}^2
+\varepsilon_N\,|\sigma_N|^{-2+2/\rho_N}\,
|\nabla_N\sigma_N(u)|^2_{h_N}\quad\hbox{where}&(6.39)\cr
&|u|_{h_{V\langle D\rangle,\varepsilon'}}^2=
|u|^2_{h_V}+\sum_{1\le j\le N-1}\varepsilon_j\,|\sigma_j|^{-2+2/\rho_j}\,
|\nabla_j\sigma_j(u)|^2_{h_j}.&(6.39')\cr}
$$
In this setting, (6.26) can be rewritten in the form of a decomposition
$$
\eqalign{\kern5pt&\kern-5pt
\langle\theta_{V\langle D\rangle,\gamma,\varepsilon}\cdot u,u\rangle
\simeq \langle\theta_{V\langle D\rangle,\gamma,\varepsilon'}\cdot u,u\rangle\cr
\noalign{\vskip8pt}  
&+\varepsilon_N\,|\sigma_N|^{-2+2/\rho_N}\,
\big(\gamma\,\Theta_{A,h_{A,\delta}}-\rho_N^{-1}\,\Theta_{\Delta_N,h_N}\big)
\,|\nabla_N\sigma_N(u)|^2\cr
&+{1\over 2\pi}{\varepsilon_N\,|\sigma_N|^{-2+2/\rho_N}\over
1+\varepsilon_N\,|\sigma_N|^{-2+2/\rho_N}\,
|\nabla_N\sigma_N|^2}\,\big|\,\nabla_N^2\sigma_N(\xi,u)-(1-1/\rho_N)\,
\sigma_N^{-1}\,\nabla_N\sigma_N(\xi)\nabla_N\sigma_N(u)\big|^2.\cr}
$$
By an iteration of our integral lower bound (6.38), we have to deal 
inductively with all intersections
$\Delta_J=\Delta_{j_1}\cap\ldots\cap\Delta_{j_m}$,
$J=\{j_1,\ldots,j_m\}\subset\{1,\ldots,N\}\,$;  we neglect the
self-intersection terms, since they are anyway non negative. We obtain
$$
\leqalignno{\kern36pt
&\lim_{\varepsilon\to 0}I_{p,\varepsilon}(\beta)\ge
\sum_{J\subset\{1,\ldots,N\}}|J|!\kern-5pt
\sum_{\scriptstyle (q_j)\in(\bN^*)^J\atop\scriptstyle \Sigma_{j\in J}q_j\le p}
{p!\,(p+r-1-\Sigma_{j\in J}q_j)!\over
(p+r-1)!\,(p-\Sigma_{j\in J}q_j)!}\,\prod_{j\in J}(1-1/\rho_j)
\int_{z\in \Delta_J}&(6.40)\cr
&\kern-36pt\int_{u\in S(V)_z}
\langle\theta_{V,\gamma}(z)\cdot u,u\rangle^{p-\Sigma_{j\in J}q_j}\wedge
\bigwedge_{j\in J}|e_j^*(u)|^{2q_j}\,
\big(\gamma\,\Theta_{A,h_{A,\delta}}-\rho_j^{-1}\,\Theta_{\Delta_j,h_j}\big)^{q_j-1}
\wedge\beta(z)\,d\mu(u)\cr}
$$
where $J=\emptyset$ corresponds to the integral taken over $X$, with a
coefficient equal to $1$ in that case.  By the Fubini theorem,
we get the following lower bound of $I_{k,p,\varepsilon}(\beta)$.

\claim 6.41. Proposition|With the same notation as above, assume that
$$
\gamma_s>\gamma_{V,D_s}:=\max(\max_j(d_j/\rho_{s,j}),\gamma_V),\quad 1\le s\le k.
$$
and consider the limit $\lim_{\varepsilon\to 0} I_{k,p,\varepsilon}(\beta)$
computed as an iterated limit
$\lim_{\varepsilon_{11}\to 0}\ldots\lim_{\varepsilon_{kN}\to 0}$ with respect
to the lexicographic order $(i,j)<(i',j')$ if $i<i'$ or $i=i'$ and $j<j'$.
Then we have the following lower bound, 
where the summation is taken over all disjoint subsets $J_1,\ldots,J_k
\subset\{1,2,\ldots,N\}\,:$
$$
\eqalign{
&\lim_{\varepsilon\to 0} I_{k,p,\varepsilon}(\beta)\ge{}\kern-8pt
\sum_{\scriptstyle J_1\amalg\ldots\amalg J_k\atop\scriptstyle\subset\{1,\ldots,N\}}\sum_{\scriptstyle(q_j)\in(\bN^*)^{J_1\amalg\ldots\amalg J_k}\atop
\scriptstyle \Sigma_{j\in J_s}q_j\le p_s}
\prod_{1\le s\le k}\kern-4pt{|J_s|!\,p_s!\,(p_s-\Sigma_{j\in J_s}q_j+r-1)!\over
(p_s+r-1)!\,(p_s-\Sigma_{j\in J_s}q_j)!}\kern-2pt
\prod_{j\in J_s}\kern-2pt\Big(1-{1\over \rho_{s,j}}\Big)\kern-3pt\cr
&\kern90pt
\int_{z\in\Delta_{J_1\amalg\ldots\amalg J_k}}\int_{(u_s)\in S(V)^k_z}~
\bigwedge_{1\le s\le k}\bigg(
\langle\theta_{V,\gamma_s}\cdot u_s,u_s\rangle^{p_s-\Sigma_{j\in J_s}q_j}
\wedge{}\cr
&\kern120pt\bigwedge_{j\in J_s}|e_j^*(u_s)|^{2q_j}\,
\big(\gamma_s\,\Theta_{A,h_{A,\delta}}-\rho_{s,j}^{-1}\,\Theta_{\Delta_j,h_j}
\big)^{q_j-1}\,d\mu(u_s)\bigg)\wedge\beta(z).\cr}
$$
\endclaim

\noindent
Our assumptions imply that we can take
$\theta_ {V,\gamma_s}>_G(\gamma_s-\gamma_V-\delta)\Theta_{A,h_A}\otimes\Id_V$
for every $\delta>0$. By Lemma~5.8~(b), we obtain the simpler and purely
cohomological lower bound
$$
\leqalignno{
&\lim_{\varepsilon\to 0} I_{k,p,\varepsilon}(\beta)\ge{}\kern-8pt
\sum_{\scriptstyle J_1\amalg\ldots\amalg J_k\atop\scriptstyle\subset\{1,\ldots,N\}}\sum_{\scriptstyle(q_j)\in(\bN^*)^{J_1\amalg\ldots\amalg J_k}\atop
\scriptstyle \Sigma_{j\in J_s}q_j\le p_s}
\prod_{1\le s\le k}\kern-4pt{|J_s|!\,p_s!\,(p_s-\Sigma_{j\in J_s}q_j+r-1)!\over
(p_s+r-1)!\,(p_s-\Sigma_{j\in J_s}q_j)!}
\prod_{j\in J_s}\kern-2pt\Big(1-{1\over \rho_{s,j}}\Big)\cr
&\kern90pt
\int_{\Delta_{J_1\amalg\ldots\amalg J_k}}~
\bigwedge_{1\le s\le k}\bigg(
\big((\gamma_s-\gamma_V)\Theta_{A,h_A}\big)^{p_s-\Sigma_{j\in J_s}q_j}
\wedge{}\cr
&\kern120pt\bigwedge_{j\in J_s}{(r-1)!\over (q_j+r-1)!}\,
\big(\gamma_s\,\Theta_{A,h_{A,\delta}}-\rho_{s,j}^{-1}\,\Theta_{\Delta_j,h_j}
\big)^{q_j-1}\bigg)\wedge\beta(z).&(6.42)\cr}
$$
What is a bit surprising in all these estimates is that, in spite of the fact
that we are integrating non closed and metric dependent forms, the
limits of the integrals as $\varepsilon\to 0$ admit rather natural
lower and upper bounds that are purely cohomological, and can be
expressed solely in terms of well understood Chern classes. This will
also be true for the related Morse integrals in~\S$\,$7. It could be
desirable to have an algebro-geometric explanation of this
phenomenon. The algebraic versions of Morse inequalities developed by
B.~Cadorel [Cad19] might possibly be used in this context.

\claim 6.43. Remark|{\rm As mentioned in the course of the proof of
$(6.41 - 6.42)$ , we have neglected certain non negative terms coming 
from self-intersections $\Delta_j^p$
of the components ($p\ge 2$), by restricting the summation to the family 
of disjoint subsets $J_1,J_2,\ldots,J_k$. It would be interesting to refine
the lower bound and to take these terms into account. This might be possible
by observing that the iterated limit process, when integrating
on $\Delta_j$, involves inductively a few extra terms in (6.30),
when we take the limit as $t_{j,\varepsilon}\to\infty$. Those terms are 
equal to
$$
\langle\theta_{V,\gamma}\cdot u'_j,u'_j\rangle+
\big(\gamma\,\Theta_{A,h_{A,\delta}}-\rho_j^{-1}\,\Theta_{\Delta_j,h_j}\big)
\,|e_j^*(u)|^2+{\ii\over 2\pi}\big\langle b_j(u'_j),b_j(u'_j)\rangle_{h_j}.
$$
One would then have to evaluate the contribution of 
$\langle b_j(u'_j),b_j(u'_j)\rangle_{h_j}$ in the integral $\int_{\Delta_j}$.}
\endclaim

\plainsection{7. Non probabilistic estimates of the Morse integrals}

The non probabilistic estimate uses more explicit curvature inequalities and
has the advantage of producing results also in the general orbifold case.
Let us fix an ample line bundle $A$ on $X$ equipped with a smooth hermitian
metric $h_A$ such that $\omega_A:=\Theta_{A,h_A}>0$, and let
$\gamma_V$ be the infimum of values $\lambda\in\bR_+$ such that
$$
\lambda\,\omega_A\otimes\Id_V
-\Theta_{V,h_V}>_G 0,\leqno(7.1)
$$
in the sense of Griffiths. For any orbifold structure
$D=\sum_j(1-1/\rho_j)\Delta_j$, Corollary 6.16 then shows that the
$s$-th directed orbifold bundle
$V_s:=V\langle D^{(s)}\rangle$ (cf.\ \S$\,$1.B)
possesses hermitian metrics $h_{V\langle D^{(s)}\rangle,\varepsilon_s}$
such that the associated curvature tensor satisfies the inequality
$$
\theta_{s,\gamma,\varepsilon}
:=\gamma_s\,\omega_A\otimes\Id_{V\langle D^{(s)}\rangle}
-\Theta_{V\langle D^{(s)}\rangle,h_{V\langle D^{(s)}\rangle,\varepsilon_s}}>_G 0,
\leqno(7.2)
$$
provided we assume $d_jA-\Delta_j$ nef and take
$$
\leqalignno{
&\gamma_s > \gamma_{V,D^{(s)}}:=\max(\max_j(d_j/\rho_j^{(s)}),\gamma_V)
\quad\hbox{where}~~ \rho_j^{(s)}=\max(\rho_j/s,1).&(7.3)
\cr}
$$
In particular, any value
$$
\gamma_s > \max\big(s\,\max_j(d_j/\rho_j),\gamma_V\big).\leqno(7.3')
$$
is admissible, and we can apply the estimates 6.41~(b) and (6.42) with these
values. Instead of exploiting a Monte Carlo
convergence process for the curvature tensor as was done in \S$\,$4.B,
we are going to use a more precise lower bound of the
curvature tensor $\Theta_{L_{\tau,k},\varepsilon}$ of the orbifold
rank $1$ sheaf associated with $F=\tau A$, $\tau\ll 1$, namely
$$
L_{\tau,k}:=\cO_{X_k(V\langle D\rangle)}(1)\otimes
\pi_k^*\cO_X(-\tau A).
\leqno(7.4)
$$
Our formulas 3.20~(a,b) become
$$
\leqalignno{
&\Theta_{L_{\tau,k},\varepsilon}=\omega_{r,k,b}(\xi)+g_{k,0,\varepsilon}(z,x,u)
-\tau\,\omega_A(z),\quad
\hbox{where}&(7.5)\cr  
\qquad~~~&g_{k,\gamma,\varepsilon}(z,x,u)=\sum_{s=1}^k{x_s\over s}\,
\theta_{s,\gamma,\varepsilon}(u_s),&(7.5')\cr
&\theta_{s,\gamma,\varepsilon}(u_s)={\ii\over 2\pi}
\sum_{i,j,\lambda,\mu}c^{(s,\gamma,\varepsilon)}_{ij\lambda\mu}(z)\,
u_{s,\lambda}\overline u_{s,\mu}\,dz_i\wedge d\overline z_j.
&(7.5'')\cr}
$$
Under the assumption $(7.3')$, we have
$g_{k,\gamma,\varepsilon}(z,x,u)\ge 0$, but in general this is not true
for $g_{k,0,\varepsilon}(z,x,u)$, so we express
$g_{k,0,\varepsilon}(z,x,u)$ as a difference of
$g_{k,\gamma,\varepsilon}(z,x,u)$ and of a multiple of $\omega_A$.
By definition $\theta_{s,\gamma,\varepsilon}=
\gamma_s\,\omega_A\otimes\Id+\theta_{s,0,\varepsilon}$, and we get
$$
\leqalignno{\quad
&\Theta_{L_{\tau,k},\varepsilon}=\omega_{r,k,b}+\alpha_\varepsilon-\beta,
\qquad\hbox{where}
&(7.6)\cr
\noalign{\vskip6pt}
&\alpha_\varepsilon=g_{k,\gamma,\varepsilon}\ge 0,\quad \beta=
\bigg(\tau+\sum_{1\le s\le k}{\gamma_sx_s\over s}\bigg)\omega_A
=\sum_{1\le q\le k}{(\gamma_q+q\tau)\,x_q\over q}\,\omega_A\ge 0.
&(7.6')\cr}
$$
Then (7.6) and the inequalities used for (4.2), especially Lemma 2.3 and
Proposition 3.10~(b), lead to
$$
\leqalignno{\qquad
&\int_{X_k(V\langle D\rangle)(L_{\tau,k},\le 1)}
\Theta_{L_{\tau,k},\varepsilon}^{n+kr-1}&(7.7)\cr
&\quad{}={(n+kr-1)!\over n!\,k!^r(kr-1)!}\int_{z\in X}
\int_{(x,u)\in\bDelta^{k-1}\times(\bS^{2r-1})^k}
\bOne_{\alpha_\varepsilon-\beta,\le 1}\;(\alpha_\varepsilon-\beta)^n\,d\nu_{k,r}(x)\,d\mu(u)\cr
&\quad{}\ge{(n+kr-1)!\over n!\,k!^r(kr-1)!}\int_{z\in X}
\int_{(x,u)\in\bDelta^{k-1}\times(\bS^{2r-1})^k}
\big(\alpha_\varepsilon^n-n\alpha_\varepsilon^{n-1}\wedge\beta\big)\,
d\nu_{k,r}(x)\,d\mu(u).\cr}
$$
The main point is thus to find a lower bound of the difference
$\alpha_\varepsilon^n-n\alpha_\varepsilon^{n-1}\wedge\beta$, hence a lower
bound of $\alpha_\varepsilon^n$ and an upper bound of
$\alpha_\varepsilon^{n-1}\wedge\beta$.
An expansion of $\alpha_\varepsilon^n$ by Newton's multinomial formula yields
$$
\alpha_\varepsilon^n=\sum_{p\in\bN^k,\,|p|=n}{n!\over p_1!\ldots\,p_k!}\,
\prod_{s=1}^k\Big({x_s\over s}\,\theta_{s,\gamma,\varepsilon}(u_s)\Big)^{p_s}.
\leqno(7.8)
$$
If we assume $k\ge n$ and retain only the monomials for which $p_s=0,1$,
we get
$$
\alpha_\varepsilon^n\ge\sum_{1\le s_1<\ldots<s_n\le k}\,{n!\over s_1\ldots s_n}\,
\prod_{\ell=1}^n x_{s_\ell}\theta_{s_\ell,\gamma,\varepsilon}(u_{s_\ell}).
\leqno(7.8')
$$
By formula 3.10~(a) and an elementary calculation (cf.\ [Dem11, Prop.~1.13]),
one gets for every $(p_1,\ldots,p_k)\in\bN^k$
$$
\int_{\bDelta^{k-1}}x_1^{p_1}\ldots x_k^{p_k}\,
d\nu_{k,r}(x)={(kr-1)!\over (r-1)!^k}\;
{\prod_{1\le s\le k}(p_s+r-1)!\over(\sum_{1\le s\le k}p_s+kr-1)!},
\leqno(7.9)
$$
and in particular, for $k\ge n$, $p_1=\ldots=p_n=1$, $p_{n+1}=\ldots=p_k=0$,
we have
$$
\int_{\bDelta^{k-1}}x_{s_1}\ldots x_{s_n}\,d\nu_{k,r}(x)
=\int_{\bDelta^{k-1}}x_1\ldots x_n\,d\nu_{k,r}(x)
={(kr-1)!\,r^n\over(n+kr-1)!}.
\leqno(7.9')
$$
As a consequence, the equality case in $(6.20 - 6.21)$ implies
$$
\leqalignno{
M_{n,k,\varepsilon}\;\rlap{:}
&=\int_{z\in X}\int_{(x,u)\in\bDelta^{k-1}
\times(\bS^{2r-1})^k}
\alpha_\varepsilon(z)^n\,d\nu_{k,r}(x)\,d\mu(u_1)\ldots\,d\mu(u_k)\cr
&\ge\sum_{1\le s_1<\ldots<s_n\le k}
\int_{\bDelta^{k-1}}
{n!\,x_{s_1}\ldots x_{s_n}\over s_1\ldots s_n}\,d\nu_{k,r}(x)\times{}\cr
&\kern100pt\int_X\int_{\raise-1pt\hbox{${\textstyle\Pi}$}\kern1pt
S(V\langle D^{(s_\ell)}\rangle)}~\bigwedge_{\ell=1}^n
\langle\theta_{s_\ell,\gamma,\varepsilon}(u_{s_\ell}),u_{s_\ell}\rangle\,
d\mu(u_{s_\ell})\cr
&\ge\sum_{1\le s_1<\ldots<s_n\le k}{(kr-1)!\over(n+kr-1)!} \,
{n!\over s_1\ldots s_n}\times{}\cr
&\kern100pt\int_X~\prod_{\ell=1}^n\bigg(r\gamma_{s_\ell}\,c_1(A)-c_1(V)
+\sum_j(1-1/\rho^{(s_\ell)}_j)c_1(\Delta_j)\bigg)\cr
&\ge{(kr-1)!\over (n+kr-1)!} \,
\int_X~\prod_{s=1}^n\bigg(r\gamma_s\,c_1(A)-c_1(V)+\sum_j
(1-1/\rho^{(s)}_j)c_1(\Delta_j)\bigg).&(7.10_1)\cr}
$$
If we assume $c_1(V^*)=\lambda_V\,c_1(A)$ and $c_1(\Delta_j)=d_j\,c_1(A)$,
the lower bound takes the simpler form
$$
M_{n,k,\varepsilon}\ge{(kr-1)!\over(n+kr-1)!}\,
\prod_{s=1}^n\bigg(r\gamma_s+\lambda_V+\sum_j
d_j(1-1/\rho^{(s)}_j)\bigg)\,A^n.
\leqno(7.10_2)
$$
In fact, our lower bounds are obtained by taking into account the 
single term $s_\ell=\ell$, $1\le\ell\le k$ (which is the unique term in the sum 
when $k=n$). A more refined method is to integrate all monomials
$x_1^{p_1}\ldots\,x_k^{p_k}$ and to use the lower bound (6.42) instead of
$(6.20 - 6.21)$. This has the advantage of eventually producing a non zero
contribution, even when $k<n$. We find
$$
\eqalign{
M_{n,k}\;\rlap{:}
&=\lim_{\varepsilon\to 0}\int_{z\in X}\int_{(x,u)\in\bDelta^{k-1}
\times(\bS^{2r-1})^k}
\alpha_\varepsilon(z)^n\,d\nu_{k,r}(x)\,d\mu(u_1)\ldots\,d\mu(u_k)\cr
&\ge\lim_{\varepsilon\to 0}\sum_{\scriptstyle p\in\bN^k\atop\scriptstyle |p|=n}
\int_{\bDelta^{k-1}}
{n!\,x_1^{p_1}\ldots x_k^{p_k}\over\prod_{s=1}^kp_s!\,s^{p_s}}\,d\nu_{k,r}(x)
\int_X\int_{\raise-1pt\hbox{${\textstyle\Pi}$}\kern1pt
S(V\langle D^{(s)}\rangle)}\bigwedge_{s=1}^k
\langle\theta_{s,\gamma,\varepsilon}(u_s),u_s\rangle^{p_s}\,
d\mu(u_s)\cr}
$$
$$
\eqalign{
&\ge\sum_{\scriptstyle p\in\bN^k\atop\scriptstyle |p|=n}
{n!\over\prod_{s=1}^kp_s!\,s^{p_s}}\,
{(kr-1)!\over (r-1)!^k}\;
{\prod_{1\le s\le k}(p_s+r-1)!\over(\sum_{1\le s\le k}p_s+kr-1)!}~
\sum_{\scriptstyle J_1\amalg\ldots\amalg J_k\atop\scriptstyle
\subset\{1,\ldots,N\}}\cr
&\sum_{\scriptstyle
(q_j)\in(\bN^*)^{J_1\amalg\ldots\amalg J_k}
\atop\scriptstyle \Sigma_{j\in J_s}q_j\le p_s}
\prod_{1\le s\le k}{|J_s|!\,p_s!\,(p_s-\Sigma_{j\in J_s}q_j+r-1)!\over(p_s+r-1)!\,
(p_s-\Sigma_{j\in J_s}q_j)!}\,\prod_{j\in J_s}\bigg(1-{1\over \rho^{(s)}_j}\bigg)
\int_{z\in\Delta_{J_1\amalg\ldots\amalg J_k}}\cr
&\kern-20pt
\bigwedge_{1\le s\le k}\big((\gamma_s-\gamma_V)\Theta_{A,h_A}\big)^{p_s
-\Sigma_{j\in J_s}q_j}\wedge\bigwedge_{j\in J_s}{(r-1)!\over (q_j+r-1)!}\,
\big(\gamma_s\,\Theta_{A,h_A}-
(\rho^{(s)}_j)^{-1}\,\Theta_{\Delta_j,h_j}\big)^{q_j-1},\cr}
$$
thus
$$
\leqalignno{
M_{n,k}&\ge{n!\,(kr-1)!\over(n+kr-1)!}\sum_{\scriptstyle p\in\bN^k\atop
\scriptstyle |p|=n}~
\prod_{1\le s\le k}{1\over s^{p_s}}\,
\sum_{\scriptstyle J_1\amalg\ldots\amalg J_k\atop\scriptstyle
\subset\{1,\ldots,N\}}\cr
&\kern-20pt
\sum_{\scriptstyle (q_j)\in(\bN^*)^{J_1\amalg\ldots\amalg J_k}
\atop\scriptstyle \Sigma_{j\in J_s}q_j\le p_s}
\prod_{1\le s\le k}{|J_s|!\,(p_s-\Sigma_{j\in J_s}q_j+r-1)!\over
  (p_s-\Sigma_{j\in J_s}q_j)!}\,\prod_{j\in J_s}
\bigg(1-{1\over\rho^{(s)}_j}\bigg)
\int_{z\in\Delta_{J_1\amalg\ldots\amalg J_k}}~\bigwedge_{1\le s\le k}\cr
&~~~
\big((\gamma_s-\gamma_V)\Theta_{A,h_A}\big)^{p_s-\Sigma_{j\in J_s}q_j}
\wedge\bigwedge_{j\in J_s}{(r-1)!\over (q_j+r-1)!}\,
\big(\gamma_s\,\Theta_{A,h_A}-
(\rho^{(s)}_j)^{-1}\,\Theta_{\Delta_j,h_j}\big)^{q_j-1}.&(7.11)\cr}
$$
In particular, if $c_1(\Delta_j)=d_j\,c_1(A)$, we infer
$$
\leqalignno{
M_{n,k}&\ge{n!\,(kr-1)!\over(n+kr-1)!}\sum_{\scriptstyle p\in\bN^k\atop
\scriptstyle |p|=n}
\prod_{1\le s\le k}{1\over s^{p_s}}\kern-5pt
\sum_{\scriptstyle J_1\amalg\ldots\amalg J_k\atop\scriptstyle
\subset\{1,\ldots,N\}}\sum_{\scriptstyle (q_j),\,q_j\ge 1\atop\scriptstyle
\Sigma_{j\in J_s}q_j\le p_s}\!\!\prod_{1\le s\le k}\!\Bigg(
{|J_s|!\,(p_s-\Sigma_{j\in J_s}q_j+r-1)!\over
(p_s-\Sigma_{j\in J_s}q_j)!}\cr
&\kern40pt
(\gamma_s-\gamma_V)^{p_s-\Sigma_{j\in J_s}q_j}\,
\prod_{j\in J_s}d_j\bigg(1-{1\over \rho^{(s)}_j}\bigg)\,
{(r-1)!\over (q_j+r-1)!}
\bigg(\gamma_s-{d_j\over\rho^{(s)}_j}\bigg)^{q_j-1}\Bigg)\,A^n.&(7.12)\cr}
$$
In the special case $k=1$ and $N\ge n$, by
taking $|J|=|J_1|=n$ and $q_j=1$ for all $j\in J$, we find
$$
M_{n,1}\ge{n!\,(r-1)!\over(n+r-1)!}
\sum_{J\subset\{1,\ldots,N\},\,|J|=n}
{n!\,(r-1)!\over r^n}\,
\prod_{j\in J}d_j\bigg(1-{1\over \rho_j}\bigg)\,A^n.
\leqno(7.12_1)
$$
Next, we turn ourselves to the evaluation of the integral of
$\alpha_\varepsilon^{n-1}\wedge\beta$. We have
$$
\alpha_\varepsilon^{n-1}\wedge\beta=\sum_{p\in\bN^k,\,|p|=n-1}
{(n-1)!\over p_1!\ldots\,p_k!}\,
\prod_{s=1}^k\bigg({x_s\over s}\,\theta_{s,\gamma,\varepsilon}(u_s)\bigg)^{p_s}
\wedge\beta,\leqno(7.13)
$$
and the upper bound given by $(6.20 - 6.21)$ provides
$$
\leqalignno{
M'_{n,k}:&=\lim_{\varepsilon\to 0}\int_{z\in X}
\int_{(x,u)\in\bDelta^{k-1}\times(\bS^{2r-1})^k}
n\,\alpha_\varepsilon(z)^{n-1}\wedge\beta\,
d\nu_{k,r}(x)\,d\mu(u_1)\ldots\,d\mu(u_k)\cr
&\le\lim_{\varepsilon\to 0}\sum_{p\in\bN^k,\,|p|=n-1}\int_{\bDelta^{k-1}}
n\,{(n-1)!\,x_1^{p_1}\ldots x_k^{p_k}\over
\prod_{s=1}^kp_s!\,s^{p_s}}\,d\nu_{k,r}(x)\times{}\cr
&\kern90pt\int_X\int_{\raise-1pt\hbox{${\textstyle\Pi}$}\kern1pt
S(V\langle D^{(s)}\rangle)}\bigwedge_{s=1}^k
\langle\theta_{s,\gamma,\varepsilon}(u_s),u_s\rangle^{p_s}\wedge\beta\,
\prod_{s=1}^kd\mu(u_s)\cr
&\le\sum_{p\in\bN^k,\,|p|=n-1}\int_{\bDelta^{k-1}}
n\,{(n-1)!\,x_1^{p_1}\ldots x_k^{p_k}\over\prod_{s=1}^kp_s!\,s^{p_s}}\,
\bigg(\sum_{q=1}^k{(\gamma_q+q\tau)\,x_q\over q}\bigg)\,d\nu_{k,r}(x)~\times\cr
&~\int_X~\bigwedge_{1\le s\le k}{p_s!\,(r-1)!\over(p_s+r-1)!}\,
\bigg(r\gamma_s\,\Theta_{A,h_A}-\Tr\Theta_{V,h_V}+\sum_j
(1-1/\rho^{(s)}_j)\Theta_{\Delta_j,h_j}\bigg)^{p_s}\wedge\Theta_{A,h_A}.\cr}
$$
By (7.9), for $|p|=\sum p_s=n-1$, we get
$$
\eqalign{
&\int_{\bDelta^{k-1}}x_1^{p_1}\ldots x_k^{p_k}\bigg(\sum_{q=1}^k{\gamma_q\over q}
+\tau x_q\bigg)d\nu_{k,r}(x)\cr
&\kern100pt{}={(kr-1)!\over (r-1)!^k}\;
{\prod_{1\le s\le k}(p_s+r-1)!\over(n-1+kr-1)!}
\bigg(\sum_{q=1}^k{\gamma_q\over q}+\tau\sum_{q=1}^k{p_q+r\over n+kr-1}\bigg)\cr
&\kern100pt{}={(kr-1)!\over (r-1)!^k}\;
{\prod_{1\le s\le k}(p_s+r-1)!\over(n+kr-2)!}
\bigg(\sum_{q=1}^k{\gamma_q\over q}+\tau\bigg).\cr}
$$
Therefore, assuming $c_1(\Delta_j)=d_j\,c_1(A)$ and
$c_1(V^*)=\lambda_Vc_1(A)$, we find
$$
\leqalignno{
M'_{n,k}&\le {n!\,(kr-1)!\over (r-1)!^k\,(n+kr-2)!}\,
\bigg(\sum_{q=1}^k{\gamma_q\over q}+\tau\bigg)
\sum_{p\in\bN^k,\,|p|=n-1}{\prod_{1\le s\le k}(p_s+r-1)!\over
\prod_{s=1}^kp_s!\,s^{p_s}}~\times\cr
&\kern25pt
\prod_{1\le s\le k}{p_s!\,(r-1)!\over(p_s+r-1)!}\,
\Big(r\gamma_s+\lambda_V+\sum_jd_j(1-1/\rho^{(s)}_j)\Big)^{p_s}\,A^n,
\cr
&\le {n!\,(kr-1)!\over(n+kr-2)!}\,
\Bigg(\sum_{q=1}^k{\gamma_q\over q}+\tau\Bigg)\times{}\cr
&\kern80pt
\sum_{p\in\bN^k,\,|p|=n-1}~
\prod_{1\le s\le k}{1\over s^{p_s}}\,\Big(r\gamma_s
+\lambda_V+\sum_jd_j(1-1/\rho^{(s)}_j)\Big)^{p_s}\,A^n.&(7.14_1)\cr}
$$
A simpler (but larger) upper bound is
$$
\leqalignno{
&&(7.14_2)\cr
&M'_{n,k}\le{n!\,(kr-1)!\over(n+kr-2)!}\,
\Bigg(\sum_{s=1}^k{\gamma_s\over s}+\tau\Bigg)
\Bigg(\sum_{1\le s\le k}{1\over s}\,\Big(r\gamma_s+
\lambda_V)+\sum_jd_j(1-1/\rho^{(s)}_j)\Big)\Bigg)^{n-1}\,A^n.\cr}
$$
Finally, inequality (7.7) translates into
$$
{1\over (n+kr-1)!}\int_{X_k(V\langle D\rangle)(L_{\tau,k},\le 1)}
\Theta_{L_{\tau,k},\varepsilon}^{n+kr-1}
\ge{1\over n!\,k!^r(kr-1)!}(M_{n,k}-M'_{n,k}).
\leqno(7.15)
$$
If we put everything together, we get the following (complicated!)
existence criterion for orbifold jet differentials.

\claim 7.16.~Existence criterion|Let $(X,V,D)$ with
$D=\sum_{1\le j\le N}(1-1/\rho_j)\Delta_j$ be a directed orbifold,
and let $A$ be an ample line bundle on~$X$. Assume that $D$ is a simple normal
crossing divisor transverse to $V$, that $c_1(\Delta_j)=d_j\,c_1(A)$,
$c_1(V^*)=\lambda_V\,c_1(A)$ and let $\gamma_V$ be the infimum
of values $\gamma>0$ such that $\Theta_A\otimes\Id_V-\Theta_V\ge_G 0$.
Take
$$
\gamma_s=\max(\max(d_j/\rho_j^{(s)}),\gamma_V),\quad
\rho_j^{(s)}=\max(\rho_j/s,1).
$$
Then, a sufficient condition for
the existence of $($many$)$ non zero holomorphic sections of multiples of
$$
L_{\tau,k}=\cO_{X_k(V\langle D\rangle)}(1)\otimes\pi_k^*\cO(-\tau A)
$$
on $X_k(V\langle D\rangle)$ is that $M_{n,k}-M'_{n,k}>0$, where $M_{n,k}$
admits the lower bounds $(7.10_2)$ or $(7.12)$, and $M'_{n,k}$ admits the
upper bound $(7.14_2)$.
\endclaim

\plainsubsection 7.B. Compact case (no boundary divisor)|

We address here the case of a compact (projective) directed manifold $(X,V)$,
with a boundary divisor $D=0$. By $(7.10_2)$ and $(7.14_2)$, we find
$$
\eqalign{
&M_{n,k}\ge {(kr-1)!\over(n+kr-1)!}(r\gamma_V+\lambda_V)^n\,A^n\quad
\hbox{if $k\ge n$},\cr
&M'_{n,k}\le {n!\,(kr-1)!\over(n+kr-1)!}
\Bigg(\tau+\gamma_V\sum_{s=1}^k{1\over s}\Bigg)\Bigg(\sum_{s=1}^k
{1\over s}\big(r\gamma_V+\lambda_V\big)\Bigg)^{n-1}.\cr}
$$
Therefore, for $\tau>0$ sufficiently small, $M_{n,k}-M'_{n,k}$ is positive as
soon as $k\ge n$ and $(r\gamma_V+\lambda_V)^n>n!\,\gamma_V(
\sum_{1\le s\le k}{1\over s})^n(r\gamma_V+\lambda_V)^{n-1}$, that is
$$
k\ge n\quad\hbox{and}\quad
\lambda_V>n!\,\Bigg(\sum_{1\le s\le k}{1\over s}\Bigg)^n\gamma_V-r\gamma_V.
\leqno(7.17)
$$

\claim 7.18. Example|{\rm 
In the case where $X$ is a smooth hypersurface of $\bP^{n+1}$ of
degree $d$ and~$V=T_X$, we have $r=n$ and $\det(V^*)=\cO(d-n-2)$.
We take $A=\cO(1)$. If $Q$ is the tautological quotient bundle on
$\bP^{n+1}$, it is well known that $T_{\bP^{n+1}}\simeq Q\otimes\cO(1)$ and
$\det Q=\cO(1)$, hence 
$T^*_{\bP^{n+1}}\otimes\cO(2)=Q^*\otimes \cO(1)=\Lambda^nQ\ge_G0$, and
the surjective morphism
$$
T^*_{\bP^{n+1}|X}\to T^*_X=V^*
$$
implies that we also have $V^*\otimes\cO(2)\ge_G 0$. Therefore,
we find $\gamma_V=2$ and $\lambda_V=d-n-2$. The above condition
(7.17) becomes $k\ge n$ and
$$
k\ge n\quad\hbox{and}\quad
d>2\,n!\,\Bigg(\sum_{1\le s\le k}{1\over s}\Bigg)^n-n+2.
$$
This lower bound improves the one stated in [Dem12], but is
unfortunately far from being optimal. Better bounds -- still probably non
optimal -- have been obtained in [Dar16] and [MTa19].}
\endclaim

\plainsubsection 7.C. Logarithmic case|

The logarithmic situation makes essentially no difference in treatment
with the compact case, except for the fact that we have to replace $V$
by the logarithmic directed structure $V\langle D\rangle$, and the numbers
$\gamma_V$, $\lambda_V$ by
$$
\leqalignno{
&\gamma_{V\langle D\rangle}=\inf\gamma~~\hbox{such that}~~
\gamma\Theta_A-\Theta_{V\langle D\rangle}\ge_G 0,&(7.19_1)\cr
&\lambda_{V\langle D\rangle}~~\hbox{such that}~~
c_1(V^*\langle D\rangle)=\lambda_{V\langle D\rangle}\,c_1(A)\quad
\hbox{(if such $\lambda_{V\langle D\rangle}$ exists)}.
&(7.19_2)\cr}
$$
We get the sufficient condition
$$
k\ge n\quad\hbox{and}\quad
\lambda_{V\langle D\rangle}>n!\,\Bigg(\sum_{1\le s\le k}{1\over s}
\Bigg)^n\gamma_{V\langle D\rangle}-r\gamma_{V\langle D\rangle}.
\leqno(7.20)
$$
For $X=\bP^n$, $V=T_{\bP^n}$, and for a divisor $D=\sum\Delta_j$ of total
degree~$d$ on $\bP^n$, we can still take $\gamma_{V\langle D\rangle}=2$ by
Proposition 5.8, and we have $\det(V^*\langle D\rangle)=\cO(d-n-1)$.
We get the degree condition
$$
k\ge n\quad\hbox{and}\quad
d>2\,n!\,\Bigg(\sum_{1\le s\le k}{1\over s}\Bigg)^n-n+1.
\leqno(7.21)
$$
Again, [Dar16] and [MTa19] gave better bounds for this particular
logarithmic situation.

\plainsubsection 7.D. Case of orbifold structures on projective $n$-space|

Let us come to our main target, namely ``genuine'' orbifolds, for which
our results are completely new. The situation
we have in mind is the case of triples $(X,V,D)$ where $X=\bP^n$, $V=T_X$,
$D=\sum (1-1/\rho_j)\Delta_j$ is a normal crossing divisor,
with components $\Delta_j$ of degree~$d_j$. Set again $A=\cO(1)$. Since
$c_1(V^*)=-(n+1)\,c_1(A)$ and $D^{(s)}=\sum_j(1-s/\rho_j)_+\Delta_j$, we have
$$
\lambda_V=-n-1,\quad \det V^*\langle D^{(s)}\rangle
=\cO_{\bP^n}\Big(-n-1+\sum_j d_j(1-s/\rho_j)_+\Big).
\leqno(7.22)
$$
Moreover, by Proposition~5.8, we get
$$
\Theta_{V^*\langle D^{(s)}\rangle}+\gamma_s\,\omega_\FS\otimes\Id>_G0
\leqno(7.23)
$$
as soon as $\gamma_s>2$ and
$\gamma_s>\max_j(d_j/\max(\rho_j/s,1))$ for all components $\Delta_j$
in $D^{(s)}$. We can take for instance
$\gamma_s>st$ where $t=\max(\max_j(d_j/\rho_j),2)$. By considering the
infimum and applying $(7.10_2)$ when $r=n$ and $k\ge n$, we find
$$
M_{n,k,\varepsilon}\ge{(kn-1)!\over(n+kn-1)!}\,
\prod_{s=1}^n\Big(ns\,t-n-1+\sum_j
d_j(1-s/\rho_j)_+\Big)\,A^n,
\leqno(7.24)
$$
while $(7.14_2)$ implies
$$
M'_{n,k}\le{n!\,(kn-1)!\over(n+kn-2)!}\,
(kt+\tau)
\Bigg(\sum_{1\le s\le k}{1\over s}\bigg(ns\,t-n-1+\sum_jd_j(1-s/\rho_j)_+
\bigg)\Bigg)^{n-1}\,A^n.
\leqno(7.25)
$$
If we take $\rho_j\ge\rho>n$, then $(1-s/\rho_j)_+\ge 1-s/\rho$ for $s\le n$,
and as $ns\,t-n-1\ge 0$ and $\sum_{1\le s\le k}{1\over s}(nst-n-1)\le nkt$,
we get for $\tau>0$ small a sufficient condition
$$
\prod_{s=1}^n\Bigg(\Big(1-{s\over \rho}\Big)\sum_jd_j\Bigg)>
kt\,(n+kn-1)\,n!\,
\Bigg(nk\,t+\Big(1+{1\over 2}+\cdots+{1\over k}\Big)\sum_jd_j\Bigg)^{n-1}.
$$
For $k=n$, the latter condition is satisfied if
$\sum_j d_j>c_nt\,\prod_{s=1}^n\big(1-{s\over\rho}\big)^{-1}$ with
$$
c_n=n(n^2+n-1)\,n!\,
\Big(1+{1\over 2}+\cdots+{1\over n}+{1\over n^3}\Big)^{n-1}.
\leqno(7.26)
$$
In fact, $c_1=1$, $c_2=32.5$ and $c_n\ge n^5$ for all $n\in\bN^*$, hence the above requirement
implies in any case the inequality $n^2t\le {1\over n^3}\sum d_j$.
The Stirling and Euler-Maclaurin formulas give
$$
c_n\sim (2\pi)^{1/2} n^{n+7/2}\,e^{-n}\,(\gamma+\log n)^{n-1}
\leqno(7.26')
$$
as $n\to+\infty$, where $\gamma=0.577215\ldots$ is the Euler constant, the
ratio being actually bounded above for $n\ge 3$ by
$\exp\big((1/2)(1-1/n)/(\gamma+\log n)+13/12n-1/n^2\big)\to 1$.
Let us observe that
$$
{1\over t}=\min\bigg(\!\min_j\bigg({\rho_j\over d_j}\bigg),{1\over 2}\bigg).
$$
In this way, we get the sufficient condition
$$
\rho_j\ge \rho>n,\quad
\sum_j d_j\cdot
\min\bigg(\!\min_j\bigg({\rho_j\over d_j}\bigg),{1\over 2}\bigg)
\,\prod_{s=1}^n\Big(1-{s\over\rho}\Big)>c_n.\leqno(7.27)
$$
For instance, if we take all components $\Delta_j$ possessing
the same degrees $d_j=d$ and ramification number $\rho_j\ge \rho$,
these numbers and
the number $N$ of components have to satisfy the sufficient condition
$$
\rho>n,\quad 
N\min(\rho,d/2)\,\prod_{s=1}^n\Big(1-{s\over\rho}\Big)>c_n.
\leqno(7.27_N)
$$
This possibly allows a single component (taking $d,\rho$ large), or
$d,\rho$ small (taking $N$ large).
Since we have neglected many terms in the above calculations, 
the ``technological constant'' $c_n$ appearing in these
estimates is probably much larger than needed. Notice that
the above estimates require jets of order $k\ge n$ and ramification
numbers $\rho>n$. Parts (${\rm a}$) and (${\rm a}'$) of Theorem~0.8
follow from (7.27) and $(7.27_N)$.
\medskip

\noindent
{\bf 7.28. Case of jet differentials of order k = 1 (symmetric differentials)}.
When $k<n$ or $\rho_j\in{}]1,+\infty]$, estimate
(7.12) still allows us to obtain an existence criterion. For instance, when
$k=1$ and $N\ge n$, $(7.12_1)$ and (7.25) give
$$
\eqalign{
M_{n,1}&\ge{n!\,(n-1)!\over(2n-1)!}
\sum_{J\subset\{1,\ldots,N\},\,|J|=n}
{n!\,(n-1)!\over n^n}\,
\prod_{j\in J}d_j\bigg(1-{1\over \rho_j}\bigg)\,A^n,\cr
M'_{n,1}&\le {n!\,(n-1)!\over(2n-2)!}\,
(t+\tau)
\Big(n\,t-n-1+\sum_jd_j(1-1/\rho_j)\Big)^{n-1}\,A^n,\cr}
$$
and we get the non void existence criterion
$$
\sum_{J\subset\{1,\ldots,N\},\,|J|=n}~
\prod_{j\in J}d_j\bigg(1-{1\over \rho_j}\bigg)>
(2n-1)\,t\,
\Big(n\,t-n-1+\sum_jd_j(1-1/\rho_j)\Big)^{n-1}
\leqno(7.29)
$$
where $t=\max(\max_j(d_j/\rho_j),2)$. For instance, if all divisors have the
same degrees $d_j=d$ and ramification numbers $\rho_j\ge\rho$, condition (7.29) 
is implied by
$$
{N\choose n}\,d^n\bigg(1-{1\over \rho}\bigg)^n>
(2n-1) \max(d/\rho,2)\,\big((N+n)d\big)^{n-1},
$$
or equivalently, by
$$
\min(\rho,d/2)\,{N\choose n}
\bigg(1-{1\over \rho}\bigg)^n>(2n-1)\,(N+n)^{n-1}.
\leqno(7.29')
$$
As $j\mapsto(N-j)/(n-j)$ is non decreasing for $0\le j<n\le N$, we have
the inequality ${N \choose n}=\prod_{0\le j<n}{N-j\over n-j}\ge (N/n)^n$,
hence
$$
{\displaystyle{N\choose n}\over (2n-1)(N+n)^{n-1}}
\ge {N^n\over n^n(2n-1)(2N)^{n-1}}={N\over 2^{n-1}\,(2n-1)\,n^n}.
$$
We finally get the sufficient condition
$$
N\ge n,\quad
N\min(\rho,d/2)\,\,\bigg(1-{1\over \rho}\bigg)^n>2^{n-1}\,(2n-1)\,n^n.
\leqno(7.29_N)
$$
Parts (${\rm b}$) and (${\rm b}'$) of Theorem~0.8
follow from (7.29) and $(7.29_N)$. Again, the constant $2^{n-1}\,(2n-1)\,n^n$
is certainly far from being optimal. Answering the problem raised in
Remark~6.43 might help to improve the bounds.

\plainsection{8. Appendix: a proof of the orbifold vanishing theorem}

The orbifold vanishing theorem is proved in [CDR20] in the case of boundary
divisors $D=\sum(1-1/\rho_j)\Delta_j$ with rational multiplicities
$\rho_j\in{}]1,\infty]$. However, the definition of
orbifold curves shows that we can replace $\rho_j$ by
$\lceil\rho_j\rceil\in\bN\cup\{\infty\}$ without
modifying the space of curves we have to deal with. On the other
hand, this replacement makes the corresponding sheaves
$E_{k,m}V^*\langle D\rangle$ larger.
Therefore, the case of arbitrary real multiplicities $\rho_j\in{}]1,\infty]$
stated in Proposition~0.7 follows from the case of integer multiplicities.
We sketch here an alternative and possibly more direct proof of
Proposition~0.7, by checking that we can still apply the Ahlfors-Schwarz
lemma argument of [Dem97] in the orbifold context. For this, we associate
to $D$ the ``logarithmic divisor''
$$
D'=\lceil D\rceil = \sum\Delta_j \ge D,
$$
and, assuming $(X,V,D')$ non singular, we make use of the tower of
logarithmic Semple bundles
$$
X^\rS_k(V\langle D'\rangle)\to
X^\rS_{k-1}(V\langle D'\rangle)\to\cdots\to
X^\rS_1(V\langle D'\rangle)\to X^\rS_0(V\langle D'\rangle):=X
\leqno(8.1)
$$
(in reference to the work of the British mathematician John Greenlees
Semple, see [Sem54]), where each stage is a smooth directed manifold
$(X^\rS_k(V\langle D'\rangle),V_k\langle D'\rangle)$ defined inductively by
$$
X^\rS_k(V\langle D'\rangle):=P(V_{k-1}\langle D'\rangle)=
\hbox{projective bundle of lines of~$V_{k-1}\langle D'\rangle$},
\leqno(8.2)
$$
and $V_k\langle D'\rangle$ is a subbundle of the logarithmic tangent bundle
of $X^\rS_k(V\langle D'\rangle)$ associated with the pull-back of $D'$.
Each of these projective bundles is equipped with a tautological line bundle
$\cO_{X^\rS_k(V\langle D'\rangle)}(-1)$ (see [Dem97] for details),
and $V_k\langle D'\rangle$ consists of the elements of the logarithmic
tangent bundle that project onto the tautological line, so that we
have an exact sequence
$$
0\to T_{X^\rS_k(V\langle D'\rangle)/X^\rS_{k-1}(V\langle D'\rangle)}
\to V_k\langle D'\rangle\to \cO_{X^\rS_k(V\langle D'\rangle)}(-1)\to 0.
\leqno(8.2')
$$
We let $\pi_{k,\ell}:X^\rS_k(V\langle D'\rangle)\to
X^\rS_\ell(V\langle D'\rangle)$ be the natural projection.
Then the top-down projection $\pi_{k,0}:X^\rS_k(V\langle D'\rangle)\to X$
yields a direct image sheaf
$$
(\pi_{k,0})_*\cO_{X^\rS_k(V\langle D'\rangle)}(m):=E^\rS_{k,m}V^*\langle D'\rangle
\subset E_{k,m}V^*\langle D'\rangle.
\leqno(8.3)
$$
Its stalk at point $x\in X$ consists of the algebraic differential 
operators $P(f_{[k]})$ acting on germs
of $k$-jets $f:(\bC,0)\to(X,x)$ tangent to $V$, satisfying the 
invariance property
$$
P((f\circ\varphi)_{[k]})=
(\varphi')^m P(f_{[k]})\circ \varphi,
\leqno(8.4)
$$
whenever $\varphi\in \bG_k$ is in the group of $k$-jets of biholomorphisms
$\varphi:(\bC,0)\to(\bC,0)$. By construction, the sheaf of
orbifold jet differentials $E_{k,m}V^*\langle D\rangle$ is
contained in $E_{k,m}V^*\langle D'\rangle$, and we have a corresponding
inclusion
$$
E^\rS_{k,m}V^*\langle D\rangle\subset E^\rS_{k,m}V^*\langle D'\rangle
\leqno(8.5)
$$
of the Semple orbifold jet differentials into the Semple
logarithmic differentials. A consideration of the algebra
$\bigoplus E^\rS_{k,m}V^*\langle D\rangle$ makes clear that there
exists a submultiplicative sequence of ideal sheaves $(\cJ_{D,k,m})_{m\in \bN}$
on $X^\rS_k(V\langle D'\rangle)$, such that
the image of $\pi_{k,0}^*\cO_X(E^\rS_{k,m}V^*\langle D\rangle)$
in $\cO_{X^\rS_k(V\langle D'\rangle)}(m)$ is a sheaf
$$
\cO_{X^\rS_k(V\langle D'\rangle)}(m)\otimes\cJ_{D,k,m}.
\leqno(8.6)
$$
It is clear that the zero variety of $V(\cJ_{D,k,m})$ projects into
the support $|D'|=|D|$ of~$D$. We consider a smooth log resolution
$$
\mu_k:\tilde X_k\to X^\rS_k(V\langle D'\rangle)
\leqno(8.7)
$$
of the ideal $\cJ_{D,k,m}$  in $X^\rS_k(V\langle D'\rangle)$, so that
$\mu_k^*(\cJ_{D,k,m})=\cO_{\tilde X_k}(-G_{D,k,m})$ for a suitable
effective simple normal crossing divisor $G_{D,k,m}$ on $\tilde X_k$
that projects into $|D|$ in~$X$. Denoting
$\cO_{\tilde X_k}(1)=\mu_k^*\cO_{X^\rS_k(V\langle D'\rangle)}(1)$, we get
$$
\mu_k^*\big(\cO_{X^\rS_k(V\langle D'\rangle)}(m)\otimes\cJ_{D,k,m}\big)=
\cO_{\tilde X_k}(m)\otimes\cO_{\tilde X_k}(-G_{D,k,m}).
\leqno(8.7')
$$
We denote by $\tilde\pi_{k,\ell}$ the composition
$$
\tilde\pi_{k,\ell}=\pi_{k,\ell}\circ\mu_k:\tilde X_k\to X^\rS_k(V\langle D'\rangle)
\to X^\rS_\ell(V\langle D'\rangle),
$$
and consider especially the projection $\tilde\pi_{k,0}:\tilde X_k\to X$.
For every entire or local orbifold entire curve
$f:\bC\supset\Omega\to(X,V,D)$, the
image $f(\Omega)$ is not entirely contained in $|D'|$, and we thus get
holomorphic $k$-jet liftings
$$
f_{[k]}:\Omega\to X^\rS_k(V\langle D'\rangle)\quad\hbox{and}\quad
\tilde f_{[k]}:\Omega\to \tilde X_k.
\leqno(8.8)
$$
Morevover, the derivative $f'_{[k-1]}$ of the $(k-1)$-jet lifting
$f_{[k-1]}$ can be seen as a meromorphic section of the logarithmic
tautological line bundle $(f_{[k]})^*\cO_{X^\rS_k(V\langle D'\rangle)}(-1)$,
since the multipli\-cities of zeroes of $f'_{[k-1]}$ are possibly less than
the ones prescribed by the logarithmic condition.
The poles are of course contained in $f^{-1}(|D'|)$. As a
consequence, $f'_{[k-1]}$ also lifts as a meromorphic section of
$(\tilde f_{[k]})^*\cO_{\tilde X_k}(-1)$, which we denote by
$\tilde f'_{[k-1]}$. If $\tau_{D'}\in H^0(X,\cO_X(D'))$ is the canonical
section of divisor equal to $D'$, we get at worst that
$$
\leqalignno{
&\tau_{D'}(f)\,f'_{[k-1]}\in H^0\big(\Omega,
(f_{[k]})^*(\cO_{X^\rS_k(V\langle D'\rangle)}(-1)\otimes
\pi_{k,0}^*\cO_X(D'))\big)
\quad\hbox{and}&(8.9)\cr
&\tau_{D'}(f)\,\tilde f'_{[k-1]}\in H^0\big(\Omega,
(\tilde f_{[k]})^*(\cO_{\tilde X_k}(-1)\otimes\tilde\pi_{k,0}^*\cO_X(D'))\big)
&(8.9\tilde{~})\cr}
$$
are holomorphic. On the other hand, every local section
$P\in H^0\big(U,E^\rS_{k,m}V^*\langle D\rangle\big)$ on
an open subset $U\subset X$ gives rise
in a one-to-one manner to a section
$$
\sigma_P\in H^0\big(U_k,\cO_{X^\rS_k(V\langle D'\rangle)}(m)
\otimes \cJ_{D,k,m}\big),\quad
U_k=\pi_{k,0}^{-1}(U)\subset X^\rS_k(V\langle D'\rangle),
$$
by the correspondence
$$
P(f_{[k]})=\sigma_P(f_{[k]})\cdot (f'_{[k-1]})^m
\leqno(8.10)
$$
for every local orbifold curve $f$ contained in $U$. By pulling back
to $\tilde X_k$, we get a section
$$
\tilde\sigma_P\in H^0\big(\tilde U_k,\cO_{\tilde X_k}(m)
\otimes \cO_{\tilde X_k}(-G_{D,k,m})\big),
\quad \tilde U_k=\mu_k^{-1}(U_k)=\tilde\pi_{k,0}^{-1}(U),
$$
such that
$$
P(f_{[k]})=\tilde\sigma_P(\tilde f_{[k]})\cdot (\tilde f'_{[k-1]})^m.
\leqno(8.10\tilde{~})
$$
However, $P(f_{[k]})$ is a holomorphic function, and we must have
a cancellation of the poles of $(\tilde f'_{[k-1]})^m$ for all sections
$\tilde\sigma_P$, which generate the sheaf
$\cO_{\tilde X_k}(m)\otimes \cO_{\tilde X_k}(-G_{D,k,m})$.
This means that
$$
\hbox{%
$\tilde f'_{[k-1]}$ is a holomorphic section of 
$(\tilde f_{[k]})^*\cO_{\tilde X_k}(-1)\otimes
\cO_\bC\big(\lfloor{1\over m}(\tilde f_{[k]})^*G_{D,k,m}\rfloor\big)$}
\leqno(8.11)
$$
For any given ample divisor $A$ over $X$, we can find $s=s_{k,m}\in\bN^*$ such
that the tensor product
$\cO_X(E^\rS_{k,m}V^*\langle D\rangle)\otimes\cO_X(sA)$ is generated by
its global sections over $X$. By taking the pull-back to $\tilde X_k$
and looking at the image in $\cO_{\tilde X_k}(m)$, we conclude that
$$
\cO_{\tilde X_k}(m)\otimes\cO_{\tilde X_k}(-G_{D,k,m})\otimes
\tilde\pi_{k,0}^*\cO_X(sA)\quad\hbox{is generated by sections
on $\tilde X_k$}.
\leqno(8.12)
$$
As in [Dem97], let us consider for every weight
$\underline a=(a_1,\ldots,a_k)\in\bZ^k$ the line bundles
$$
\cO_{X^\rS_k(V\langle D'\rangle)}(\underline a)=
\bigotimes_{1\le\ell\le k}\pi_{k,\ell}^*\cO_{X^\rS_\ell(V\langle D'\rangle)}(a_\ell),
\qquad
\cO_{\tilde X_k}(\underline a)=\mu_k^*
\cO_{X^\rS_k(V\langle D'\rangle)}(\underline a).
\leqno(8.13)
$$
Since each factor $\cO_{X^\rS_\ell(V\langle D'\rangle)}(1)$ is relatively
ample with respect to $\pi_{\ell,\ell-1}$, it is easy to see by induction
on $k$ that thee exists a weight $\underline a\in(\bN^*)^k$ and $b\in\bN^*$
such that the line bundle
$\cO_{X^\rS_k(V\langle D'\rangle}(\underline a)\otimes\pi_{k,0}^*\cO_X(bA)$
is~ample. After possibly replacing $(\underline a,b)$ by a multiple,
we can find a $\mu_k$-exceptional divisor $H_{D,k}$ on $\tilde X_k$ such that
$$
\cO_{\tilde X_k}(\underline a)\otimes
\cO_{\tilde X_k}(-H_{D,k})\otimes\tilde\pi_{k,0}^*\cO_X(bA)
\leqno(8.14)
$$
is very ample on $\tilde X_k$. Finally, we select $c\in\bN^*$ such that
$$
\cO_X(cA-D')~~\hbox{is very ample on $X$}.\leqno(8.15)
$$
By taking the tensor product of $(8.12-8.15)$, (8.15) being raised to a
power $t\in\bN^*$, we find that
$$
L_{k,m}:=\cO_{\tilde X_k}(m)\otimes\cO_{\tilde X_k}(\underline a)\otimes
\cO_{\tilde X_k}(-G_{D,k,m}-H_{D,k})\otimes
\tilde\pi_{k,0}^*\cO_X((s+b+tc)A-tD')
\leqno(8.16)
$$
is very ample on $\tilde X_k$. We will later need to take
$t=|\underline a|=\sum_\ell a_\ell$, which is of course an admissible choice.

\claim 8.17. Lemma|Let $(X,V,D)$ be a projective non singular
directed orbifold, and $A$~an ample divisor on $X$. Then, for every
orbifold entire curve $f:\bC\to(X,V,D)$ and every section
$$
P\in H^0\big(X,E^\rS_{k,m}V^*\langle D\rangle\otimes\cO_X(-A)\big),
$$
we have $P(f_{[k]})=P(f,f',\ldots,f^{(k)})=0$.
\endclaim

\plainproof. As we have already seen for local sections, every global jet
differential $P$ in $H^0\big(X,E^\rS_{k,m}V^*\langle D\rangle
\otimes\cO_X(-A)\big)$ gives rise to sections
$$
\eqalign{
&\sigma_P\in H^0\big(X^\rS_k(V\langle D'\rangle),\cO_{X^\rS_k(V\langle D'\rangle)}(m)
\otimes \cJ_{D,k,m}\otimes\pi_{k,0}^*\cO_X(-A)\big),\cr
&\tilde\sigma_P\in H^0\big(\tilde X_k,\cO_{\tilde X_k}(m)
\otimes \cO_{\tilde X_k}(-G_{D,k,m})\otimes\tilde\pi_{k,0}^*\cO_X(-A)\big)\cr}
$$
such that
$$
P(f_{[k]})=\tilde\sigma_P(\tilde f_{[k]})\cdot (\tilde f'_{[k-1]})^m
\in H^0(\bC,f^*\cO_X(-A)).
$$
Assume that $P(f_{[k]})\ne 0$ (so that, in particular $\tilde\sigma_P\ne 0$).
We consider a basis $(g_j)$ of sections of $L_{k,m}$ in (8.16), the canonical
section $\eta_{D,k}\in H^0(\tilde X_k,\cO_{\tilde X_k}(H_{D,k}))$ and take
the products
$$
h_j=g_j\,(\tilde\sigma_P)^{q-1}\,(\tau_{D'})^t~\eta_{D,k}
\in H^0\big(\tilde X_k,
\cO_{\tilde X_k}(mq)\otimes\cO_{\tilde X_k}(\underline a)\otimes
\cO_{\tilde X_k}(-qG_{D,k,m})\big)\leqno(8.18)
$$
where $q=s+b+tc+1$. We now observe, thanks to our choice
$t=|\underline a|=\sum a_\ell$, that
$$
\leqalignno{
&\kern35pt
h_j(\tilde f_{[k]})\cdot (\tilde f'_{k-1})^{mq}\cdot \prod_{1\le\ell\le k}
(d\tilde\pi_{k,\ell}(\tilde f'_{k-1}))^{a_\ell}
&(8.19)\cr
&=\Big(\tilde\sigma_P(\tilde f_{[k]})\cdot (\tilde f'_{k-1})^m\Big)^{q-1}\times
\Big(g_j(\tilde f_{[k]})\cdot (\tilde f'_{k-1})^m\cdot
\prod_{1\le\ell\le k}d\tilde\pi_{k,\ell}(\tau_{D'}(f)\,\tilde f'_{k-1})^{a_\ell}\Big)
\times \eta_{D,k}(\tilde f_{[k]})\cr}
$$
is a product of holomorphic sections on $\bC$, by (8.11) and
(8.9$\tilde{~}$) combined with (8.16) and (8.18), and the fact that
$P(f_{[k]})=\tilde\sigma_P(\tilde f_{[k]})\cdot (\tilde f'_{k-1})^m$ 
is holomorphic with values in~$f^*\cO_X(-A)$. The product also takes
value in the trivial bundle over $\bC$, and can thus be seen as a 
holomorphic function. As $j$ varies, these functions are not all 
equal to zero, and we define a
hermitian metric $\gamma(t)=\gamma_0(t)\,|dt|^2$ on the
complex line $\bC$ by putting
$$
\gamma_0=\Bigg(\sum_je^{\psi(\tilde f_{[k]})}\bigg|
h_j(\tilde f_{[k]})\cdot(\tilde f'_{k-1})^{mq}\cdot \prod_{1\le\ell\le k}
d\tilde\pi_{k,\ell}(\tilde f'_{k-1})^{a_\ell}\bigg|^2~
\Bigg)^{1\over mq+|\underline a|},
\leqno(8.20)
$$
where $\psi$ is a quasi plurisubharmonic potential on $\tilde X_k$ which
will be chosen later. Notice that $\gamma_0(t)$ is locally bounded from above
and almost everywhere non zero.
Since (8.19) only involves holomorphic factors in the right hand side, we get
$$
\ii\ddbar\log\gamma_0\ge
{1\over mq+|\underline a|}\,(\tilde f_{[k]})^*(\tilde\omega_k+\ii\ddbar\psi)
\leqno(8.21)
$$
where $\tilde\omega_k=\ii\ddbar\log|g_j|^2$ is a K\"ahler metric on $\tilde X_k$,
equal to the curvature of the very ample line bundle $L_{k,m}$
for the projective embedding provided by $(g_j)$. (In fact, (8.21) could be
turned into an equality by adding a suitable sum of Dirac masses in
the right hand side). Of course, $\psi$ will be taken
to be an $\omega$-plurisuharmonic potential on $\tilde X_k$.
We wish to get a contradiction by means of
the Ahlfors-Schwarz lemma (see e.g.\ [Dem97, Lemma~3.2]), by showing
that $\ii\ddbar\log\gamma_0\ge A\gamma$ for some $A>0$, an
impossibility for a hermitian metric on the entire complex line.
Since $\psi$ is locally bounded from above,
by (8.19) and the inequality between geometric and arithmetic means,
we have
$$
\gamma_0(t)\le C\,\Big(\sum|h_j(\tilde f_{[k]}(t))|^2\Big)^{1\over mq+|\underline a|}
\,|f'_{[k-1]}(t)|_{\log}^2\leqno(8.22)
$$
where $C>0$ and the norms $|h_j|^2$ and $|f'_{[k-1]}(t)|_{\log}^2$ are computed
with respect to smooth metrics on
$\cO_{\tilde X_k}(mq)\otimes\cO_{\tilde X_k}(\underline a)\otimes
\cO_{\tilde X_k}(-qG_{D,k,m})$ and on the logarithmic
tautological line bundle $\cO_{X^\rS_k(V\langle D'\rangle)}(-1)$, respectively.
The term $|h_j|^2$ is bounded, but one has to pay attention to the fact that
$|f'_{[k-1]}(t)|_{\log}^2$ has poles on $f^{-1}(|D'|)$. If we use local
coordinates $(z_1,\ldots,z_n)$ on $X$ such that $\Delta_j=\{z_j=0\}$, we have
$$
|f'_{[k-1]}|_{\log}^2\sim
|f'_{[k-1]}|_{\omega_k}^2+\sum_j|f_j|^{-2}\,|f'_j|^2
$$
in terms of a smooth K\"ahler metric $\omega_{k-1}$ on
$X^\rS_k(V\langle D'\rangle)$. What saves us is that $h_j$ contains
a factor $\tau_{D'}(f)^t$ that vanishes along all components $\Delta_j$.
Therefore (8.22) implies the existence of a number $\delta>0$ such that
$$
\gamma_0(t)\le C'\,\Big(\,|f'_{[k-1]}(t)|_{\omega_{k-1}}^2
+\sum_j|f_j|^{-2+2\delta}\,|f'_j|^2\Big).
\leqno(8.22')
$$
Since the morphism $\tilde\pi_{k,k-1}$ has a bounded differential
and $f'_{[k-1]}(t)=d\tilde\pi_{k,k-1}(\tilde f'_{[k]}(t))$,
we infer
$$
\gamma_0(t)\le C''\,\Big(\,|\tilde f'_{[k]}(t)|_{ \tilde\omega_k}^2
+\sum_j|f_j|^{-2+2\delta}\,|f'_j|^2\Big).
\leqno(8.22'')
$$
By (8.21) and $(8.22'')$, in order to get a lower bound
$\ii\ddbar\log\gamma_0\ge A\gamma$, we only need to choose
the potential $\psi$ so that
$$
\sum_j|f_j|^{-2+2\delta}\,|f'_j|^2\leq C'''
\,(\tilde f_{[k]})^*(\tilde\omega_k+\ii\ddbar\psi).\leqno(8.23)
$$
If $\tau_j\in H^0(X,\cO_X(\Delta_j))$ is the canonical section of divisor
$\Delta_j$, (8.23) is achieved by taking $\psi=\varepsilon\sum_j
|\tau_j\circ\tilde\pi_{k,0}|^{2\delta}$, for any
choice of a smooth hermitian metric on $\cO_X(\Delta_j)$ and
$\varepsilon>0$ small enough. In some sense, we have to take a suitable
orbifold K\"ahler metric $\tilde\omega_k+\ii\ddbar\psi$ on $\tilde X_k$
to be able to apply the Ahlfors-Schwarz lemma. It might be interesting
to find the optimal choice of $\delta>0$, but this is not
needed in our proof.\qed

\noindent
{\it End of the proof of Proposition~0.7.} We still have to extend the
vanishing result to the case of non necessarily $\bG_k$-invariant
orbifold jet differentials
$$
P\in H^0(X,E_{k,m}V^*\langle D\rangle\otimes\cO_X(-A)).
$$
One can then argue by using the $\bG_k$-action on jet differentials
$$
(\varphi,P)\mapsto\varphi^*P,\qquad
(\varphi^*P)(f_{[k]}):=P((f\circ\varphi)_{[k]})\circ\varphi^{-1},\quad
\varphi\in\bG_k.
\leqno(8.24)
$$
This action yields a decomposition
$$
(\varphi^*P)(f_{[k]})=\sum_{
\scriptstyle\alpha\in\bN^k\atop\scriptstyle|\alpha|_w=m
\phantom{\textstyle|}}(\varphi^{(\alpha)}
\circ\varphi^{-1})\,P_\alpha(f_{[k]}),\quad
P_\alpha\in H^0(X,E_{k,m_\alpha}V^*\langle D\rangle\otimes\cO_X(-A))
\leqno(8.25)
$$
where $\alpha=(\alpha_1,\ldots,\alpha_k)\in\bN^k$,
$\varphi^{(\alpha)}=(\varphi')^{\alpha_1}(\varphi'')^{\alpha_2}\ldots
(\varphi^{(k)})^{\alpha_k}$,
$|\alpha|_w=\alpha_1+2\alpha_2+\ldots+k\alpha_k$ is the weighted degree,
and $P_\alpha$ is a homogeneous polynomial of degree
$$
m_\alpha:=\deg P_\alpha=m-(\alpha_2+2\alpha_3+\ldots+(k-1)\alpha_k)=
\alpha_1+\alpha_2+\ldots+\alpha_k.
$$
In particular $\deg P_\alpha<m$ unless $\alpha=(m,0,\ldots,0)$, in which case
$P_\alpha=P$. If the result is known for degrees${}<m$, then all
$P_\alpha(f_{[k]})$ vanish for $P_\alpha\ne P$ and one can reduce the proof
to the invariant case by induction, as the term $P_\alpha$ of minimal degree
is invariant. The proof makes use of induced directed
structures, and is purely formal and group theoretic. Essentially,
the argument is that $P$ becomes an invariant jet differential when
restricted to the subvariety of the Semple $k$-jet bundle consisting
of germs $g_{[k]}$ of $k$-jets such that $P_\alpha(g_{[k]})=0$ for 
$P_\alpha\ne P$.
Singularities may appear in this subvariety, but this does not affect the proof
since the induced directed structure is embedded in the non singular
logarithmic Semple tower.  We refer the reader to [Dem20,~\S$\,$7.E]
and [Dem20,~Theorem~8.15]
for details.\qed
\bigskip

\centerline{\twelvebf References}
\medskip

\Bibitem[BD18]&Brotbek, D., Deng, Y.&: Kobayashi hyperbolicity of the 
complements of general hypersurfaces of high degree.& Geom.\ Funct.\
Anal.\ {\bf 29 }(2019), no. 3, 690--750&

\Bibitem[Cad17]&Cadorel, B.:& Jet differentials on toroidal compactifications
of ball quotients.& arXiv: math.AG/1707.07875&

\Bibitem[Cad19]&Cadorel, B.:& Generalized algebraic Morse inequalities
and jet differentials.& arXiv: math.AG/1912.03952&
  
\Bibitem[CDR20]&Campana, F., Darondeau, L., Rousseau, E.:& Orbifold 
hyperbolicity.& Compos.\ Math.\ {\bf 156} (2020), 1664--1698&

\Bibitem[Dar16]&Darondeau, L.:& On the logarithmic Green-Griffiths conjecture.&
Int.\ Math.\ Res.\ Not.\ IMRN {\bf 6} (2016), 1871--1923&

\Bibitem[DR20]&Darondeau, L., Rousseau, E.:& Quasi-positive orbifold cotangent
bundles; pushing further an example by Junjiro Noguchi.& arXiv:2006.13515&
  
\Bibitem[Dem80]&Demailly, J.-P.:& Relations entre les diff\'erentes notions
de fibr\'es et de courants positifs.& S\'em.\ P.~Lelong-H.~Skoda (Analyse)
1980/81, Lecture Notes in Math.\ n${}^\circ\,$919, Springer-Verlag, 56--76&

\Bibitem[Dem95]&Demailly, J.-P.:& Propri\'et\'es de semi-continuit\'e de la cohomologie
et de la dimension de Kodaira-Iitaka.& C.~R.\ Acad.\ Sci.\ Paris S\'er.~I Math.\
{\bf 320} (1995), 341--346&

\Bibitem[Dem97]&Demailly, J.-P.:& Algebraic criteria for Kobayashi
hyperbolic projective varieties and jet differentials.& AMS Summer
School on Algebraic Geometry, Santa Cruz 1995, Proc.\ Symposia in
Pure Math., ed.\ by J.~Koll\'ar and R.~Lazarsfeld, Amer.\ Math.\ Soc.,
Providence, RI (1997), 285–-360&

\Bibitem[Dem11]&Demailly, J.-P.:& 
Holomorphic Morse Inequalities and the Green-Griffiths-Lang Conjecture.&
Pure and Applied Math.\ Quarterly {\bf 7} (2011), 1165--1208&

\Bibitem[Dem12]&Demailly, J.-P.:& Hyperbolic algebraic varieties and
holomorphic differential equations.& Expanded version of the lectures
given at the annual meeting of VIASM, Acta Math.\ Vietnam.\ {\bf 37}
(2012), 441-–512&

\Bibitem[Dem20]&Demailly, J.-P.:& Recent results on the Kobayashi and
Green-Griffiths-Lang conjectures.& Expanded version of talks given at
the 16th Takagi Lectures in Tokyo, November 29, 2015, Japanese Journal
of Mathematics volume {\bf 15} (2020),  1--120&
  
\Bibitem[GrGr80]&Green, M., Griffiths, P.:& Two applications of algebraic
geometry to entire holomorphic mappings.& The Chern Symposium 1979,
Proc.\ Internal.\ Sympos.\ Berkeley, CA, 1979, Springer-Verlag, New York
(1980), 41--74&

\Bibitem[Lan05]&Landau, E.:& Sur quelques th\'eor\`emes de M.~Petrovitch
relatifs aux z\'eros des fonctions analytiques.& Bull.\ Soc.\ Math.\ France
{\bf 33} (1905), 251--261&

\Bibitem[MTa19]&Merker, J., Ta, The-Anh :&
Degrees $d\ge (\sqrt{n}\,\log n)^n$ and $d\ge (n\,\log n)^n$ in the Conjectures
of Green-Griffiths and of Kobayashi.& arXiv:1901.04042, [math.AG]&

\Bibitem[Sem54]&Semple, J.G.:& Some investigations in the geometry of
curves and surface elements.& Proc.\ London Math.\ Soc.\ (3) {\bf 4}
(1954), 24--49&
\vskip1cm

\parindent=0cm
(Version of November 2, 2021, printed on \today, \timeofday)
\bigskip

Fr\'ed\'eric Campana\\
Institut de Math\'ematiques \'Elie Cartan,
Universit\'e de Lorraine, B.P.\ 70239\\
54506 Vand{\oe}uvre-l\`es-Nancy, France\\
E-mail : frederic.campana@univ-lorraine.fr
\medskip

Lionel Darondeau\\
Universit\'e Montpellier II, 
Institut Montpelli\'erain Alexander Grothendieck,\\
Case courrier 051, Place Eug\`ene Bataillon, 34090 Montpellier, France\\
E-mail : lionel.darondeau@normalesup.org
\medskip

Jean-Pierre Demailly\\
Universit\'e Grenoble Alpes, Institut Fourier\\
100 rue des Maths, 38610 Gi\`eres, France\\
E-mail : jean-pierre.demailly@univ-grenoble-alpes.fr
\medskip

Erwan Rousseau\\
Institut Universitaire de France, Universit\'e de Bretagne Occidentale\\
6, avenue Victor Le Gorgeu, 29238 Brest Cedex 3, France\\
E-mail : erwan.rousseau@univ-brest.fr

\end{document}